\documentclass[reqno]{amsart}
\usepackage{amssymb,amscd,verbatim, amsthm,graphicx, color}

\usepackage[mathscr]{eucal}

\usepackage{longtable}

\newcommand \fk[1]{{{\mathfrak #1}}}
\newcommand \C[1]{{\mathcal #1}}
\newcommand \ovl[1]{{\overline {#1}}}
\newcommand \ch[1]{{\check{#1}}}

\newcommand \ti[1]{{\tilde #1}}
\newcommand \wti[1]{{\widetilde {#1}}}
\newcommand \wht[1]{{\widehat {#1}}}

\newcommand\fg{\mathfrak g}
\newcommand\fn{\mathfrak n}

\newcommand\cCO{{\check{\CO} }}
\newcommand \vO{{\check \CO}}

\newcommand \bA{{\mathbb A}}
\newcommand \bC{{\mathbb C}}
\newcommand \bF{{\mathbb F}}
\newcommand \bH{{\mathbb H}}
\newcommand \bR{{\mathbb R}}
\newcommand \bZ{{\mathbb Z}}

\newcommand\cha{{\check \alpha}}
\newcommand\oal{{\ovl \alpha}}

\newcommand\one{1\!\!1}

\newcommand\CA{{\C A}}
\newcommand\CB{{\C B}}
\newcommand\CH{{\C H}}

\newcommand\CO{{\C O}}

\newcommand\CX{{\C X}}

\newcommand\CR{{\C R}}

\newcommand\ic{infinitesimal character }

\newcommand\ie{{\it i.e.~}}

\newcommand\eg{{\it e.g.~ }}

\newcommand\ep{{\epsilon}}

\newcommand\om{{\omega}}
\newcommand\al{{\alpha}}
\newcommand\sig{{\sigma}}

\newcommand \vg{\check{\fk g}}
\newcommand \vm{\check{\fk m}}

\newcommand \ve{{\check e}}

\newcommand \vG{{\check{G}}}

\newcommand\oX{\ovl{X}}

\newcommand\fa{{\mathfrak a}}
\newcommand\fh{{\mathfrak h}}
\newcommand\fm{{\mathfrak m}}
\newcommand\ft{{\mathfrak t}}
\newcommand\fz{{\mathfrak z}}

\newtheorem*{theorem}{Theorem}
\newtheorem*{theorem 2}{Theorem 2}
\newtheorem*{corollary}{Corollary}
\newtheorem*{lemma}{Lemma}
\newtheorem*{proposition}{Proposition}
\newtheorem*{definition}{Definition}
\newtheorem*{remark}{Remark}

\newcommand\Hom{\operatorname{Hom}}
\newcommand\Ind{\operatorname{Ind}}

\numberwithin{equation}{subsection}

\begin{document}

\title[Whittaker unitary dual]{Whittaker unitary dual of affine graded
Hecke algebras of type $E$}

\author{Dan Barbasch}
       \address[D. Barbasch]{Dept. of Mathematics\\
               Cornell University\\Ithaca, NY 14850}
       \email{barbasch@math.cornell.edu}

\author{Dan Ciubotaru}
        \address[D. Ciubotaru]{Dept. of Mathematics\\ University of
          Utah\\ SLC, UT 84112}
        \email{ciubo@math.utah.edu}

\date\today

\begin{abstract}
This paper gives the classification of the Whittaker unitary dual for
affine graded Hecke algebras of type $E$. By the Iwahori-Matsumoto
involution, this is equivalent also to the classification of the
spherical unitary dual for type $E$. Together with \cite{BM3}, \cite{Ba2}, and
\cite{C}, this work
completes the classification of the Whittaker Iwahori-spherical
unitary dual, or equivalently, the spherical unitary dual of any
split $p$-adic group. 
\end{abstract}

\maketitle

\setcounter{tocdepth}{1}


\begin{small}
\tableofcontents 
\end{small}



\section{Introduction}\label{sec:1}

\subsection{}\label{1.1}

The present paper completes the classification of the unitary
representations with Iwahori-fixed vectors and {\it generic} (\ie
admitting Whittaker models)  
for split linear algebraic groups over $p$-adic fields by treating the groups of type
$E$. 

The full unitary dual for $GL(n)$ was obtained
in \cite{T}, and for $G_2$ in \cite{M}. The Whittaker unitary dual
with Iwahori-fixed vectors for classical split groups was determined
in \cite{BM3} and \cite{Ba2}. For $F_4$, this is part of  \cite{C}.

It is well-known that the category of representations with Iwahori
fixed vectors admits an involution called the Iwahori-Matsumoto
involution, denoted by $IM$, which takes hermitian modules to hermitian modules,
and unitary modules to unitary modules (\cite{BM1}). In particular it
interchanges spherical modules with generic modules. For example, $IM$
takes the trivial representation into the Steinberg representation.  Thus 
this paper also gives a classification of the spherical unitary dual
of split p-adic groups of type $E,$  completing the classification
of the spherical unitary dual as well.
 
\subsection{}
Let $\mathcal G$ be a split reductive linear algebraic group over a
$p$-adic field $\bF$ of 
characteristic zero.   Recall that $\bF\supset\C R\supset\C P,$ where
$\C R$ is the ring of integers and  $\C P$ the maximal prime ideal. We
fix a maximal compact subgroup $\mathcal K=\mathcal G(\C R).$ Let $\C
I$ be an Iwahori subgroup, $\C I\subset \C K.$ 
Fix also a  rational Borel subgroup $\mathcal B=\mathcal H\mathcal N .$
Then $\mathcal G=\mathcal K\mathcal B.$
An admissible representation $(\pi,V)$  is called  {\it spherical}
if $V^{\mathcal K}\ne (0).$  It is called {\it Iwahori-spherical} if
$V^{\mathcal I}\ne (0).$ 

We recall the well known classification of irreducible admissible 
spherical modules. For every irreducible spherical representation $\pi,$
there is a character $\chi\in\widehat{\mathcal H}$ such that
$\chi |_{\mathcal H\cap\mathcal K}=triv,$ and $\pi$ is the unique
spherical subquotient $L(\chi)$  
of $X(\chi)=\Ind_{\mathcal B}^{\mathcal G}[\chi\otimes\one].$ A
character $\chi$ whose restriction to  $\mathcal H\cap\mathcal K$ is
trivial is called \textit{unramified}. Two modules $L(\chi)$ and
$L(\chi')$ are equivalent if and only if there is an element in the
Weyl group $W$ such that $w\chi=\chi'.$ A module $L(\chi)$ admits an
invariant hermitian form if
and only if there exists $w$ such that {$w\chi=\ovl{\chi}^{-1}.$} 

More generally, by a theorem of Casselman, every irreducible
Iwahori-spherical representation of 
$\C G$ is a  subquotient of an $X(\chi)$. Furthermore each $X(\chi)$
has a unique irreducible subquotient which is {\it generic}. 

When $\chi$ is dominant, the spherical module $L(\chi)$ is the unique
irreducible quotient of $X(\chi).$ In this case, it is known that the
generic subquotient is  a submodule of $X(\chi)$. By \cite{BM4}, if
$\C G$ is adjoint, a
subquotient is both generic and spherical if and only if it is the full
$X(\chi)$ (in other words, if $X(\chi)$ is irreducible).

\medskip

The results in \cite{BM1} show that
in the $p-$adic case the classification of the Iwahori-spherical
unitary dual is 
equivalent to the corresponding problem for the Iwahori-Hecke
algebra. In \cite{BM2}, the problem is reduced to computing the
unitary dual of the case of the affine graded (Iwahori-)Hecke algebra
of a possibly smaller group, and 
\textit{real infinitesimal character}. We will review these notions
later in the paper, for now we recall the notion of \textit{real infinitesimal
character}.   
A character $\chi$ is called \textit{real} if it takes only positive
real values.  An irreducible representation $\pi$ is said to have real
infinitesimal character if it is the subquotient of an $X(\chi)$ with
$\chi$ real. So we start by parameterizing real unramified characters
of $\mathcal H.$  Since $\mathcal G$ is split, $\mathcal H\cong
(\bF^\times)^n$ where $n$ is the rank.  Define
\begin{equation}
  \label{eq:1.1.1}
  \mathcal L(\mathcal H)_\bR^*:= X^*(\mathcal H)\otimes_{\bZ}\bR,
\end{equation}
where $X^*(\mathcal H)$ is the lattice of characters of the algebraic
torus $\mathcal H.$ 
Each element $\nu\in\mathcal L(\mathcal H)_\bR^*$ defines an
unramified character 
$\chi_\nu$ of $\mathcal H$, determined by the formula
\begin{equation}
  \label{eq:1.1.2}
  \chi_\nu(\tau(a))=|a|^{\langle \tau,\nu\rangle},\qquad a\in\bF^\times,
\end{equation}
where  $\tau$ is an element of the lattice of one
parameter subgroups $X_*(\mathcal H).$ Since the torus is split, each element of
$X^*(\mathcal H)$ can be regarded as a homomorphism of $\bF^\times$
into $\mathcal H.$
The pairing in the exponent in (\ref{eq:1.1.2}) corresponds to the
natural identification of $\mathcal L(\mathcal H)_\bR^*$ with
$\Hom[X_*(\mathcal H),\bR].$
The map $\nu\longrightarrow \chi_\nu$ from $\mathcal L(\mathcal H)_\bR^*$ to real
unramified characters of $\CH$ is an isomorphism. We will often
identify the two sets writing simply $\chi\in\mathcal L(\mathcal
H)_\bR^*.$ 

\medskip
Because we will be dealing exclusively with the graded affine Hecke
algebra $\bH$ (introduced in \cite{L2}) which is defined in terms of the
complex dual group, we will denote by $G$ the complex group dual to $\mathcal G$, and let
$H$ be the torus dual to $\mathcal H.$ Then  the real unramified characters $\chi$ are
naturally identified with hyperbolic
elements of the Lie algebra $\fh.$ The infinitesimal characters
are identified with orbits of hyperbolic  elements (section \ref{sec:2.1}). 
We will  assume that all characters are real. 

\subsection{} Next we explain the nature of our classification of the
Whittaker unitary dual.

We attach to each $\chi$ a nilpotent orbit $\CO=\CO(\chi)$ satisfying
the following properties.  Fix a
Lie triple $\{e,h,f\}$ corresponding to $\CO$ such that
$h\in\fh.$ We write $Z(e,h,f)$, respectively $\fk
z(e,h,f)$ for the  centralizer of $\{e,h,f\}$ 
in $G$, respectively $\fg,$ and  abbreviate them $Z(\CO)$,
respectively 
$\fk z(\CO).$ Then $\CO$ is such that 
\begin{align}\label{eq:1.1.3}
&\text{(1) there exists $w\in W$ such that $w\chi=\frac12 h +\nu$ with
  $\nu\in\fk z(\CO),$}\\\notag
&\text{(2) if $\chi$ satisfies property (1) for any other $\CO',$ then
  $\CO'\subset\ovl{\CO}.$}
\end{align}
The results in \cite{BM1} guarantee that for any $\chi$ there is a
unique $\CO(\chi)$ satisfying (1) and (2) above.
Another characterization of the orbit $\CO=\CO(\chi)$ is as follows.
Set
\begin{equation}
\fg_1:=\{\ x\in \fg\ :\ [\chi,x]=x \ \},\qquad 
\fg_0:=\{\ x\in \fg\ :\ [\chi,x]=0 \ \}.
\end{equation}
Then $G_0$, the complex Lie group corresponding to the Lie algebra $\fg_0$
has an open dense orbit in $\fg_1.$ The $G-$saturation in $\fg$
of this orbit is $\CO.$  

\smallskip

Every generic module of $\bH$ (and every spherical module
of $\bH$) is uniquely parametrized by a pair
$(\CO,\nu)$, $\CO=\CO(\chi)$ as in (\ref{eq:1.1.3}). In order to make
this connection more precise, we will need to recall the geometric
classification of irreducible $\bH$-modules (\cite{KL},\cite{L5}), and we postpone
this to section \ref{sec:2.2}. (See in particular remark \ref{r:2.2}.)  

\smallskip

\noindent{\bf Remark.} The pair $(\CO,\nu)$ has remarkable
properties. For example, if  $\nu=0$ ($\chi=h/2$), then the generic
representation parametrized by $(\CO,0)$  is tempered, therefore
unitary. The corresponding spherical module $L(h/2)$ is one of the
parameters that the conjectures of Arthur predict to play a role in
the residual spectrum. In particular, $L(h/2)$ should be unitary. This is true
because it is the Iwahori-Matsumoto dual (definition \ref{sec:2.2.b})
of the generic tempered module.

\begin{definition}[1]
  \label{d:1.3}
The spherical modules $L(h/2)$ will be called spherical unipotent
representations.
\end{definition}
In our main result, theorem \ref{main}, the tempered generic modules
can be regarded as the building blocks of the Whittaker unitary
dual. In the spherical unitary dual, this role is played by
the spherical unipotent modules.

\medskip

We partition the Whittaker (equivalently, the spherical)
unitary dual into {\it complementary series} attached to nilpotent
orbits. We say that an infinitesimal character $\chi$ as above is {\it
  unitary} if the generic module parametrized by $\chi$ (equivalently,
the spherical $L(\chi)$) is unitary.
 
\begin{definition}[2]\label{d:1.2}
The (generic or spherical) {\it $\CO$-complementary series} is the set
of unitary parameters $\chi$ such that $\CO=\CO(\chi)$
as in (\ref{eq:1.1.3}). The complementary series for the trivial
nilpotent orbit is called the {\it 0-complementary series}. 
\end{definition}

Our first result is the identification of $0$-complementary series for
type $E$ in section \ref{zero}. These are the irreducible principal
series $X(\chi)$ which are unitary. 
(For a summary of the relevant results for
classical groups from \cite{Ba2}, and $G_2,\ F_4$, from \cite{C},  see sections
\ref{classical},\ref{p:g2}.) The $0$-complementary series have a nice
explicit combinatorial description: they can be viewed as a union of alcoves in
the dominant Weyl chamber of $\fh$, where the number of alcoves is a
power of $2$, \eg in $G_2$ there are $2$, and in $E_7,$ $E_8$, there
are $8$, respectively $16.$
The explicit description of the alcoves is in sections
\ref{sec:genspherE6}-\ref{sec:genspherE8}.

\smallskip

The main result of the paper is the description of the complementary
series for all $\CO$ in type $E$, and can be summarized as follows. 
We use the Bala-Carter notation for nilpotent orbits in exceptional
complex semisimple Lie algebras (see \cite{Ca}). 
\begin{definition}[3]
  \label{d:1.1.5}
Set
\begin{equation}\label{eq:1.1.5}
Exc=
\{\underbrace{A_1\wti A_1}_{\text{in }F_4},\underbrace{A_2+3A_1}_{\text{in }E_7},
\underbrace{A_4A_2A_1,A_4A_2,D_4(a_1)A_2,A_3+2A_1,A_2+2A_1,4A_1}_{\text{in }E_8}\}. 
\end{equation}

\end{definition}

Recall that $\fz(\CO)$ denotes 
the reductive algebra which is the centralizer in $\fg$ of a fixed Lie
triple for $\CO.$ 

\begin{theorem}\label{main} Let $\bH$ be the affine graded Hecke
  algebra for $G$ (definition in section \ref{sec:2.1}), and $\CO$ be
  a nilpotent $G$-orbit in the complex Lie 
  algebra $\fg.$ Denote by $\bH(\fz(\CO))$  the affine graded Hecke algebra constructed from the root system of $\fz(\CO).$ 

Assume $\CO\notin Exc$ (definition (\ref{eq:1.1.5})).
A (real) parameter $\chi=\frac 12 h+\nu$ is in the
complementary series of $\CO$ (definition \ref{d:1.2}) if and only if
$\nu$ is in the $0$-complementary series of $\bH(\fz(\CO)).$ 

The explicit description of the complementary series, including when
$\CO\in Exc,$ are tabulated in section \ref{tables}. 
\end{theorem}

The complementary series for $\CO\in Exc$ are smaller than the
corresponding $0$-complementary series for $\bH(\fz(\CO))$, except
when $\CO=4A_1$ in $E_8.$ For this one orbit, the complementary
series is larger (see section \ref{sec:4A1}). 

The proof of the theorem for $G$ of classical type is in
\cite{Ba2}. For types $G_2$ and $F_4$, this is part of 
\cite{C}. In the present paper, we treat the case of groups of type
$E.$ The method is different from the above mentioned papers.

\medskip

The main method of the proof (proposition \ref{p:4.2})  consists of a
{\it direct} 
comparison between the signature of hermitian forms on the generic modules
for $\bH$ parametrized in the geometric classification (see
section \ref{sec:2.2}) by
$\CO$, and the signature of hermitian forms on the
spherical principal series which are irreducible (that is
representations which are both spherical and generic) for the Hecke
algebra $\bH(\fz(\CO)).$

This method of comparing signatures has the advantage that it explains
the match-up  of complementary series in theorem \ref{main}. It often
extends to non-generic modules (\eg \cite{Ci2} for non-generic modules
of $E_6$).   

\subsection{} If one assumes the infinitesimal character (the $\chi$
above) to be real, {one can use the same set for the
  parameter spaces for the spherical dual of a real and 
$p$-adic split group (attached to the same root datum).} The
main criterion for ruling out nonunitary modules is the computation of
signature characters: in the real case on $K-$types, and in the
$p$-adic case on $W-$types. So it is natural to try to compare
signatures on $K-$types and $W-$types. In \cite{Ba2} and \cite{Ba1}, the notion
of {\it petite $K-$types} was used to transfer the results about
signatures from the $p-$adic split group to the corresponding real
split group. The methods employed there are very different from this
paper. More precisely, to every petite $K-$type there corresponds a
Weyl group representation such that the signature characters are the
same. But only a small number of Weyl group representations correspond in
this way to petite $K-$types.
In this paper, we inherently use signature computations for {\it all}
Weyl group representations. Therefore, the results here cannot be used
directly towards the spherical unitary dual of the corresponding split real
groups. 
 
In different work however, we studied the signature of petite $K-$types for
exceptional groups of type $E$. The main consequence of that work is that the
spherical unitary dual for a split real group $\mathcal G(\mathbb R)$ is
a {\it subset} 
of the spherical unitary dual for the corresponding $p$-adic group
$\mathcal G(\mathbb F)$ (conjecturally they are the same). Details will
appear elsewhere.

\subsection{}\label{1.1a}
To obtain the results of this paper, we made a
minimal use of computer calculations, essentially for linear
algebra, \eg conjugation of semisimple elements by the Weyl group, or
multiplication of matrices in a variable $\nu$ for some of the
``maximal parabolic'' cases in section \ref{sec:3.1}.

However, by the machinery presented in sections \ref{sec:2.11} and
\ref{sec:alcoves}, for every given Hecke algebra $\bH$, one can reduce the
identification of the unitary parameters $\chi$ to a brute force computer
calculation. More precisely, one considers sample points with rational
coordinates for every facet in the arrangement of hyperplanes given by
coroots equal to $1$ in the dominant Weyl chamber. (These are the
hyperplanes where $X(\chi)$ is reducible.) It is known (see
\cite{BC}) that the signature of the long intertwining operator is
constant on each facet of this arrangement.

 One can then run a
computer calculation of the long intertwining spherical operator
(section \ref{sec:2.11}) on {\it each} representation of the Weyl 
group at {\it every} sample point. Then one finds the signature of the
resulting hermitian matrices. The unitary parameters $\chi$ correspond to those
facets for which these matrices are positive semidefinite for all Weyl
group representations. The size of the
calculation can be reduced significantly by making use of
some ideas in this paper.   This is {\it
  not} the approach of this paper, but we did carry out this
calculation independently for exceptional groups in order to confirm
the results of this paper.

\subsection{}\label{sec:1.2}
We give an outline of the paper. In section \ref{sec:2}, we review the
relevant notions
about the affine graded Hecke algebra, and its representations. We
introduce the construction of intertwining operators that we need for
the signature computations.
In section \ref{zero}, we restrict to the setting of modules which
are both generic and spherical, and determine the $0$-complementary
series. 
In section \ref{sec:extended}, we describe a construction of extended
Hecke algebras for disconnected groups, and apply it to the setting of
centralizers of nilpotent orbits.
Section \ref{sec:4} contains the main results of the paper, theorem
\ref{t:4.1}, propositions \ref{p:3.1} and \ref{p:4.2}, and presents the
main ingredients of the method for matching signatures of intertwining
operators. 
The explicit details and calculations needed for the proofs are in
section \ref{sec:5}. For the convenience of the
reader, the results, including the 
exact description of the complementary series for $\CO\in Exc$, and
of the $0$-complementary series
are tabulated in section \ref{tables}.

\smallskip
This research was supported by NSF grants DMS-9706758, 0070561,
03001712, and FRG-0554278.  

\smallskip

\noindent{\bf Notation.} If $\fk G$ is an algebraic group, we will
denote by $\fk G^0$ its identity component. The center will be denoted
$Z(\fk G).$ For every set of elements $\C E,$ we will denote by
$Z_{\fk G}(\C E)$ the simultaneous centralizer in $\fk G$ of all
elements of $\C E,$ and by $A_{\fk G}(\C E)$ the group of
components of $Z_{\fk G}(\C E).$  

\section{Intertwining operators}\label{sec:2}

\subsection{}\label{sec:2.1}

As mentioned in the introduction, we will work only
with the Hecke algebras and the $p$-adic group will not
play a role. Therefore, in order to simplify notation, we will call
the dual complex group $G$, its Lie algebra $\fg$ etc.

 Let $H$ be a
maximal torus $G$ and $B\supset H$ be a Borel subgroup.
The affine Hecke algebra $\mathscr H$ can be described by generators and relations. Let
$z$ be an indeterminate (which can then be specialized to $q^{1/2}$).  
Let $\Pi\subset \Delta^+\subset \Delta$
be the simple roots, positive roots, respectively roots corresponding to
$H\subset B,$ and $S$ be the simple root reflections. Let
$G_m:=GL(1,\bF)$, $\ch{\C X}:=\Hom(G_m,H)$ be the (algebraic) lattice
of 1-parameter subgroups,
and $\C X:=\Hom(H,G_m)$ the lattice of algebraic characters.
Then $\mathscr H$ will denote the Hecke algebra over $\bC[z,z^{-1}]$ 
attached to the root datum $(\C X,\check{\C X},\Delta,\ch \Delta,\Pi).$ The set
of generators we will use is the one first introduced by
Bernstein.  

The algebra $\mathscr H$ is generated over $\bC[z,z^{-1}]$ 
by $\{ T_w\}_{w\in W}$ and $\{\theta_x\}_{x\in
\CX}$, subject to the relations
\begin{equation}\begin{aligned}
&T_wT_{w'}=T_{ww'}\ (l(w)+l(w')=l(ww')),\quad
\theta_x\theta_y=\theta_{x+y},\\
&T_s^2=(z^2-1)T_s+z^2,\quad \theta_xT_s=T_s\theta_{sx}+(z^2-1)\frac{\theta_x-\theta_{sx}}{1-
\theta_{\alpha}}. 
\end{aligned} \label{2.1.2}\end{equation}
This realization is very convenient for determining the center of
$\mathscr H$ and
thus computing infinitesimal characters of representations. Let
$\mathscr A$
be the subalgebra over $\bC[z,z^{-1}]$ generated by the $\theta_x.$
The Weyl group acts on this via $w\cdot\theta_x:=\theta_{wx}.$  
\begin{proposition} [Bernstein-Lusztig]\label{p:center} The center of
  $\mathscr H$ is
given by $\mathscr A^W,$ the Weyl group invariants in $\mathscr A.$ 
\end{proposition}
Infinitesimal characters are parametrized by  $W-$orbits
$\chi=(q,t)\in \bC^*\times H.$ We always assume that $q$ is real
or at least not a root of unity. In particular, such an infinitesimal character
is called {\it real} if $t$ is hyperbolic. 
Unless indicated otherwise, we will assume from here on that \textbf{the infinitesimal character is
always real}. The study of representations of $\mathscr H$ can be simplified
by using the graded Hecke algebra $\bH$ introduced in \cite{L2}. One can
identify $\mathscr A$ with the algebra of regular functions on $\mathbb
C^*\times H.$ Define
\begin{equation}
\mathscr J=\{f\in \mathscr A\ :\ f(1,1)=0\}.
\label{eq:2.1.3}\end{equation}
This is an ideal in $\mathscr A$ and it satisfies $\mathscr H \mathscr
J=\mathscr J \mathscr H.$ Set
$\mathscr H^i=\mathscr H\cdot \mathscr J^i$ (the ideal $\mathscr J^i$
in $\mathscr H$ consists of the functions
which vanish to order at least $i$ at $(1,1)$). Thus $\mathscr H$ has a 
filtration  
\begin{equation}
\mathscr H=\mathscr H^0\supset \dots \supset \mathscr H^i\supset
\mathscr H^{i+1}\supset \dots,
\label{eq:2.1.4}\end{equation}
and denote the graded object by $\bH.$ It can be written as
\begin{equation}
\bH=\bC[W]\otimes\bC[r]\otimes \bA,
\label{eq:2.1.5}\end{equation}
where $r\equiv z-1\ (mod\ \mathscr J),$ and  $\bA$ is the symmetric algebra
of $\fk h^*=\CX\otimes_{\mathbb Z} \bC.$ The previous relations become
\begin{equation}
\aligned
t_wt_{w'}=&t_{ww'}, \ w,w'\in W,\\
t_s^2=&1,\ s\in S,\\
\omega t_s=&t_ss(\omega) + 2r\langle\omega,\cha\rangle,\quad
s=s_\alpha\in S,\ \omega\in\fk h^*.
\endaligned \label{eq:2.1.6}
\end{equation}
The center of $\bH$ is $\bC[r]\otimes\bA^W$ (\cite{L2}). Thus infinitesimal
characters are parametrized by $W-$orbits of elements 
$\ovl \chi=(r,t)\in \bC\times \fk h.$ 
\begin{theorem}\label{t:2.1} (\cite{L2}) There is a matching $\chi
\longleftrightarrow \ovl{\chi}$ between real infinitesimal characters
$\chi$ of $\mathscr H$ and infinitesimal characters $\ovl \chi$ of $\bH$ so
that if $\mathscr H_\chi$ and $\bH_{\ovl{\chi}}$ are the quotients by the
corresponding ideals, then 
$$
\mathscr H_\chi\cong \bH_{\ovl{\chi}}.
$$
\end{theorem}
We fix  $r=1/2$, and  transfer the study of
the representation theory of $\mathscr H$ to $\bH.$

In order to consider unitary representations for $\bH,$ we also need a
$*$ operation. This is given in section 5 of \cite{BM2}:
\begin{equation}
  \label{eq:2.1.7}
  \begin{aligned}
    &t_w^*=t_{w^{-1}},\quad w\in W,\\
    &\om^*=-\overline\om +\sum_{\al\in \Delta^+}\langle\om,\cha\rangle t_{s_\al},\quad \om\in\fh^*.
  \end{aligned}
\end{equation}


\subsection{}\label{sec:2.2.b} 
\begin{definition}[1] A $\mathbb H$-module $V$ is called {\it
  spherical} if $V|_W$ contains the trivial $W-$type. 
The module $V$ is called {\it generic} if $V|_W$ contains the sign
$W-$type (\cite{BM4}).
\end{definition}

\begin{definition}[2]
The {\it Iwahori-Matsumoto involution} $IM$ is defined as
\begin{equation}
  \label{eq:2.11.1}
  \begin{aligned}
    &IM(t_w):=(-1)^{l(w)}t_w,\\
    &IM(\om):=-\om,\quad \om\in\fk h^*.
  \end{aligned}
\end{equation}
\end{definition}

$IM$ takes spherical modules into generic modules and it preserves
unitarity. In particular, $IM(triv)$ is the Steinberg module $St.$

\subsection{}\label{sec:2.2}

We parameterize irreducible representations of $\bH$ as in \cite{L4}
and \cite{L5} by $G$-conjugacy classes $(\chi,e,\psi),$ where
$\chi\in \fg$ is semisimple, $e\in \fg$ is nilpotent such that 
$[\chi,e]=e$, and $(\psi,V_\psi)$ are certain irreducible representations of
${A(e,\chi)},$ the component group of the centralizer in $G$ of $e$ and $\chi.$

Embed $e$ into a Lie triple $\{e,h,f\}$. Write $\chi=h/2+\nu$ where
$\nu$  centralizes $\{e,h,f\}$.

\medskip
The results in \cite{KL} and \cite{L5} attach to each ($G$-conjugacy class) 
$(e,\chi)$ a module
$X(e,\chi)$ which decomposes under the action of $A(e,\chi)$ as
a sum of standard modules $X(e,\chi,\psi)$: 
\begin{equation}
  \label{eq:2.2.1}
X(e,\chi)=
\bigoplus_{(\psi,V_\psi)\in\widehat{A(e,\chi)}_0}X(e,\chi,\psi)\otimes V_\psi,  
\end{equation}
where $\widehat{A(e,\chi)}_0$ will be defined below.
  
As a  $\bC[W]$-module, 
\begin{equation}
  \label{eq:2.2.2}
X(e,\chi)\cong H^*(\CB_e), 
\end{equation}
where $\CB_e$ is the variety of Borel subalgebras of $\fg$ containing
$e$. The action of $W$ is the generalization of the one defined  by Springer.
The component group $A(e,\chi)$ is naturally a subgroup of $A(e)$ because
in a connected algebraic group, the centralizer of a torus is
connected. The group $A(e)$ acts on the right hand side of
(\ref{eq:2.2.2}), and the action of $A(e,\chi)$ on $X(e,\chi)$ is
compatible with its inclusion into $A(e),$ and the isomorphism in
(\ref{eq:2.2.2}). Let $\CO$ be the
$G$ orbit of $e.$ According to the {\it Springer correspondence} (\cite{Sp}), 
\begin{equation}
  \label{eq:2.2.3}
  H^*(\C B_e)=\bigoplus_{\phi\in\widehat{A(e)}} H^*(\C
  B_e)^\phi\otimes V_\phi.
\end{equation}
 Furthermore,  $H^{2\dim(\CB_e)}(\CB_e)^\phi$ 
is either zero, or it is an irreducible representation of $W$. Denote
\begin{equation}\label{eq:2.2.3a}
\widehat{A(e)}_0=\{\phi\in
\widehat{A(e)}:H^{2\dim(\CB_e)}(\CB_e)^\phi\neq 0\}, 
\end{equation}
and define $\widehat {A(e,\chi)}_0$ to be the representations  of
$A(e,\chi)$ which are restrictions of representations of $A(e)$ in
$\wht{A(e)}_0.$ 
 
For $\phi\in \widehat{A(e)}_0$, we will denote the Springer
 representation by $\mu(\CO,\phi)$. Each representation of $W$ is
 uniquely of the form $\mu(\CO,\phi)$ for some $(\CO,\phi).$  
The correspondence is normalized so that if $e$ is
the principal nilpotent, and $\phi$ is trivial, then
$\mu(\CO,\phi)=sgn.$

Moreover,   the $\mu(\CO',\phi)$ occurring in
$H^*(\CB_e)$ all correspond to $\CO'$ such that
$\CO\subset\ovl{\CO'}.$

Comparing with (\ref{eq:2.2.1}) and (\ref{eq:2.2.2}), we conclude that
\begin{equation}
  \label{eq:2.2.4}
\Hom_W[\mu(\CO,\phi):X(e,\chi,\psi)]=[\phi\mid_{A(e,\chi)}\ : \psi].
\end{equation}

\begin{definition}\label{d:lwt}
Following \cite{BM1}, the $W-$representations in the set 
\begin{equation}\label{eq:2.2.4a}
\{\mu(\CO,\phi):~ [\phi\mid_{A(e,\chi)}\ : \psi]\neq 0\}
\end{equation}
will be called the {\it lowest $W-$types} of $X(e,\chi,\psi).$ When
$\psi=triv,$ we call the $W-$type $\mu(\CO,triv)$ in (\ref{eq:2.2.4}) the
{\it generic} lowest $W-$type.
\end{definition}
\noindent Clearly, the generic lowest $W$-type always appears with multiplicity
one in $X(e,\chi,triv).$

\smallskip

By \cite{KL} and \cite{L5}, if $\nu=0,$
then $X(e,\chi,\psi)$  is \textit{tempered irreducible}, and it has a
unique lowest $W$-type, $\mu(\CO,\psi)$, whose multiplicity is one. If, in
addition, $e$ is an element of a distinguished nilpotent orbit,
$X(e,\chi,\psi)$ is a \textit{discrete series} module.

\smallskip

By \cite{L5}, the module $X(e,\chi,\psi)$
has a unique irreducible subquotient $\oX(e,\chi,\psi)$ characterized
by the fact that it contains each lowest $W-$type $\mu(\CO,\phi)$ with
full multiplicity $[\phi\mid_{A(e,\chi)}\ : \psi]$.

\begin{remark}\label{r:2.2} In the geometric classification, the spherical modules are those of the form
$\overline X(0,\chi,triv).$ The generic modules are $X(e,\chi,triv),$ such that
$X(e,\chi,triv)$ is irreducible (\cite{BM4},\cite{Re}). For the generic
modules, the semisimple element $\chi$
determines the nilpotent orbit $\CO=G\cdot e$ uniquely, according to
(\ref{eq:1.1.3}).
\end{remark}

\subsection{}\label{sec:2.2a}The analogous formula to (\ref{eq:2.2.1}) holds whenever the data $(e,\chi)$
factor through a Levi component $M.$ Let $A_M(e,\chi)$ denote the
component group of the centralizer in $M$ of $e$ and $\chi$. The
following lemma is well known.

\begin{lemma}[1]\label{l:2.2a} The natural map $A_M(e,\chi)\to A(e,\chi)$ is an injection.
\end{lemma}
\begin{proof}
If $T=Z(M)^0$,
then $M=Z_G(T).$ We have that $Z_G(T)\cap Z_G(e)^0=Z_{Z_G(e)^0}(T)$ is
connected, since the centralizer of a torus in a connected algebraic
group is connected. Therefore
$M\cap Z_G(e)^0=Z_M(e)^0.$ 
\end{proof}

 We have:
\begin{equation}
  \label{eq:2.2.5}
\begin{aligned}
&X(e,\chi)=\bH\otimes_{\bH_M} X_M(e,\chi),\ \text{and}\\
&\bH\otimes_{\bH_M} X_M(e,\chi,\tau)= \bigoplus_{\psi\in \widehat{A(e,\chi)}_0}\  [\psi |_{A_M(e,\chi)}\ :\
\tau]\    X(e,\chi,\psi).
\end{aligned} 
\end{equation}

\noindent{\bf Notation.} We write $\Ind_M^G[\pi]$ for the module $\bH\otimes_{\bH_M} \pi.$ 


Define 
\begin{equation}\label{eq:2.4.1a}
M(\nu):=Z_G(\nu),
\end{equation}
 and $P=M(\nu) N$ is such that $\langle\nu,\al\rangle >0$ for all
 roots $\al\in \Delta(\fk n).$  Write $M(\nu)= M_0(\nu) Z(M(\nu)),$ where
 $Z(M(\nu))$ is the center.

\begin{lemma}[2]\label{l:2.2a.1}
$A_{M(\nu)}(e,\chi)=A_G(e,\chi).$ 
\end{lemma}
\begin{proof}
This is because the centralizer of
$e$ is of the form $LU$ with $U$ connected unipotent, and $L$ is the
centralizer of both $e$ and $h.$ It follows that every component of
$A_G(e,\chi)$ meets $L,$ and therefore 
\begin{equation}
  \label{eq:2.2.6}
 A_G(e,\chi)=A_G(e,h,\nu)=A_{M(\nu)}(e,h)=A_{M(\nu)}(e,\chi). 
\end{equation}
\end{proof}

For $\tau\in\widehat{A_{M(\nu)}(e,\chi)}=\widehat{A_{M_0(\nu)}(e)},$
\begin{equation}
  \label{eq:2.2.7}
X_{M(\nu)}(e,\chi,\tau)=X_{M_0(\nu)}(e,h/2,\tau)\otimes\bC_\nu.
\end{equation}
The representation \begin{equation}\sig:=X_{M_0(\nu)}(e,h/2,\tau)\end{equation} is a tempered
irreducible module. Let $\psi\in \widehat{A_G(e,\chi)}$ be the
representation corresponding to $\tau\in\widehat{A_{M(\nu)}(e,\chi)}$
via the identification in lemma \ref{l:2.2a.1}(2). Then
\begin{equation}X(e,\chi,\psi)=\bH\otimes_{\bH_{M(\nu)}} X_{M(\nu)}(e,\chi,\tau).\end{equation}

\begin{definition}In general, whenever $\sig$ is a tempered representation of $\bH_M$ corresponding to the
parameter $(e,h/2,\tau)$, $\tau\in\widehat{A_M(e)}$, and $\nu\in \fk z(\fk m)$, $\chi=h/2+\nu,$
we will write  
\begin{equation}
  \label{eq:2.2.8}
X(M,\sig,\nu):=\Ind_M^G[\sigma\otimes\nu],  
\end{equation}
and call it  a \textit{standard} module also. By equation
(\ref{eq:2.2.5}), it decomposes as 
\begin{equation}
X(M,\sig,\nu)=\bigoplus_{\psi\in\widehat {A_G(e,\nu)}_0}
[\psi|_{A_M(e)}:\tau]~X(e,\chi,\psi).
\end{equation}
If $M=M(\nu)$, then 
\begin{equation}
X(M,\sigma,\nu)=X(e,\chi,\psi),
\end{equation}
where $\psi$ corresponds to $\tau$ as in lemma \ref{l:2.2a.1}(2).
\end{definition}
The terminology is justified
by the fact that $X(M,\sigma,\nu)$ is (via the Borel-Casselman
correspondence)  the $\C I$-fixed vectors of an induced
(standard) module in the classical form of Langlands classification for the
$p-$adic group. If $\langle\nu,\al\rangle\ge 0$ for all
positive roots, then $X(M,\sig,\nu)$, with $M=M(\nu)$, has a unique irreducible
 quotient $\ovl{X}(M,\sig,\nu).$ If $\langle\nu,\al\rangle\le 0$ for
 all positive roots, then
  $X(M,\sig,\nu)$, $M=M(\nu)$ has a unique irreducible submodule
  $\ovl{X}(M,\sig,\nu).$ In the setting of graded Hecke algebras, this
 form of the 
 classification 
  is proved in \cite{E}.

\subsection{}\label{sec:2.3}
Let $\mathfrak z(e,h,f)$ be the centralizer of the triple $\{e,h,f\}$, and
$\mathfrak a_{BC}\subset \fz(e,h,f)$ a Cartan subalgebra such that
$\nu\in\mathfrak a.$ Let $\mathfrak m_{BC}$ be the centralizer of
$\mathfrak a,$ with decomposition  
\begin{equation}
  \label{eq:2.3.1}
 \mathfrak m_{BC}=\mathfrak m_{BC,0}+\mathfrak a_{BC}.
\end{equation}
Then the Lie triple is contained in $[\fk m_{BC},\fk m_{BC}]\subset
[\fk m_{BC,0},\fk m_{BC,0}].$ Thus $\fk m_{BC,0}$ is semisimple (its
center centralizes the 
triple, so must be contained in $\fk a_{BC}.$) So $\fk m_{BC,0}$ is the derived
algebra of $\fk m_{BC},$ and the {nilpotent element  $e$
  is distinguished in $\fk m_{BC,0}.$} The Levi component $\fk m_{BC}$ is the one used in the Bala-Carter
classification of nilpotent orbits, hence the notation. Let $M_{BC},
M_{BC,0}$ be the corresponding groups. The
triple $(e,h/2,\psi)$ with $\psi\in\wht{A_{M_{BC,0}}(e,\chi)}$ determines a
discrete series parameter on $M_{BC,0}.$ Clearly, for any
$\nu\in\fa_{BC},$ if $M(\nu)$ is as in (\ref{eq:2.4.1a}), then $M_{BC}\subset M(\nu).$ 

We are interested in the question of reducibility for the induced
modules $X(M,\sigma,\nu)$, where $M_{BC}\subset M$, and
$\sigma$ is generic.

\begin{proposition}\label{p:2.3}
Let $M$ be a Levi subgroup such 
that $M_{BC}\subset M.$
 Assume $\sig$ is a tempered generic module of $\bH_{M,0}$ corresponding
 to $(e,h/2,triv)$, and consider the standard
 module $X(M,\sig,\nu).$  Set $\chi=h/2+\nu$ for the
infinitesimal character and $\CO=G\cdot e.$ Then $X(M,\sigma,\nu)$ is
reducible if and only if one of the following two conditions is satisfied:
\begin{enumerate}
\item there exists $\CO'$ satisfying $\ovl{\CO'}\supset \CO,$ such that $\CO'$
has a representative $e'$ satisfying $[\chi,e']=e'.$
\item there is no $\CO'$ as in (1), but $A_M(e,\chi)\ne A_G(e,\chi),$ and there
exists a nontrivial character $\psi\in\wht{A_G(e,\chi)}_0$ satisfying 
$[\psi\mid_{A_M(e,\chi)}\ :\ triv]\ne 0.$
\end{enumerate}

\end{proposition}

\begin{proof}
Condition (1) follows from \cite{BM1}. Condition (2) is an
immediate consequence of formula (\ref{eq:2.2.5}).
\end{proof}

\begin{remark}
By equation (\ref{eq:2.2.6}), when $M\supset M(\nu)$, (1) in proposition \ref{p:2.3} is necessary and sufficient.
\end{remark}


\subsection{}\label{sec:2.4} In the next sections we will construct
intertwining operators associated to elements which preserve the data
$(M,\sig).$ 

\medskip
Assume first that $M$ is the Levi component of an arbitrary standard
parabolic subgroup, and $\sig$ a
representation of
$\bH_M$. Let 
\begin{equation}\label{eq:2.4.0}
\fk m=\fk m_0 +\fk a
\end{equation}
 be the Lie algebra of $M,$ with
center $\fa,$ and derived algebra $\fm_0.$ Write
$\fh=\fk t +\fk a$ for the Cartan subalgebra. If 
$w\in W=W(\fg,\fk h)$ is such that $w(\fk m)=\fk m'$ is another Levi
subalgebra (of a standard parabolic subalgebra), choose $w$ to be
{\it minimal} in the double coset $W(M)wW(M').$ Let $w=s_{\al_1}\dots
s_{\al_k}$ be a reduced decomposition. In \cite{BM3} the elements 
\begin{equation}\label{eq:2.4.1}
r_\al=t_{s_\al}\al-1
\end{equation}
are introduced. 
Set 
\begin{equation}\label{eq:2.4.1b}
r_w:=r_{\al_1}\dotsb r_{\al_k}.
\end{equation}
 By
lemma 1.6 in \cite{BM3}, the definition does not depend on the choice
of reduced expression. Because $w$ is minimal in its double coset, it
defines an isomorphism of the root data, and therefore an isomorphism
$a_w:\bH_M\longrightarrow \bH_{M'}.$ Let $w\cdot\sig$ be the
representation of $\bH_{M'}$ obtained from $\sig$ by composing with
$a_w^{-1}.$ Then $r_w$ defines an intertwining operator
\begin{align}
  \label{eq:2.4.2}
&A_w(\sig,\nu):\Ind_M^G(\sig,\nu)\longrightarrow 
\Ind_{M'}^G(w\cdot\sig,w\cdot\nu),\\\notag 
&t\otimes v\mapsto tr_{w}\otimes a_w^{-1}(v),\quad t\in W,v\in\sigma.
\end{align}
For each $\mu\in\wht W,$ $A_w(\sigma,\nu)$ induces an intertwining operator
\begin{equation}
  \label{eq:2.4.3}
  A_{w,\mu}(\sig,\nu):\Hom_W[V_\mu :\Ind_M^G(\sig,\nu)]\longrightarrow
\Hom_W[V_\mu :\Ind_{M'}^G(w\cdot\sig,w\cdot\nu)], 
\end{equation}
which by Frobenius reciprocity can be written as
\begin{equation}
  \label{eq:2.4.3a}
  A_{w,\mu}(\sig,\nu):\Hom_{W(M)}[V_\mu :\sig]\longrightarrow
\Hom_{W(M')}[V_\mu :w\cdot\sig]. 
\end{equation}
These operators are defined for all $\nu$ not just real ones. We
assume that $\nu$ is complex for the rest of the section.
Define the element $\kappa_w\in\bH$:
\begin{equation}\label{eq:2.4.3b}
\kappa_w=\displaystyle{\prod_{\beta>0,\ w\beta<0}
  (\beta^2-1)}.
\end{equation}

\begin{proposition}\label{p:2.4} 
The operators $A_{w}(\sig,\nu)$ have the following properties: 
\begin{enumerate}
\item $A_{w}(\sig,\nu)$ is polynomial in $\nu.$ 
\item $A_{w^{-1}}(w\cdot\sig,w\cdot\nu)\circ
  A_{w}(\sig,\nu)=(\sig\otimes\nu) (\kappa_w)$, where $\kappa_w$ is
  defined in (\ref{eq:2.4.3b}). Furthermore $\kappa_w$ is an element of the center of $\bH_M,$ so the right
hand side is a scalar multiple of the identity.
\item Assume $\sig$ is hermitian. Then the hermitian dual of $A_w(\sigma,\nu)$
  is $A_w(\sig,\nu)^*=A_{w^{-1}}(w\cdot\sig,-w\cdot\ovl{\nu}).$
\end{enumerate}  
\end{proposition}
\begin{proof}
Part (1) is clear from the definition.
Part (2) follows from the fact that in the Hecke
algebra, \begin{equation}r_\al^2=(t_\al\al-1)^2=\al^2-1,\text{ for }\al\text{
    a simple root.}\end{equation}
 The fact
that $\kappa_w$ is in the center of $\bH_M$ follows from the
fact that $w$ is shortest in the double coset. Let $\fk p=\fm+\fn$ and
$\fk p'=\fm'+\fn'$ be the standard parabolic subalgebras. If
$\beta\notin\Delta(\fm)^+$ is such that $w\beta<0,$ then $-w\beta\in\Delta(\fn').$ 
If $\al\in\Pi(\fm)$ is a simple root, then
$-w(s_\al(\beta))=-s_{w(\al)}(w\beta)\in\Delta(\fn'),$ because $w\al$
is a simple root of $\fm',$ and so preserves $\fn'.$ 

For part
(3) we recall from \cite{BM3} that the hermitian dual of
$\Ind_M^G(\sig,\nu)$ is $\Ind_M^G(\sig,-\ovl{\nu})$ with the pairing
given by
\begin{equation}
  \label{eq:2.4.4}
  \langle x\otimes v_1\one_\nu,y\otimes v_2\one_{-\ovl\nu}\rangle= \langle 
\sig(\ep_M(y^{-1}x)v),w\rangle.
\end{equation}
In this formula, $x,y\in W/W(M),$ and $\ep_M$ is the projection of
$\bC[W]$ onto $\bC[W(M)].$ We omit the rest of the proof.
\end{proof}
\begin{remark}
The formula for the pairing in (\ref{eq:2.4.4}) follows from
\cite{BM3} which uses the fact that $\bH$ has the $*$ operation given by
equation (\ref{eq:2.1.7}).
\end{remark}
\subsection{}\label{sec:disc} We still assume that $\nu$
is complex. We specialize to the case when $\sig$ is
a generic discrete series for $\bH_{M_0}.$ So in this case $M=M_{BC},$ and the
induced module in the previous section is $X(M,\sig,\nu).$ As in equation (\ref{eq:2.2.5}), this decomposes into a
direct sum of standard modules $X(M,\sig,\nu,\psi),$ with $\psi\in
\widehat{A(e,\nu)}_0.$ In particular, $X(M,\sigma,\nu,triv)$ is the
only generic summand.

In this case we normalize the operators $A_w$ so that
they are $Id$ on $\mu(\CO,triv),$ and restrict them to the subspace
$X(M,\sig,\nu,triv)$.  We denote the normalized operators
$\CA_w(\sigma,\nu)$ (on $X(M,\sigma,\nu)$), respectively
$\CA_w(\sigma,\nu,triv)$ (on $X(M,\sigma,\nu,triv)$). They define, by
restriction to $\Hom$ spaces as in (\ref{eq:2.4.3a}), operators
$\CA_{w,\mu}$ for $\mu\in\widehat W.$   

\smallskip

Assume that $w$ decomposes into $w=w_1w_2$,  such that
$\ell(w)=\ell(w_1)+\ell(w_2)$ ($\ell(w)$ is the length of
$w$). {The fact that $r_w=r_{w_1}r_{w_2}$ (\cite{BM3}) implies  one of the most
important properties of the operators $\CA_w$, the factorization}: 
\begin{align}
  \label{eq:2.4.6}
\CA_{w_1w_2}(\sig,\nu)&=
\CA_{w_1}(w_2\cdot\sig,w_2\cdot\nu)\circ\CA_{w_2}(\sig,\nu),\text{ and
  similarly}\\\notag
\CA_{w_1w_2}(\sig,\nu,triv)&=
\CA_{w_1}(w_2\cdot\sig,w_2\cdot\nu,triv)\circ\CA_{w_2}(\sig,\nu,triv).
\end{align}

\begin{proposition}
  \label{p:normalized}
Assume that $\fm$ is the Levi component of a maximal standard parabolic
subalgebra. Then $\CA_w(\sigma,\nu,triv)$ does not have any poles in the region of
$\nu$ satisfying $\langle\nu,\beta\rangle\ge 0$ for all $\beta>0$ such that
$w\beta<0$.
\end{proposition}
\begin{proof}
Either $w\nu=\nu,$ or else $\langle w\nu,\beta\rangle<0$ for
all $\beta\in \Delta(\fn)$.
By \cite{E}, if $\langle \beta,\nu\rangle >0$ for all
$\beta\in\Delta(\fn),$ then $X(M,\sig,\nu,triv)$
has a unique irreducible quotient, while if
$\langle\nu',\beta\rangle<0$ for all $\beta\in\Delta(\fn'),$ then
$X(M',\sig',\nu',triv)$ has a unique irreducible
submodule. 
By the results in \cite{KL} and \cite{L5} this is the unique
subquotient containing $\mu(\CO,triv).$ Thus $\CA_w$ maps
$X(M,\sig,\nu,triv)$ onto $\ovl{X}(M,\sig,\nu,triv).$ 

Assume that $\CA_w$ has a pole
of order $k>0$ at $\nu_0$ with $Re\nu_0>0.$ Then
$(\nu-\nu_0)^k\CA(\sig,\nu,triv)$ extends analytically to
$\nu=\nu_0,$ and is nonzero. Its image is disjoint from
$\ovl{X}(M,\sig,\nu,triv),$ which contradicts the fact that
$\oX(w\cdot M, w\cdot\sig,w\cdot\nu)$ is the unique irreducible submodule of
$X(w\cdot M,w\cdot\sig,w\cdot\nu,triv).$ 

Now suppose $\CA_w$ has a pole at $\nu_0$
with $Re\nu_0=0.$ We use the analogues of (1)-(3) from proposition
\ref{p:2.4}; the relation (2) implies that for the normalized
operators we have:
\begin{equation}\CA_{w^{-1}}\circ \CA_w=Id.\end{equation} Write
\begin{equation}
  \label{eq:2.4.7}
  \begin{aligned}
&\CA_w(\sig,\nu,triv)=(\nu-\nu_0)^k[A_0+(\nu-\nu_0)A_1+\dots ],\quad\text{where }
A_0\neq 0, \text{ and}\\
&\CA_{w^{-1}}(w\cdot\sig,w\cdot\nu,triv)=
(\CA_{w}(\sig,-\ovl\nu,triv))^*=(-\nu+\nu_0)^k[A_0^*+(-\nu+\nu_0)A_1^*+\dots ].   
  \end{aligned}
\end{equation}
Then if $k<0$ relation (2) in proposition \ref{p:2.4} implies 
$A_0^*A_0=0,$ which is a contradiction.
\end{proof}

\subsection{}\label{sec:factoring} We present a standard technique for
factorizing intertwining operators (see \cite{SV} for the setting of
real reductive groups).

\begin{definition}\label{d:adjacent}
We say that two Levi components $\fm,\ \fm'$ are adjacent, if either
$\fm=\fm',$ or there is a
Levi component $\Sigma$ such that $\fm,\fm'\subset\Sigma$ are maximal Levi
components conjugate by $W(\Sigma).$
\end{definition}

\begin{lemma}\label{l:factoring}
Let $w$ be such that $w(\fm)=\fm'$, and $w$ minimal in the double
coset $W(M)wW(M').$ Then there is a chain of adjacent Levi components
$\fm_0=\fm,\dots, \fm_k=\fm'$.    
\end{lemma}

\begin{proof}
We do an induction on the length of $w.$ If $\fm =\fm'$ and $w=1,$
there is nothing to prove. Otherwise there is $\al$ simple such that
$w\al<0.$ Then let $\Sigma_1$ be the Levi component with simple roots
$\Delta(\fm)\cup\{\al\}.$ Then $ww_1^{-1}$ has shorter length, and the
induction hypothesis applies. 
\end{proof}

We will always consider minimal length chains of Levi subalgebras. The main reason for
these notions is the
following. Let $w_{\Sigma_i}^0$ be the longest element
in $W(\Sigma_i),$ and $w_i$ be the shortest element in
$W(\fm_{i-1})w_{\Sigma_i}^0W(\fm_i).$ Then
we can write
\begin{equation}
  \label{eq:2.4.8}
  w=w_k\cdot\ \dots\ \cdot w_1,\quad
  \CA_w=\CA_{w_k}\circ\dots \circ \CA_{w_1}.
\end{equation}
The $\CA_{w_i}$ are induced from the corresponding operators for
maximal Levi components, and so proposition \ref{p:normalized}
applies.

\begin{theorem}
  \label{t:normalized}
The intertwining operators $\CA_w$ have the following properties.
\begin{enumerate}
\item $\CA_w(\sigma,\nu,triv)$ is analytic for $\nu$ such that
  $\langle Re~\nu,\beta\rangle\ge 0$ for all $\beta>0$ such that
  $w\beta<0.$
\item $\CA_{w^{-1}}(w\cdot\sig,w\cdot\nu,triv)\circ\CA_w(\sig,\nu,triv)=Id,$
\item  $\CA_w(\sig,\nu,triv)^*=\CA_{w^{-1}}(w\cdot\sig,-w\cdot\ovl{\nu},triv).$
\end{enumerate}
\end{theorem}
\begin{proof}
  This follows from proposition \ref{p:normalized} and lemma \ref{l:factoring}.
\end{proof}
\begin{remark}\label{r:normalized}
If there exists an isomorphism
$\tau:w\cdot\sig\longrightarrow\sig,$ we compose the intertwining
operators $\CA_w$ with $(1\otimes\tau).$ For simplicity, we will
denote these operators by $\CA_w$ also. 
 
If in fact
$$w\cdot\sig\cong\sig,\ w\cdot\nu=-\ovl\nu,$$ the operator $\CA_w$ gives rise
to a hermitian form. This is  because
$\Ind_M^G(\sigma,-\ovl\nu)$ is the hermitian dual of $\Ind_M^G(\sigma,\nu).$   
\end{remark}

\subsection{}\label{sec:int}
We assume  that $\nu$ is \textbf{real}. Let $x\in G$ stabilize
$\{e,h,f\}.$ Then we can choose the Cartan subalgebra $\mathfrak a_{BC}$ of
$\fz(e,h,f)$ so that it is stabilized by
$x.$ Furthermore, since $x$ stabilizes $\mathfrak m_{BC}$ and
$\mathfrak m_{BC,0},$ there is a Cartan subalgebra $\mathfrak
t\subset\mathfrak m_{BC,0},$ stabilized by $x.$ Let
\begin{equation}
  \label{eq:int.1}
  \mathfrak h:=\mathfrak t + \mathfrak a_{BC}
\end{equation}
be the Cartan subalgebra of $\fk m_{BC}.$
We can also choose a Borel subalgebra of $\mathfrak m_{BC}$ containing
$\fh$ which is stabilized
by $x.$ So $x$ gives rise to a Weyl group element $w_x,$ the shortest
element in the double coset $W_{M_{BC}}xW_{M_{BC}}.$ Thus we get an
intertwining operators $\CA_{w_x}$ by the construction in sections
\ref{sec:2.4}-\ref{sec:factoring}. 

If $x\cdot\nu=-\nu$ and
$\tau:x\cdot\sig\overset{\cong}\longrightarrow \sig$, by remark
\ref{r:normalized}, $\CA_{w_x}$  
 gives rise to a hermitian form.

\subsection{}\label{sec:2.6}
We apply the construction of section \ref{sec:int} in the following special
case.  
Let $\oal$ be a simple root of $\fk a_{BC}\subset\fz(e,h,f).$ Let $x_{\oal}\in
Z(e,h,f)^0$ be an element inducing the reflection $s_{\oal}$ on $\fk
h.$ 
Then $x_\oal$ stabilizes $\fm_{BC}.$ The element $x_\oal$ may need to be
modified by an element in $M_{BC,0}$ so as to stabilize $\ft$ as well. Then
it gives rise to a Weyl group element $w_\oal,$ shortest in the double
coset $W_{M_{BC}}x_\oal W_{M_{BC}},$ and to an intertwining operator $\CA_{w_{\bar\al}}$. The new $x_\oal$ may
not fix the Lie triple. But since it modified the original element by
one in $M_{BC,0},$ there is an isomorphism 
$\tau_\oal:w_\oal\sig\overset{\cong}\longrightarrow \sig.$ 

Then, as in remark \ref{r:normalized}, we have a normalized intertwining
operator 
\begin{equation}
  \label{eq:2.6.1}
  \CA_\oal: X(M_{BC},\sig,\nu,triv)\longrightarrow
  X(M_{BC},\sig,w_\oal\nu,triv). 
\end{equation}

\subsection{}\label{sec:2.9}
We construct intertwining operators for another class of elements
normalizing $\sig.$  We consider an $M\supset M_{BC}$, and write 
  $\fk m=\fk m_0+\fk a,$ $\fk a\subset\fk a_{BC},$ as in equation
  \ref{eq:2.4.0}. Let $A$ and $H$ be the
Cartan groups corresponding to $\fa,$ $\fh$, and let 
$\sig$ be a tempered representation of  $\bH_{M,0}.$  Define
\begin{equation}
  \label{eq:2.9.1}
  \begin{aligned}
&  N(\fa):=\{ w\in W : w\fa=\fa \},\\ 
&  C(\fa,M):=\{w\in N(\fa) : w(\Delta^+(\fm))=\Delta^+(\fm) \}.
  \end{aligned}
\end{equation}

The following formula is a particular case (which we need here for the
construction of intertwining operators) of a more general result that  we
postpone to section \ref{sec:2.8}. 

\begin{lemma}\label{l:2.9}
\begin{equation}
  \label{eq:2.9.4}
\text{(1)}\quad  N(\fa)=C(\fa,M)\ltimes W(M).
\end{equation}
 \label{l:2.10}
  \begin{equation}
\text{(2)}\quad    N_G(\fa)/M\cong C(\fa,M).
  \end{equation}  
\end{lemma}

\begin{proof}
(1) From (\ref{eq:2.9.1}), we see that
\begin{equation}
  \label{eq:2.9.2}
  N(\fa)=C(\fa,M)\cdot W(M).
\end{equation}
In fact, as in the proof of lemma \ref{l:2.8},
\begin{equation}
  \label{eq:2.9.3}
C(\fa,M)\cap W(M)=\{1\},  
\end{equation}
and $W(M)$ is a normal subgroup, because any element $xmx^{-1}$ with $x\in N_G(\fa)$ centralizes $\fa,$ so must be in $M.$

\smallskip

(2) Clearly $M$ is normal in $N_G(\fa).$ 
  Let $n\in N_G(\fa).$ Then $n\fa=\fa,$ and $n\fh=\fh'=\ft'+\fa.$ There
  is an element $m\in M$ such that 
  \begin{equation}
    \label{eq:2.10.5}
    mn\fh=\fh,\qquad mn(\Delta^+(\fm))=\Delta^+(\fm).
  \end{equation}
Thus the $M$-coset of $mn$ is in $C(\fa,M).$ This map is a group
homomorphism, and an isomorphism onto $C(\fa,M).$ 
\end{proof}

If $c\in C(\fa,M)$ is such that $c\cdot\sig\cong\sig,$ then by the
construction in section \ref{sec:disc}, in particular remark
\ref{r:normalized}, there is a normalized intertwining operator 
\begin{equation}\label{eq:2.9.5a}
\CA_c(\sigma,\nu):X(M,\sigma,\nu)\to X(M,\sigma,c\cdot\nu).
\end{equation}
and for every $(\mu,V_\mu)\in \wht W$, this induces an operator
$\CA_{c,\mu}(\sigma,\nu)$ as in equations (\ref{eq:2.4.3}) and (\ref{eq:2.4.3a}).

\subsection{}\label{sec:2.10}
We put the constructions in the previous sections together. We
consider the case when $M=M_{BC}.$

Let $W(\fz,\fa_{BC})$ denote the Weyl group of $\fa_{BC}$ in
$\fz:=\fz(e,h,f).$  Denote by $\bH(\fz)$ the graded Hecke algebra
constructed from the root system of $\fz.$ 
In this section we study the relation of $W(\fz,\fa_{BC})$ with $C(\fa_{BC},M_{BC})$,
in particular we show that $C(\fa_{BC}, M_{BC})$ contains naturally a subgroup
isomorphic to $W(\fz,\fa).$  Elements in this subgroup give
rise to $\bH(\fz)$-intertwining operators of the (spherical) principal series
$X_{\bH(\fz)}(0,\nu)$ of $\bH(\fz)$, as well as $\bH$-intertwining
operators for $X(M,\sig,\nu)$ by equation (\ref{eq:2.9.5a}).

Set
\begin{equation}
  \label{eq:2.10.1}
  \wti A_{BC}:=Z_{Z(e,h,f)}(\fa_{BC}).
\end{equation}
Then $A_{BC}\subset \wti A_{BC},$ so there is a surjection
\begin{equation}
  \label{eq:2.10.2}
  N_{Z(e,h,f)}(\fa_{BC})/A_{BC}\longrightarrow N_{Z(e,h,f)}(\fa_{BC})/\ti A_{BC}.
\end{equation}
Furthermore, there is an injective group homomorphism,
\begin{equation}
  \label{eq:2.10.3}
  W(\fz,\fa_{BC})=N_{Z(e,h,f)}(\fa_{BC})/\ti A_{BC}\longrightarrow N_G(\fa_{BC})/Z_G(\fa_{BC})=N_G(\fa_{BC})/M_{BC}.
\end{equation}

\begin{proposition}
  \label{p:2.10}\ 

\begin{enumerate}
\item The composition of the map in lemma \ref{l:2.10} with
  the map in (\ref{eq:2.10.3}) gives an injective homomorphism
$$W(\fz,\fa_{BC})\hookrightarrow C(\fa_{BC},M_{BC}).$$
\item The composition of the map in lemma \ref{l:2.10} with the map in
(\ref{eq:2.10.2}) 
$$A_G(e)\ltimes W(\fz,\fa_{BC})\cong N_{Z(e,h,f)}(\fa_{BC})/A_{BC}\longrightarrow
N_G(\fa_{BC})/M_{BC}\stackrel{\cong}\longrightarrow C(\fa_{BC},M_{BC})
$$ is onto.
\end{enumerate}
\end{proposition}

\begin{proof}
Part (1) is clear. For part (2), let $n\in N_G(\fa_{BC})$ be given. Then $n$ induces an automorphism of $\fm_{BC}.$ So it
maps the Lie triple $\{e,h,f\}$ into another Lie triple $\{e',h',f'\}.$
The Levi component is of the form
\begin{equation}
  \label{eq:2.10.6}
  \fm_{BC}\cong \fm_1\times gl(a_1)\times\dots \times gl(a_r),
\end{equation}
with $\fm_1$ simple, not type A. The nilpotent orbit is a
distinguished one on $\fm_1,$ and the principal nilpotent on the
$gl(a_i)$ factors. Since any automorphism  of a simple (or even a
reductive algebra with simple derived algebra) maps a distinguished
orbit into itself, there is $m\in M_{BC},$ such that $mn$ stabilizes the
triple $\{e,h,f\}.$ Thus every $M_{BC}$ coset of $N_G(\fa_{BC})$ contains a
representative in $N_{Z(e,h,f)}(\fa_{BC}),$ which is the claim of the proposition.
\end{proof}

The image of the map in part (2) consists of elements which stabilize
$\sig.$ Thus each element in $x\in A_G(e)\ltimes W(\fz,\fa_{BC})$ gives
rise to an intertwining operator
\begin{equation}
  \label{eq:2.10.7}
  \CA_x(\sig,\nu):X(M,\sig,\nu)\longrightarrow
  X(M,\sig,w_x\cdot\nu). 
\end{equation}
normalized to be $Id$ on $\mu(\CO,triv)$. In particular we get an
action of $A_G(e,\nu)$ on $X(M,\sig,\nu)$. This action should coincide
with the one defined geometrically, but we have not been able to verify this.

\subsection{}\label{sec:2.10b} 
Denote by $W(\fz(\CO))\cong W(\fz,\fa_{BC})$ the abstract Weyl group of
$\fz(e,h,f)$, and similarly $A(\CO)$ for the component group,  and set $W(Z(\CO)):=A(\CO)\ltimes W(\fz(\CO)).$

We will restrict now to the case of
hermitian Langlands parameters, $(M,\sigma,\nu),$ where $\sigma$ is a
discrete series for $M$. Recall that this means that $M=M_{BC},$ but
in order to simplify notation, we drop the subscript in this section.  As before, there
must exist $w\in W$ such that
\begin{equation}\label{eq:2.10b.1}
wM=M,\quad w\sigma\cong\sigma, \quad{\text{and }} w\nu=-\nu.
\end{equation}
For $(\mu,V_\mu)\in\widehat W$, section \ref{sec:2.9} defines an operator
$\CA_\mu(\sigma,\nu)$ 
(by Frobenius reciprocity)  on the
  space $\text{Hom}_{W(M)}(V_\mu,\sigma)$. The group $C(\fa,M)$ acts on $W(M)$,
  and therefore on $\widehat {W(M)}$, and  preserves $\sigma.$   

 Let
  $\mu_{M}(\CO,triv)$ be the unique lowest $W(M)$-type of $\sigma$. Then
  \begin{equation*}
    \begin{aligned}
&\Hom_W[\mu(\CO,triv):X(M,\sig,\nu)]=\Hom_{W(M)}[\mu(\CO,triv):\sig]=\\
&\Hom_{W(M)}[\mu(\CO,triv):\mu_M(\CO,triv)]=1.     
   \end{aligned}
  \end{equation*}
In the calculations in section
  \ref{sec:5}, we will only consider $W-$types $\mu$
  in  $X(M,\sigma,\nu)$ with the property that
  $$\Hom_{W(M)}[\mu:\sigma]=\Hom_{W(M)}[\mu:\mu_{M}(\CO,triv)].$$
  We need the fact that $C(\fa,M)$ preserves
  $\mu_{M}(\CO,triv)$. Since $\sigma$ is tempered, this is equivalent to
  the fact that $C(\fa,M)$ preserves $\sigma$.

\begin{definition}\label{d:2.10b} Let $\sigma$ be a discrete series
  for $\bH_M$ parametrized by $\CO$, where $M=M_{BC}$ of $\CO.$ The
  space $\text{Hom}_{W(M)}(\mu,\sigma)$ has the structure of  a 
  representation of $C(\fa,M)$ and via the map from proposition
  \ref{p:2.10}, it is a $W(\fz(\CO))$-representation, and a
  $W(Z(\CO))$-representation, which we will denote $\rho(\mu)$,
  respectively $\rho'(\mu).$   
  \end{definition}

\subsection{}\label{sec:2.10c} 
In view of lemma \ref{l:2.2a}, for every Levi subgroup $M_{BC}\subset
M\subset G,$ one has
$A_M(e)\subset A_G(e).$ 
In a large number of cases, $A_G(e)=A_{M_{BC}}(e)$, and
analyzing the standard modules $X(M_{BC},\sig,\nu)$ with $\sig$ a discrete series 
is sufficient. In the other cases, we also need intermediate Levi
components $M'$ with the property that $A_G(e,\nu)=A_{M'}(e).$

\medskip 

Consider the Levi subgroups $M$ with Lie algebras $\fk m$ subject to the
  conditions:
\begin{enumerate}
\item $e\in \fk m;$
\item $A_{G}(e)=A_{M}(e).$
\end{enumerate}
We call the nilpotent orbit $\CO$ {\it quasi-distinguished} if the 
minimal subalgebra with respect to conditions (1) and (2) is $\fg.$
Note that every distinguished $\CO$ is also quasi-distinguished.

\begin{proposition}
If $\CO$ is a quasi-distinguished nilpotent orbit, then $\fz(\CO)$ is
a torus.
\end{proposition}

\begin{proof} It is easy to verify the statement case by case using the
  Bala-Carter (\cite{Ca}) classification of nilpotent orbits.

\end{proof}

\begin{definition}
If $\sigma$ is a tempered irreducible module parametrized by 
a quasi-distinguished $\CO,$ we call $\sigma$ a {\it limit of discrete
  series}. 
\end{definition}
With this
definition, any discrete series is a limit of discrete series. We list
next the limits of discrete series, which are not discrete series,
and appear for various Levi subalgebras of
$E_6,$ $E_7,$ and $E_8.$ Clearly, if $\sigma$ is a limit of discrete
series for $\fm$ in $E_6,$ it will also be considered in $E_7$ and
$E_8.$ Therefore, to eliminate redundancy, we will list a pair
$(\fm,\CO)$ only for the smallest algebra for which appears. For $\fm$
of type $D$, we also give the notation of the orbit as a partition. In
type $A$, the only quasi-distinguished orbit is the principal orbit.

\begin{tiny}
\begin{center}
\begin{longtable}{|c|c|c|}
\caption{Limits of discrete series}\label{table:lds}\\
\hline
\multicolumn{1}{|c|}{Type of $\fg$} &
\multicolumn{1}{c|}{Levi subalgebra $\fm\subset\fg$} 
& \multicolumn{1}{c|}{Nilpotent in $\fm$} \\ \hline 
\endfirsthead

\multicolumn{3}{c}%
{{  \tablename\ \thetable{} -- continued from previous page}} \\
\hline
\multicolumn{1}{|c|}{Type of $\fg$} &
\multicolumn{1}{c|}{Levi subalgebra $\fm\subset \fg$} 
& \multicolumn{1}{c|}{Nilpotent in $\fm$}
 \\ \hline 
\endhead


\hline \hline
\endlastfoot

$E_6$ &$E_6$  &$D_4(a_1)$\\
\hline
     &$D_4$  &$A_2=(3311)$\\
\hline
\hline
$E_7$ &$E_7$  &$E_6(a_1)$\\
\hline
     &$E_7$  &$A_4+A_1$\\
\hline
     &$D_6$  &$D_5(a_1)=(7311)$\\
\hline
     &$D_6$  &$A_4=(5511)$\\
\hline
     &$D_6$  &$A_3+A_2=(5331)$\\
\hline
\hline
$E_8$ &$E_8$  &$D_7(a_1)$\\
\hline
     &$E_8$  &$E_6(a_1)+A_1$\\
\hline
     &$E_8$  &$D_7(a_2)$\\
\hline
     &$E_8$  &$D_5+A_2$\\
\hline
     &$D_7$  &$D_6(a_1)=(9311)$\\
\hline
     &$D_7$  &$D_6(a_2)=(7511)$\\ 
\hline
     &$D_7$  &$D_4+A_2=(7331)$\\
\hline
     &$D_7$  &$A_4+2A_1=(5531)$
\end{longtable}
\end{center}
\end{tiny}

As before, consider the module $X(M_{BC},\sigma,\nu),$ $\sigma$ generic discrete
series. For the calculations in section \ref{sec:5}, whenever
$A_{M_{BC}}(e,\nu)\neq A_G(e,\nu),$ we can find a pair $(M',\sigma')$, where
$M'$ is a Levi
component $M'\supset M_{BC},$ with the following properties: 
\begin{enumerate}
\item $A_{M'}(e,\nu)=A_G(e,\nu),$
\item $\sigma'$ is the generic summand of $\Ind_{M_{BC}}^{M'}[\sigma]$ and
  $\sigma'$ is a limit of discrete series for $M'$,
\item $X(M_{BC},\sigma,\nu,triv)=X(M',\sigma',\nu)$.
\end{enumerate}

\section{The $0$-complementary series}\label{zero}

\subsection{}\label{sec:2.11}
We specialize to the case of spherical principal series. Some of
these results were 
already presented in the introduction in the setting of the split
$p$-adic group.

Consider the {\it principal series module} 
\begin{equation}\label{2.11.2}
X(\chi)=\mathbb H\otimes_{\mathbb H}  \mathbb C_\chi,\ \chi\in\fk h.
\end{equation}

As a $W-$representation, $X(\chi)$ is isomorphic to $\bC[W].$ In
particular, the module $X(\chi)$ has a unique generic subquotient and
a unique spherical subquotient
$\overline X(\chi)$. We will refer to a semisimple element
$\chi$ as unitary if $\overline X(\chi)$ is unitary. 

The construction of intertwining operators as presented in sections
\ref{sec:2.4}-\ref{sec:2.9} {becomes simpler} in this setting.
Consider the {intertwining operator} given by $r_{w_0}$, where
$w_0$ is the longest
element in the Weyl group, and normalized so that it is $Id$ on the
trivial $W-$type. Since
the operator only depends on  $\chi$, we will simply denote it by
$\CA_{w_0}(\chi): X(\chi)\to X(w_0\chi)$. 

If $\chi$ is dominant (i.e., $\langle
  \chi,\alpha\rangle \ge 0$ for all  roots $\alpha\in \Delta^+$) the
  image of $\CA_{w_0}(\chi)$ is $\overline X(\chi)$.
Moreover, $\overline X(\chi)$ is hermitian if and only if $w_0
  \chi=-\chi$. It is reducible if and only if
  $\langle\al,\chi\rangle=1$ for some $\al\in\Delta^+.$ The generic
  subquotient is also spherical if and only if $X(\chi)$ is irreducible.

Note that $r_{w_0}=r_{\al_1}\dotsb r_{\al_k}$ acts on the right and
therefore, each $\al_j$ in the decomposition into $r_{\al_j}$'s can be
replaced by the scalar $\langle \al_j, w_j\chi\rangle$, where $w_j=s_{j+1}s_{j+2}\dotsb
s_k$ in the intertwining operator
$\CA_{w_0}(\chi)$. For every $(\mu,V_{\mu})\in\widehat W,$ denote 
\begin{equation}\label{eq:2.11.2a}
a_\mu(\chi)=\CA_{w_0,\mu}(\chi): V_\mu^*\longrightarrow V_\mu^*.
\end{equation}

\begin{remark}\label{longint} Assume $w_0 \chi=-\chi$. 
The
  hermitian form on $\overline X(\chi)$ is positive definite 
  if and only if all the operators $a_\mu(\chi)$
  are positive semidefinite. 
\end{remark}

More precisely, the operators $a_\mu(\chi)$ are characterized by the fact that, in
the decomposition $a_\mu(\chi)=a_{\mu,\al_1}(w_1\chi)\dotsb
a_{\mu,\al_k}(w_k\chi)$ coming from the reduced expression for $w_0$ as
above (see also section \ref{sec:factoring}) ,
\begin{equation}\label{2.12.3}
a_{\mu,\al_j}(\nu)=\left\{\begin{matrix} 1, &\text{on the }
    (+1)\text{-eigenspace of } s_{\al_j} \text{ on } V_\mu^*\\
\frac {1-\langle\al_j,\nu\rangle}{1+\langle\al_j,\nu\rangle},&\text{on the }
    (-1)\text{-eigenspace of } s_{\al_j} \text{ on } V_\mu^*.\end{matrix}\right. 
\end{equation}
\noindent If $\al$ is a simple root, we have the formula (\cite{BM3})
$t_{s_\al} r_w=r_wt_{s_{w^{-1}\al}}.$
>From this, since $s_{w^{-1}\al}=w^{-1}s_\al w$, it follows that 
\begin{equation}
t_{w}r_{w}=r_{w}t_{w},\quad\text{for any }w\in W.
\end{equation}
In particular, for $w=w_0$, we obtain  that every
$a_\mu(\chi)$ preserves the $(+1)$, respectively $(-1)$,
eigenspaces of $w_0$ on $\mu^*.$ 

\subsection{}\label{sec:alcoves}
Consider $\chi$ in the $(-1)$-eigenspace of
$w_0$. In order to determine if $\chi$ is unitary, one would have to
compute the operators $a_\mu(\chi)$ on the W-type $\mu$.  An
operator $a_\mu(\chi)$ has constant signature on any facet in the
arrangement of hyperplanes 
\begin{equation}\label{eq:hyp}
\langle \chi,\al\rangle=1,\ \al\in \Delta^+\text{ and }\langle
\chi,\al\rangle=0,\ \al\in \Pi,
\end{equation} 
in the dominant Weyl chamber $\C C$ of $\fh$ (see theorem 2.4 in \cite{BC}). 

The
$0$-complementary series (definition \ref{d:1.3}(2)) is a union of open regions in this
arrangement of hyperplanes.

Recall that the {\it fundamental alcove} $\C C_0$ is the set
\begin{equation}
\C C_0=\{\chi\in\C C:\langle \al,\chi\rangle<1,\text{ for all }\al\in
\Delta^+\}.
\end{equation}
If $W_\text{aff}$ denotes the affine Weyl group, an {\it alcove} is, by
definition, any open region in $\C C$ which is  $W_\text{aff}$-conjugate
with $\C C_0.$ Clearly, any alcove is a simplex.

\medskip

The main results of this section are summarized next.

\begin{theorem}\label{p:4.1} The $0$-complementary series are:
\noindent\begin{itemize}
\item as in theorem \ref{classical} for types $A,\ B,\ C,\ D$,
\item as in proposition \ref{p:g2} for types $G_2,\ F_4$,
\item the hermitian $\chi$ ($w_0\chi=-\chi$) in the union of the $2$ alcoves in section \ref{sec:genspherE6} for
  type $E_6,$
\item the union of the $8$ alcoves in section \ref{sec:genspherE7} for
  type $E_7,$
\item  the union of the $16$ alcoves in section
  \ref{sec:genspherE8} for type $E_8.$
\end{itemize}
\end{theorem}

\subsection{}

 We recall the description of the $0$-complementary series for Hecke
 algebras of classical types.

\begin{theorem}(\cite{BM3},\cite{Ba2})\label{classical} 
The  parameters $\chi=(\nu_1,\nu_2,\dotsc,\nu_n)$  in
  the 0-complementary series are:
\begin{itemize}
\item[{\bf A}]: $\chi=(\nu_1,\dotsc,\nu_k,-\nu_k,\dotsc,-\nu_1)$  or
  $(\nu_1,\dotsc,\nu_k,0,-\nu_k,\dotsc,-\nu_1)$, with
  $0\le\nu_1\le\nu_2\le\dots\le\nu_k<\frac 12$. 
\item[{\bf C}]: $0\le\nu_1\le\nu_2\le\dotsb\le\nu_n<\frac 12.$
\item[{\bf B,D}]: there exists $i$ such that
  $0\le\nu_1\le\dotsb\le\nu_i<1-\nu_{i-1}<\nu_{i+1}<\dotsb<\nu_n<1$,
  and between any $\nu_j<\nu_{j+1}$, $i\le j<n$, there is an odd
  number of $(1-\nu_l)$, $1\le l<i.$
\end{itemize}
\end{theorem}

\subsection{} 

We will also need the description of the $0$-complementary series for
the Hecke algebras of type $G_2$ and $F_4$.  We use the roots $\al_1=(\frac 23,-\frac 13,-\frac 13)$ and
$\al_2=(-1,1,0)$ for $G_2$ and $\al_1=(1,-1,-1,-1)$,
$\al_2=(0,0,0,2)$, $\al_3=(0,0,1,-1)$, $\al_4=(0,1,-1,0)$ for $F_4$.

\begin{proposition}[\cite{C}]\label{p:g2} 
\begin{enumerate}
\item If $\bH$ is of type $G_2$ and
  $\chi=(\nu_1,\nu_1+\nu_2,-2\nu_1-\nu_2)$, $\nu_1\ge 0$, $\nu_2\ge 0$,
  is a spherical parameter, the 0-complementary series is
  \begin{equation}\{3\nu_1+2\nu_2< 1\}\cup\{2\nu_1+\nu_2<1<3\nu_1+\nu_2\}.\end{equation}

\item If $\bH$ be of type $F_4$ and
  $\chi=(\nu_1,\nu_2,\nu_3,\nu_4)$, $\nu_1-\nu_2-\nu_3-\nu_4\ge 0$,
  $\nu_2\ge \nu_3\ge\nu_4\ge 0$,
  is a spherical parameter, 
 the 0-complementary series  is
  \begin{equation}\{2\nu_1< 1\}\cup\{\nu_1+\nu_2+\nu_3-\nu_4<1<\nu_1+\nu_2+\nu_3+\nu_4\}.\end{equation}

\end{enumerate}
\end{proposition} 

Part (1) of proposition \ref{p:g2} was first established in \cite{M}.

\subsection{}\label{sec:7.1}

In the rest of this section, we determine the $0$-complementary series for types $E_7$ and $E_8$. (The method also applies in
type $D_n$, where we recover known results of \cite{BM3} and
\cite{Ba2}). For $E_6,$ the argument needs to be modified slightly due
to the fact that $w_0\neq -1,$ but it is essentially the same. It is
presented in detail in section 3.5 of \cite{Ci2}.

\medskip

Assume $G$ is
of type $D_{2m}$, $E_7$ or $E_8$. The notation for $W$-types is as in
\cite{Ca}. One important nonunitarity criterion
that we will use is the following. 
Let $M$ be a Levi subgroup of type
$A_2$. The nilpotent orbit
$A_2$ has two lowest $W-$types, $\mu(A_2,triv)$ and $\mu(A_2,sgn)$ as follows:
\begin{align}\label{eq:lwt}
& D_{2m}:\ &(2m-2,1)\times (1),\ (2m-2)\times (11)\notag\\
& E_7: \ &56_a',\ 21_a\notag\\
& E_8: \ &112_z,\ 28_x,
\end{align}
on which operators $\CA(St,\nu)$  for the standard module
$X(M,St,\nu),$ $M=A_2,$ have opposite signature whenever $Z_G(\nu)=M$.
(The details for this type of calculation are in lemma \ref{l:3.1} and
in section \ref{sec:5}.)  
This means that for all $\nu$ such that $Z_G(\nu)=M=A_2,$ the module
$\ovl X(A_2,St,\nu)$ is not unitary. Therefore:

\begin{lemma}\label{l:7.1}
The generic module $X(A_2,St,\nu)$ is not unitary for all
parameters $\nu$ such that $Z_G(\nu)=A_2.$
\end{lemma}

\subsection{} \label{sec:7.2} 
Recall the hyperplane arrangement (\ref{eq:hyp}). The
        connected components of the complement of this hyperplane
        arrangement in $\C C$ will be called {\it regions}.  
Inside any region $\C F$, the
intertwining operators $a_\mu(\chi)$ are isomorphisms, therefore
their signature is constant in $\C F$.

We recall first that the unbounded (open) regions are not
unitary. This is a well-known result. A proof in the setting of the
Hecke algebra can be found in \cite{BC}, 3.3.  

\begin{lemma} \label{l:7.2} If the open region $\C F$ is unbounded, and
  $\chi\in\C F$, then the operator
  $a_\mu(\chi)$, for $\mu$ the reflection representation,  is not
  positive definite.
\end{lemma}

\subsection{}\label{sec:7.3} 
Recall the relation
of partial order on $\Delta^+$: 
\begin{equation}\label{eq:7.3.1}
\beta_1>\beta_2\text{ if }
\beta_1-\beta_2\text{ is a sum of positive roots.}
\end{equation}
If $\beta_1\ge\beta_2$ or
$\beta_2>\beta_1$, then $\beta_1,\beta_2$ are said to be {\it
  comparable}, otherwise they are {\it incomparable}. A subset of
incomparable positive roots is called an {\it antichain}. Two roots in
an antichain, being incomparable, must have nonpositive inner product.

\smallskip

If $\Pi=\{\al_1,\dots,\al_n\}$ are
the simple roots and a positive root $\beta$ is
$\beta=\sum_{i=1}m_i\al_i$, call $\sum_{i=1}m_i$ the {\it height of}
$\beta$.  We consider the positive roots ordered in (\ref{eq:7.3.1}) on {\it levels} given by
the height. The simple roots are level $1$ and the highest root is
level $h-1$, where $h$ is the Coxeter number ($h=2(n-1)$ in $D_n$,
$h=18$ in $E_7$ and $h=30$ in $E_8$).

\smallskip

Any region $\C F$ is an intersection of
half-spaces $\langle\beta,\chi\rangle>1$ or
$\langle\beta,\chi\rangle<1,$ for all $\beta\in \Delta^+$, and
$\langle\al,\chi\rangle\ge 0$, for all $\al\in \Pi$. Let 
\begin{align}
&\delta(\C F)
\text{ be the set of maximal roots among the roots }\beta<1\text{ on
}\C F,\text{ and}\\\notag
&\delta'(\C F)\text{ be the set of minimal roots among the roots }\beta'>1\text{ on }\C
F.
\end{align}
The following proposition is clear (and well-known).

\begin{proposition}
For every region $\C F,$  both $\delta(\C F)$ and $\delta'(\C F)$ are antichains in $\Delta^+$.
Moreover, the correspondences $\C F\to\delta(\C F)$ and $\C F\to
\delta'(\C F)$ are bijections between the set of regions and the set of antichains of positive roots. 
\end{proposition}

\begin{remark}A region $\C F$ is infinite if and
only if $\delta'(\C F)\cap\Pi\neq \emptyset.$ 
\end{remark}
\begin{proof}
Let
$\chi\in\C F$ and assume $\langle\al,\chi\rangle>1,$ for some simple
root $\al$. If $\omega_\al$ is the corresponding coweight, for all
$t\ge 0$, $\langle\beta,\chi+t\omega_\al\rangle>1+t\ge 1,$ if
$\beta>\al$, and
$\langle\beta',\chi+t\omega_\al\rangle=\langle\beta',\chi\rangle$, for
all $\beta'$ incomparable to $\al$. This implies that
$\chi+t\omega_\al$ is in $\C F$, for all $t\ge 0$.
\end{proof}

The walls of the region $\C F$ (regarded as a convex polytope) are
given by the hyperplanes $\beta=1$, for $\beta\in \delta(\C F)\cup
\delta'(\C F)$, and possibly by $\al=0$, for some simple roots
$\al$. 

Note that a simple root $\al$ does not give a wall $\al=0$ of $\C F$ if and only there
exists a root $\beta\in \delta(\C F)$ such that $\beta+\al$ is also a
root (in the simply-laced, equivalently,
$\langle\beta,\al\rangle=-1$). This is because in this case, for all
$\chi\in\C F$, $\langle\beta,\chi\rangle<1<\langle\beta+\al,\chi\rangle,$
so one cannot set $\langle\al,\chi\rangle=0$ without crossing a
hyperplane $\beta=1$. Similarly, one can formulate such a condition
with the roots in $\delta'(\C F).$

\subsection{} \label{sec:7.4} The signature of
intertwining operators $a_\mu(\chi)$ on the walls of the dominant
Weyl chamber is known by unitary induction from smaller groups. In
$D_{2n}$, by setting a simple root equal to $0$, we get a parameter
unitarily induced irreducible from $D_{2n-2}+A_1$, in $E_7$ from
$D_6$, and in $E_8$ from $E_7$. In particular, a region $\C F$, which
has a wall $\al=0$, for some simple root $\al$, is unitary if and only
if the parameters on the wall $\al=0$ are induced from a unitary
region in the smaller group. This is a well-known argument, see lemma
\ref{l:4.5}. 

We will need the following information about the antichains formed of
mutually orthogonal roots. We call such subsets {\it orthogonal antichains}.

\begin{lemma}\label{l:7.4}
If $\Delta$ is a simply laced root system, the maximal cardinality of an
orthogonal antichain in $\Delta^+$ equals the number of positive
roots at level $[\frac {h(\Delta)+1}2],$ where $h(\Delta)$ is the Coxeter number.
\end{lemma}

\begin{proof} We verified this assertion case-by-case. It also follows
  from the main theorem in 
  \cite{So}, which states that every antichain is $W-$conjugate to a 
  subset of the Dynkin diagram of $\Delta$.
\end{proof}

\begin{proposition}\label{p:7.4} Any unitary region $\C F$ has a wall
  of the form $\al=0$, for some simple root $\al$.
\end{proposition}

\begin{proof} In view of lemma \ref{l:7.2}, we may assume that $\C F$
  is a finite region, that is, a convex polytope. Assume by contradictions
  that all the walls of $\C F$ are $\beta=1$, for
  $\beta\in \delta(\C F)\cup \delta'(\C F).$ 

There are two cases which we treat separately:

a) $\C F$ has a dihedral angle of $\frac {2\pi}3.$

b) All dihedral angles of $\C F$ are non-obtuse.

\smallskip

a) Let $\beta_1\in\delta(\C F),\beta_2\in\delta'(\C F)$ be such that
$\langle\beta_1,\beta_2\rangle=-1$ and they give adjacent walls of $\C
F.$ Let $\chi_0$ be a parameter such that $\chi_0\in
(\beta_1=1)\cap(\beta_2=1)\cap \overline{\C F}$, but
$\langle\beta,\chi_0\rangle\neq 1$, for any
$\beta\notin\{\beta_1,\beta_2\}$. This is possible, otherwise there
should exist a positive root $\beta$  such that $\beta_1=1,\beta_2=1$
implies necessarily $\beta=1$. In particular, $\{\beta_1,\beta_2,\beta\}$ are
linearly dependent over $\mathbb Z$. Since we are in the simply laced
case, one must be a sum of the other two roots, but then they cannot
all be equal to $1$ simultaneously.

The principal series $X(\chi_0)$ is reducible. The generic factor  is
parametrized by the nilpotent orbit $A_2$.  
By lemma \ref{l:7.1}, this factor is not unitary, and therefore the region $\C F$
is also nonunitary.

\medskip

b) Assume that all dihedral angles of $\C F$ are non-obtuse. A
classical theorem of Coxeter implies in our case that $\C F$ must be 
in fact a simplex. 

We are therefore in the following situation:
\begin{align}\label{eq:7.4.2}
&\langle\beta_1,\beta_2\rangle=0, &\text{if }\beta_1,\beta_2\in\delta(\C
  F)\text{ or }\beta_1,\beta_2\in\delta'(\C F),\\\notag
&\langle\beta,\beta'\rangle\in\{0,1\}, &\text{if }\beta\in\delta(\C
F)\text{ and }\beta'\in \delta'(\C F).
\end{align}
The antichains $\delta(\C F)$ and $\delta'(\C F)$ are orthogonal.  Set
$k=|\delta(\C F)|$, $k'=|\delta'(\C F)|,$ and $k+k'=n+1$, where $n$ is
the rank of $\Delta.$ By lemma \ref{l:7.4}, $k\le m+1$ for $D_{n}$
($n=2m$) and $k\le 4$ for $E_7$, $E_8,$ and same for $k'$. This
immediately gives a contradiction for $E_8$ ($k+k'\le 8<9$). In $E_7$,
the only possibility is $k=k'=4$, and in $D_{2m}$, $k=m+1,k'=m$
($k=m,k'=m+1$ is analogous). It remains to analyze these cases.

\smallskip

Fix $\beta'\in\delta'(\C F)$. For all $\beta\in \delta(\C F)$,
$\langle\beta,\beta'\rangle\in\{0,1\}$. If for all $\beta\in\delta(\C F),$
$\beta'$ is not comparable to $\beta$ (in particular
$\langle\beta,\beta'\rangle=0$), $\{\beta_1,\dots,\beta_k,\beta'\}$ would be
an antichain of $k+1$ orthogonal roots, contradiction. Thus, there exists
$\beta$ such that $\beta'>\beta$. Let $\al$ be a simple root such that
$\langle\beta,\al\rangle=-1$, and $\beta'\ge\beta+\al>\beta$ (this is
always possible in the simply-laced case). Since $\beta<1$ is a wall,
$\beta+\al>1$, so necessarily $\beta'=\beta+\al$ (otherwise $\beta'>1$
would not be a wall).

To summarize, for each $\beta'\in\delta'(\C F)$, there exists
$\beta\in\delta(\C F)$
such that $\beta'-\beta$ is a simple root. Similarly, for each
$\beta\in\delta(\C F)$ there exists $\beta'\in\delta'(\C F)$ with
$\beta'-\beta$ a simple root. 

If $\al$ is a simple root, $\al=0$ is not a wall of $\C F$ if and
only if there exists $\beta\in\delta$ such that $\beta<1<\beta+\al$ in
$\C F$. From the discussion above, the region $\C F$ is not adjacent to
the walls of the dominant chamber if and only if for any $\al$
simple root, there exists $\beta\in\delta(\C F)$ and
$\beta'\in\delta'(\C F)$ such
that $\beta'-\beta=\al$.  

If this is the case, we are looking at a
bipartite graph with $k+k'$ vertices (roots) $\delta(\C
F)\cup\delta'(\C F)$ and at least
$n=k+k'-1$ edges (simple roots), such that any vertex has degree $\ge 1$. We
would like to claim that this graph is connected. The only way to fail
connectedness is if there exists a complete (bipartite) subgraph
$\{\beta_1,\beta_2\}\cup\{\beta_1',\beta_2'\}$. This means that there exist
simple roots $\al_1,\dots,\al_4$ such that 
\begin{equation}
\beta_1'=\beta_1+\al_1=\beta_2+\al_2,\quad
\beta_2'=\beta_1+\al_3=\beta_2+\al_4.
\end{equation}
Then
\begin{equation}
1=\langle\beta_1',\beta_2\rangle=\langle\beta_1+\al_1,\beta_2\rangle=\langle\al_1,\beta_2\rangle,
\end{equation}
and similarly $\langle\al_4,\beta\rangle=1$. But then
\begin{equation}
0=\langle\beta_1',\beta_2'\rangle=\langle\beta_1+\al_1,\beta_2+\al_4\rangle=2+\langle\al_1,\al_4\rangle,
\end{equation}
so $\langle\al_1,\al_4\rangle=-2$, which gives a contradiction (simply-laced case).

If the graph is connected, it means that $\delta(\C F)$, respectively
$\delta'(\C F)$ are formed of the positive roots on the same level of
the root system, and moreover the two levels are consecutive. But this
is false by inspection.

\end{proof}

\begin{corollary}\label{c:7.4} A parameter $\chi$ is in the
  $0$-complementary series if and only if $\chi$ can be deformed
  irreducibly to a point $\chi_0,$ such that $X(\chi_0)$ is unitarily
  and irreducibly induced from a parameter in the $0$-complementary
  series on a proper Levi component.
\end{corollary}

We also remark that part b) of the proof of proposition \ref{p:7.4}
can be applied to the regions $\C F$ for which the antichains
$\delta(\C F)$ and $\delta'(\C F)$ are formed only of roots at levels
greater than or equal to $\frac {h(\Delta)}2$ (since the sum of two such
roots cannot be a root, their inner product is non-negative). Then,
all such regions are adjacent to the walls of the dominant Weyl
chamber. By induction, we will find that all unitary regions are of
this form.

\subsection{}\label{sec:relevant} An important fact is that for the determination of the 0-complementary
series, one only needs to know the signature of intertwining operators
on a small number of $W-$types (and not on all of $\bC[W]$). In
addition to its intrinsic interest, we will need
this information in section \ref{sec:4} and for the calculations in
section \ref{sec:5}. (See section \ref{sec:4.4} for the explanation.)

\begin{definition}\label{d:relevant}  Assume the root system of $\bH$
  is simple. 
The following $W-$types are called {\it 0-relevant}:

\begin{tiny}
\begin{center}
\begin{longtable}{|c|c|}\caption{$0$-relevant
    $W-$types}\label{table:relevant}\\ 
\hline
{\bf Type} &{\bf $0$-relevant $W-$types}\\
\hline

$A$  &$\{(n-1,1)\}$\\
\hline
$B,C,D$ &$\{(n-1)\times (1),\ 
  (1,n-1)\times (0)\}$ or $\{(n-1)\times (1), (n-2)\times (2) \}$\\
\hline
$G_2$ &$\{2_1,2_2\}$ or $\{2_1,1_2,1_3\}$\\
\hline
$F_4$ &$\{4_2,9_1\}$\\
\hline
$E_6$ &$\{6_p,20_p\}$\\
\hline
$E_7$ &$\{7_a',27_a\}$ or $\{7_a',21_b'\}$\\
\hline
$E_8$ &$\{8_z,35_x\}$\\
\hline
\end{longtable}
\end{center}
\end{tiny}
\end{definition}

\begin{proposition}\label{r:relevant} A parameter $\chi$ is in the $0$-complementary series
  if and only if the operators $a_\mu(\chi)$ are positive definite on
  all $0$-relevant $\mu.$
\end{proposition}

In every list of $0$-relevant $W-$types, the reflection
representation, $refl$ is present. Recall lemma \ref{l:7.2} which says
that the signature of $refl$ in any infinite region is not positive
definite. Note also, that for every type of $W$ not type $A$, the
second $W-$type in a possible list of $0$-relevant, appears in $Sym^2(refl).$ (In fact
for exceptional groups, this is the unique nontrivial $W-$type in
$Sym^2(refl).$) 

\begin{proof} For type $A$, the claim follows easily from the fact
  that, in this case, every region (\ref{sec:7.2}) is adjacent to a
  wall of the dominant Weyl chamber. 

For types $B,C,D$, the proof is in \cite{BC}. The proof is
  conceptual, and it is based on a some simple calculations of
  determinants of intertwining operators. An essential step in  the
  proof is the fact that the centralizer $\fz(\CO)$ of the nilpotent
  orbit $\CO=A_1=(2,2,1,\dots,1)$ has a factor of type $A_1.$

Types $G_2$ and $F_4$ can be found in \cite{C}, and type $E_6$ is in
\cite{Ci2}. A similar argument as in the classical types works here as
well; the argument uses the fact that the centralizer of
$\CO=A_1$ is of type $A$, more precisely, $A_1$ for $G_2,$ $A_3$ for $F_4$, and $A_5$
for $E_6.$

{For $E_7$ and $E_8$ one cannot use the same argument. The
  difference is that the 
  centralizers $\fz(\CO)$ for $\CO=A_1$ do not contain a factor of type
$A$. 
The proof of the proposition and corollary in section \ref{sec:7.4},
shows that a spherical parameter $\chi$ is in the $0$-complementary if
and only if the operators $a_\mu(\chi)$ are positive definite of
\begin{enumerate}
\item[(i)] $\mu\in\{7_a',27_a,56'_a,21_a\}$ or $\mu\in\{7_a',21_b',56'_a,21_a\}$ for $E_7$,
\item[(ii)] $\mu\in\{8_z,35_x,112_z,28_x\}$  for
  $E_8,$ 
\end{enumerate}
in other words, on a strictly larger set than what we called
$0$-relevant in table \ref{table:relevant}. In order to show that in
fact, it is sufficient to consider only the signatures of the
$0$-relevant $W$-types for $E_7$, $E_8$, we used a
computer calculation.   
We will only need to use this finer information for $E_7$ at one place in this
paper, namely in section \ref{sec:match1}, for the
nilpotent $A_1\subset E_8$ (whose centralizer is $E_7$). Proposition
\ref{r:relevant} for $E_8$  will not be needed in the sequel.}
\end{proof}

\section{Extended Hecke algebras}\label{sec:extended}


\subsection{}\label{sec:2.8}
The goal is to construct graded Hecke algebras for certain disconnected groups.

\medskip
Suppose $\fk G$ is an arbitrary linear algebraic group with connected
component $\fk G^0,$ and component group $R:=\fk G/\fk G^0.$   
Let $\bH$ denote the graded Hecke algebra
associated to $\fk G^0.$ Choose a pair  $(B,H),$ where
$B$ is a Borel subgroup, and $H\subset B$ a Cartan subgroup in $\fk G^0.$
Denote by $W:=N_{\fk G^0}(H)/H,$ the Weyl group of $\fk G^0$.

\begin{lemma}\label{l:2.8}
  \begin{equation*}
    N_{\fk G}(H)/H\cong R\ltimes W.
  \end{equation*}
\end{lemma}
\begin{proof}
Let 
\begin{equation}
  \label{eq:2.8.1}
 R':=\{ g\in \fk G : gH=H,\ gB=B\}. 
\end{equation}
We will show that $R\cong R'/H.$ 
It is clear that $H\subset R'$ and $R'\cap \fk G^0=H.$ Furthermore, 
 \begin{equation}
   \label{eq:2.8.2}
   N_{\fk G}(H)=R\cdot N_{\fk G^0}(H).
 \end{equation}
Finally, $R'$ meets every component of $\fk G.$ Indeed, if $g\in \fk G,$ then
$g\cdot B=B',g\cdot H=H',$ where $(B',H')$ is another pair of the same
type as $(B,H).$ Then there is $g_0\in \fk G^0$ such that
$(g_0B,g_0H)=(B,H).$ Then $g_0g\in R',$ and belongs to the same
component as $g.$ The proof follows.
\end{proof}

If $g\in \fk G$, then $(g\cdot B,g\cdot H)$ is another pair of Borel
and Cartam subgroups. Thus
there exist an element $x\in \fk G^0$ such that $xg$ stabilizes the pair
$(B,H).$ Then $xg$ determines an automorphism $a_g$ of the based root
datum. If $g\in \fk G^0,$ then $a_g=Id.$ Suppose $g_1,\ g_2\in \fk G,$ and
$x_1,\ x_2\in \fk G^0$ are such that $x_1g_1,\ x_2g_2$ stabilize the pair
$(B,H).$ Then the fact that
\begin{equation}
  \label{eq:2.8.3}
  x_1g_1x_2g_2=(x_1g_1x_2g_1^{-1})(g_1g_2),\qquad
  x_1g_1x_2g_1^{-1}\in \fk G^0,
\end{equation}
implies that
\begin{equation}
  \label{eq:2.8.4}
  a_{g_1}a_{g_2}=a_{g_1g_2}.
\end{equation}
Thus the group $R\cong R'/H$ maps to the group of automorphisms of
the root datum for $\fk G^0,$ and therefore maps to the automorphism group
of $\bH,$ the corresponding affine graded Hecke algebra. We will
identify $R$ with this automorphism group.
\begin{definition}\label{d:2.8.1} Let $\bH$ denote the graded Hecke
algebra for the   root datum of $\fk G^0$ (as in (\ref{eq:2.1.5})). We
define $\bH'$ to be the semidirect product 
\begin{equation}
  \label{eq:2.8.4a}
  \bH':=\bC[R]\ltimes \bH,
\end{equation}
where the action of $R$ on $\bH$ is induced by the $a_g$ defined
earlier.
\end{definition}

\subsection{}\label{sec:2.8a}
We are interested in the spherical representations of $\bH'.$ This is a 
special case of Mackey induction. Set \begin{equation}\C K':=R\ltimes W,\text{ and }\C
K:=W.\end{equation} 
 A representation of $\bH'$ is called 
{\it spherical}, if it contains the trivial representation of $\C K'.$ 

\begin{lemma}\label{l:2.8.1}
The center of $\bH'$ is $\bA^{\C K'}.$
\end{lemma}
\begin{proof}
  This is clear from proposition \ref{p:center}.
\end{proof}

For every $\nu\in\fh^*,$ we make the following notation:
\begin{align}\label{eq:2.8b.1}
&R(\nu)=\text{the centralizer of }\nu\text{ in }R,\\\notag
&\bA'(\nu)=\bC[R(\nu)]\ltimes\bA,\\\notag
&\bH'(\nu)=\bC[R(\nu)]\ltimes\bH,
\end{align}
where $\bA$ is the abelian part of $\bH$ (as in (\ref{eq:2.1.5})), and
the action of $R(\nu)\subset R$ is as in definition \ref{d:2.8.1}. 

Consider 
\begin{equation}  \label{eq:2.8b.3}
X'(\nu)=\bH'\otimes_{\bA'(\nu)}\bC_\nu.
\end{equation}

\begin{proposition}\label{l:2.8.2}
Assume $(\pi,V)$ is a spherical irreducible representation of $\bH'.$
The multiplicity of the trivial representation of $\C K'$ is 1.   
\end{proposition}

\begin{proof}
Let $\nu$ be a weight of $V$ under $\bA,$ spanned by $v_\nu,$ and
define
\begin{equation}
  \label{eq:2.8a2}
  R_\nu:=\{r\in R\ :\ \pi(r)v_\nu=v_\nu\}.
\end{equation}
Set $\bA_\nu':=\bC[R_\nu]\ltimes\bA.$ 
Then $V$ is a quotient of $\bH'\otimes_{\bA_\nu}\bC_\nu,$ via the map
$x\otimes\one_\nu\mapsto \pi(x)v_\nu.$  But as a $\C K'$ module,
\begin{equation}
  \label{eq:2.8.a3}
\bH'\otimes_{\bA_\nu}\bC_\nu=\sum_{\mu\in\wht{\C K}} V_\mu\otimes
  (V_\mu^*)^{R_\nu}. 
\end{equation}
Thus the trivial representation occurs exactly once in
$\bH'\otimes_{\bA_\nu}\bC_\nu,$ and the claim follows.
\end{proof}
\begin{corollary}
  \label{c:2.8a}
  \begin{equation*}
    R_\nu=R(\nu).
  \end{equation*}
\end{corollary}
\begin{proof}
Let $V$ denote the spherical irreducible quotient of
$\bH'\otimes_{\bA_\nu}\bC_\nu,$ as in the proof of proposition \ref{l:2.8.2}. Consider the subspace
\begin{equation}
  \label{eq:2.8a.4}
  \{ \sum_{y\in R(\nu)} ky\otimes \one_\nu\}_{k\in\C K'}\subset 
\bH'\otimes_{\bA_\nu}\bC_\nu,
\end{equation}
This is $\bH'-$invariant, and isomorphic to
$X'(\nu)$ from (\ref{eq:2.8b.3}). Since by the analogue of
(\ref{eq:2.8.a3}) $X'(\nu)$ is spherical, we get a nontrivial homomorphism
(hence surjective)
\begin{equation}
  \label{eq:2.8a.5}
X'(\nu)\longrightarrow V.  
\end{equation}
The claim follows from the fact that the stabilizer of $\one_\nu$ in
$R$ is $R(\nu).$ 
\end{proof}

\subsection{}\label{sec:2.8b}

There is a natural extension of the Langlands classification for spherical
modules to $\bH'$. We will not make use of it in an essential way in this
paper, rather it is listed here in order to make clearer the analogy
between the description of $\CO$-complementary series (section
\ref{sec:4}, especially \ref{sec:4.2a}-\ref{sec:4.3}) and the
spherical unitary dual of the extended Hecke algebra constructed from
the centralizer $Z(\CO)$ (see \ref{sec:2.8e}). 

\begin{proposition}
  \label{p:2.8}
Every irreducible spherical module of $\bH'$ is of the form 
$$
L'(\nu):=\bH'\otimes_{\bH'(\nu)} L(\nu)
$$
Two such modules $L'(\nu)$ and $L'(\nu')$ are equivalent, if
and only if $\nu$ and $\nu'$ are in the same orbit under $\C K'$.
 
If
$\nu\ge 0,$ then $X'(\nu)$ has a unique irreducible quotient $L'(\nu)$,
if $\nu\le 0,$ then $X'(\nu)$ has a unique irreducible submodule $L'(\nu)$.
\end{proposition}

\begin{proof}
The proof is based on the Langlands classification for $\bH$ and the
restriction formulas listed below. We will omit the details of
the proof. Corollary \ref{c:2.8a} implies that the restriction to $\bH$ of
$X'(\nu)$ is
\begin{equation}
  \label{eq:2.8a.6}
  X'(\nu)|_{\bH}=\sum_{r\in R/R(\nu)}
\bH\otimes_{\bA}\bC_{r\nu},
\end{equation}
Moreover, if $L_0$ is any spherical $\bH$-module in the restriction
$L'(\nu)|_{\bH}$, then
\begin{equation}\label{eq:2.8a.7}
L'(\nu)|_{\bH}=\sum_{r\in R/R(\nu)} r\cdot L_0.
\end{equation}  
\end{proof}

\begin{corollary}\label{c:2.8b}
\begin{enumerate}
\item If $L'(\nu)$ is hermitian, but $L(\nu)$ is not, then the form on
  $L'(\nu)$ is indefinite.
\item $L'(\nu)$ is unitary if and only if $L(\nu)$ is unitary.
\end{enumerate}
\end{corollary}
\begin{proof}
  If $L'(\nu)$ is unitary, then so is every factor of its
  restriction to $\bH;$ these are the $L(k\nu)$ with $k\in\C K'.$ 
Also, if a factor $L(k\nu)$ is not hermitian, its hermitian dual
occurs in the decomposition, and necessarily the hermitian form on
$L'(\nu)$ cannot   be positive definite. If on the other hand $L(\nu)$
is unitary, then all the $L(k\nu)$ occurring in the decomposition
(\ref{eq:2.8a.6}) are unitary as well. 
\end{proof}
\subsection{}\label{sec:2.8c} 
We can extend the definition of intertwining operators
to this setting. Assume $\xi w\in R\ltimes W.$
Then, similarly to section \ref{sec:2.11}, we define a spherical $\bH'$-operator 
\begin{equation}\label{eq:2.8.8}
A'_{\xi w}(\nu):X'(\nu)\to X'(\xi w\nu),\ x\otimes
\one_\nu\mapsto x\xi r_w\otimes \one_{\xi w\nu}.
\end{equation} 
The operator $\CA'_{\xi w}$ is $A'_{\xi w}$ normalized to be the
identity on the trivial $\C K'-$type. 
For every $\C K'$-type $\mu'$,
this induces an operator 
\begin{equation}\label{eq:2.8.9}
a'_{\xi w,\mu'}(\nu):\Hom_{\C K'}[\mu':X'(\nu)]\longrightarrow\Hom_{\C
  K'}[\mu':X'(\xi w\nu)]. 
\end{equation}

\begin{remark}\label{r:2.8c}
When $w\nu=-\nu,$  the $\bH'$-operator $\CA'_w(\nu)$ gives rise to a hermitian
form on $\Hom_{\C K'}[\mu':X'(\pm\nu)]$ which can be naturally
identified with the form induced by the $\bH$-operator $\CA_w(\nu)$ on
$\Hom_{R(\nu)}[\mu':triv]=((\mu')^*)^{R(\nu)},$
\end{remark}

\subsection{}\label{sec:2.8e} The definitions in the previous sections
can be applied 
 to centralizers of nilpotent orbits. Let $\CO$ be a nilpotent orbit in $\fg$,
  and $Z(\CO)$ be the centralizer in $G$ of a Lie triple $\{e,h,f\}$ of $\CO,$
  with identity component $Z(\CO)^0.$ We will denote by
  $\bH(Z(\CO))$, respectively $\bH(\fz(\CO))$, the Hecke algebras
  $\bH',$ respectively $\bH$, from definition \ref{d:2.8.1}. In this
  particular case, we have:
\begin{align}\label{eq:2.8e.1}
\C K&=W(\fz(\CO)),\\\notag
\C K'&=W(Z(\CO)),\\\notag
R&=A_G(e),\\\notag
R(\nu)&=A_G(e,\nu).
\end{align}

By corollary \ref{c:2.8b}, one can identify the spherical unitary dual
of $\bH(Z(\CO))$ with that of $\bH(\fz(\CO)).$ 

\subsection{}\label{sec:2.8d} We present an interesting instance of
the construction. Assume the root system $\Delta$ is 
simple and it has roots of two lengths. Let $c:\Pi\to \bZ_{\ge 0}$ be
a function, such that $c(\al)=c(\al')$ whenever $\al$ and $\al'$ are
$W-$conjugate. One defines the graded Hecke algebra $\bH_c$ with
parameter $c$ as in section \ref{sec:2.1}, in particular
(\ref{eq:2.1.6}), but with commutation relation  
\begin{equation}\label{eq:2.8d.1}
\om t_s=t_ss(\om)+c(\al)\langle\om,\check\al\rangle,\quad s=s_\al,\
\om\in\fh^*. 
\end{equation}
Consider the case 
\begin{equation}
c(\al)=\left\{\begin{matrix} 1, &\al \text{ long root}\\
                             0, &\al \text{ short root}.
              \end{matrix}\right.
\end{equation}
Denote the corresponding graded Hecke algebra by $\bH_{1,0},$ and let
$\Delta_l\subset \Delta$ denote the subset of long roots, which is a root
(sub)system, and $\Pi_l$ be the simple roots in $\Delta_l$. (Note that
$\Pi_l\not\subset \Pi$, in fact $\text{rank }\Delta_l=\text{rank }
\Delta.$)  Let $W(\Delta_l)$ be  the corresponding Weyl group, and let 
$W_s$ denote the reflection subgroup of $W$ generated by the simple
short roots in $\Pi.$ Then $W_s$ acts on $\Delta_l$, and on
$W(\Delta_l)$ by conjugation. 

\begin{lemma}
$W=W_s\ltimes W(\Delta_l).$
\end{lemma}

\begin{proof}
This follows from the classification of simple root systems.
\end{proof}

Let $\bH(\Delta_l)$ denote the graded Hecke algebra corresponding to the
root datum $(\C X,\check{\C X}, \Delta_l,\check \Delta_l, \Pi_l).$  We can
apply the construction in (\ref{eq:2.8.4a}) with $\bH=\bH(\Delta_l)$
and $R=W_s.$ 

\begin{proposition}
$\bH_{1,0}\cong \bC[W_s]\ltimes \bH(\Delta_l).$
\end{proposition} 

\begin{proof}
In view of the definitions with generators and relations, one only
needs to check  that if $\beta\in \Pi_l,$ then equation
(\ref{eq:2.8d.1}) holds with $s=s_\beta.$ There exists a reflection
$s$ in a simple short root  and $\al\in \Pi$ (long root) such that
$\beta=s(\al),$ therefore $t_\beta=t_st_{s_\al} t_s.$ Using this, it
is straight-forward to check that $\om
t_{s_\beta}=t_{s_\beta}s_\beta(\om)+\langle\om,\check\beta\rangle.$ 
\end{proof}

\begin{remark}
If $\Delta$ is simple, the possible cases are:
\begin{enumerate}
\item $\bH(C_n)_{1,0}=\bC[S_n]\ltimes \bH(A_1^n);$
\item $\bH(B_n)_{1,0}=\bC[S_2]\ltimes \bH(D_n);$
\item $\bH(G_2)_{1,0}=\bC[S_2]\ltimes \bH(A_2);$
\item $\bH(F_4)_{1,0}=\bC[S_3]\ltimes \bH(D_4).$ 
\end{enumerate}
The cases (1), with $n\le 3$, and (2)-(4) all appear as Hecke algebras
$\bH(Z(\CO)).$ 
\end{remark}

\section{Main results}\label{sec:4} In this section we present
the main results of this paper. The explicit calculations (for type $E_8$)
are presented in sections \ref{sec:5} and \ref{exceptions}. We only
consider modules with real infinitesimal characters.

\subsection{}\label{sec:4.1} Recall $\CO\subset \fg,$ where $\fg$ is
of type $E_6,\ E_7,\ E_8.$ Let $\{e,h,f\}$ be a Lie triple for $\CO$,
and let $X(e,\chi,triv)$ be a generic hermitian
representation. Recall the centralizer $Z(\CO)$ with Lie algebra
$\fz(\CO),$ and the decomposition $\chi=h/2+\nu.$ The algebra $\fz(\CO)$
is a product of simple algebras and a torus.

By definition \ref{d:1.2}, the complementary series attached to
$\CO$ is the set of all $\chi=\frac 12h+\nu$ such that the generic module
$X(e,\chi,triv)$ is unitary (and irreducible). The parameter
$\nu\in\fz(\CO)$ parametrizes a spherical module for the Hecke algebra
$\bH(\fz(\CO))$, and by section \ref{sec:extended}, also a spherical module
for the Hecke algebra $\bH(Z(\CO)).$

\begin{theorem}\label{t:4.1}
The parameter $\chi=h/2+\nu$ is in the complementary series attached to
$\CO$ if and only if the corresponding
parameter $\nu$ is in the $0$-complementary series of
$\bH(\fz(\CO))$. 

The $0$-complementary series for the Hecke algebras of simple types are
listed in proposition \ref{p:4.1}. 

The following exceptions occur:
\begin{itemize}
\item $\CO=A_1+\wti A_1$ in $F_4,$
\item $\CO=A_2+3A_1$ in $E_7,$
\item $\CO\in\{A_4+A_2+A_1, A_4+A_2, D_4(a_1)+A_2, A_3+2A_1, A_2+2A_1,
  4A_1\}$ in $E_8.$
\end{itemize}
In all the exceptions, but $\CO=4A_1$ in $E_8$, the complementary
series attached to $\CO$ is smaller than the $0$-complementary series
of $\bH(\fz(\CO)).$ The explicit description is recorded in section \ref{tables}.
\end{theorem}

In the rest of this section, we present the elements of the proof.

\subsection{}\label{sec:3.3} 
The starting case is that of intertwining
operators for induced modules from Levi components of {\it maximal}
parabolic subalgebras. We would like to relate these operators with
operators for Hecke algebras of rank one. 

First we need to record some
results about the reducibility of standard modules. Let $P=MN$
($\fk p=\fm+\fn$) be a maximal parabolic, and
$X(M,\sigma,\nu)$ be a standard module. Using proposition \ref{p:2.3},
we can easily find  the reducibility points of
$X(M,\sigma,\nu),\ \nu>0.$ The answer is given in theorem \ref{t:3.3} below.
Its nature is related to conjectures of Langlands.

Let $\{e,h,f\}\subset\fm$ be a Lie triple parameterizing the tempered
module $\sigma.$ Then $\fk n$ is
a module for the  $\fk{sl}(2,\bC)$ generated by $\{e,h,f\}.$ Let $\al$
be the unique simple root not in $\Delta(\fm),$ and $\check\om$ the
corresponding coweight, which commutes with $\{e,h,f\}$. The
eigenvalues of $\check\om$ on $\fn$ are of the form $1,2,\dots,k$,
where $k$ is the multiplicity of $\al$ in the highest root. (For classical groups, $k\le 2.$)
Let 
\begin{equation}\label{eq:3.3.0}
\fk n=\oplus_{i=1}^k \fk n_i
\end{equation}
 be the corresponding decomposition
into eigenspaces, and decompose each $\fn_i$ into simple $\fk{sl}(2)$
modules. 

The following statement follows from the geometric classification (and
proposition \ref{p:2.3}), and it is also known as a consequence of the main result
of \cite{MS}.

\begin{theorem}[\cite{MS}]\label{t:3.3} Assume $\sigma$ is a generic
  tempered module. Let $\fk n=\oplus_j (d_{ij})$ be the
  decomposition of $\fk n_i$, $i=1,k,$ into simple $\fk{sl}(2)=\bC\langle
  e,h,f\rangle$ modules, where $(d)$ denotes the simple module of
  dimension $d$. Then the reducibility points of $X(M,\sigma,\nu)$,
  with $\nu>0$, are $$\left\{\frac {d_{ij}+1}{2i}\right\}_{i,j}.$$
\end{theorem}

Now we restrict to the case when $\sig$ is a generic discrete series,
and set $\CO:=G\cdot e.$ Moreover, since
$\fk m=\fk m_{BC}$ is a maximal Levi component, the algebra $\fz(\CO)$ is
either $sl(2)$ or a one dimensional torus (\cite{Ca}). If the trivial $\fk{sl}(2)$ module appears in the
decomposition (\ref{eq:3.3.0}), let $i(\sigma)$ denote the eigenvalue $i$ for which it
appears. This is the case precisely when $\fz(\CO)=A_1.$ It turns out that $i(\sigma)\in\{1,2\}.$

\begin{proposition}\label{p:3.3}
\begin{enumerate}
\item If $\fz(\CO)=T_1$ (\ie it is a one-dimensional torus),
  then $X(M,\sigma,\nu)$ is reducible at $\nu=0.$
\item If $\fz(\CO)=A_1,$ then $X(M,\sigma,\nu)$ is
  irreducible at $\nu=0.$
\item When $\fz (\CO)=A_1,$ and $\CO\ne A_4+A_2+A_1$ in $E_8$,  the first
  reducibility point of $X(M,\sigma,\nu),\ \nu\ge 0$ is  
\begin{equation}
  \label{eq:3.3.1}
\nu_0=\frac1{i(\sigma)}.
\end{equation}
\item When $\CO=A_4+A_2+A_1$ in $E_8,$ the first reducibility point of
  $X(M,\sigma,\nu),$ $\nu\ge 0,$ is $\nu_0=3/10$ (while $1/i(\sigma)=1/2)$.
\end{enumerate}  
\end{proposition}

\begin{proof}
This follows from  the conditions in
proposition \ref{p:2.3}. Alternatively, for the reducibility points
$\nu>0$,  one can use theorem \ref{t:3.3}
which has a different proof. When $\CO=A_4+A_2+A_1,$  we have $k=6.$
The trivial $\fk{sl}(2)$-module appears in $\fk n_2,$ so
$\nu_0=1/i(\sigma)=1/2.$ But in this case, $X(M,\sig,\nu)$ is
reducible at $\frac 3{10}$ because $\CO'=A_4+A_3$ is as in (1) of
proposition \ref{p:2.3}. Equivalently, because there is 
a $2$-dimensional $sl(2)$-module in $\fn_5,$ theorem \ref{t:3.3} gives
a reducibility point $3/10.$
\end{proof}

\subsection{}\label{sec:3.1}

Assume  $(M,\sigma,\nu)$ is hermitian with $\sig$ a generic discrete
series,  and let $w\in W$ be such that $w(M)=M,\ w\sig\cong\sig,\ w\nu=-\nu.$
Recall that $\mu_M(\CO,triv)$ is the lowest $W(M)$-type of $\sig,$ and
$\mu(\CO,triv)$ is the generic lowest $W-$type of $X(M,\sigma,\nu).$
As in section \ref{sec:2}, the element $w$ gives rise to intertwining
operators $\CA_{w,\mu}(\sigma,\nu)$ on each $W-$type $\mu$ appearing in
$X(M,\sigma,\nu).$ Recall that these operators are normalized so that
$\CA_{w,\mu(\CO,triv)}(\sigma,\nu)$ is the identity operator.
  
The following result  is proposition 2.4 in \cite{BM3}. For $\nu>>0,$  $X(M,\sig,\nu)$ is irreducible, so the signature on any
$W-$type is constant. We call this the {\it signature at $\infty.$}

\begin{lemma}\label{l:3.1} Assume the W-type $\mu$
satisfies the conditions:
\begin{align}\notag
&\dim\Hom_W[\mu:X(M,\sigma,\nu)]=1\ \text{and}\\\notag
&\Hom_{W(M)}[\mu:\sigma]=\Hom_{W(M)}[\mu:\mu_M(\CO,triv)].
\end{align}
Then the signature at $\infty$ of the operator $\CA_{w,\mu}(\sigma,\nu)$
is
\begin{equation*}
  d(\mu)=(-1)^{\deg\mu+\deg\mu(\CO,triv)},
\end{equation*}
where $\deg\mu$ denotes the lowest harmonic degree of $\mu.$
\end{lemma}

Now we turn to the unitarity of $X(M,\sigma,\nu).$ 

\begin{proposition}\label{p:3.1} Let $(M,\sigma,\nu)$, $\nu>0,$ be hermitian
  maximal parabolic data attached to a
  nilpotent orbit $\CO$, with $\sigma$ a generic discrete series, and $e\in\CO.$  

\begin{itemize}
\item[(1)] Assume that $\fz(\CO)=T_1.$ Then there exists a lowest
  $W-$type $\mu(\CO,\psi),$ $\psi\neq triv$, of $X(M,\sigma,\nu)$, occurring with multiplicity
  1, such that 
\begin{equation}\label{eq:3.1.1}
\CA_{w,\mu(\CO,triv)}(\sigma,\nu)=+Id
 \text{ and } \CA_{w,\mu(\CO,\psi)}(\sigma,\nu)=-Id, \text{ for }\nu>0.
\end{equation} 

\item[(2)] Assume that $\fz(\CO)=A_1.$ Let $\CO'$ be the nilpotent
  orbit in $\fg$ which meets $\fm\times\fz(\CO)$ in $e$ on $\fm$ and
  the principal orbit on $\fz(\CO).$  
  Then $\mu(\CO',triv)$ occurs
  with multiplicity 1 in $X(M,\sigma,\nu)$ , and
\begin{equation}\label{eq:3.1.2}
\CA_{w,\mu(\CO,triv)}(\sigma,\nu)=Id\text{ and
}\CA_{w,\mu(\CO',triv)}(\sigma,\nu)=\left(\frac 
  {1-{i(\sigma)}\nu}{ 1+{i(\sigma)}\nu}\right)^\ell Id,\text{ for }\nu\ge 0,
\end{equation} 
where $\ell$ is some odd positive integer (which may depend on $(M,\sigma)$).
\end{itemize}

\end{proposition}
For uniformity, in case (1) of the proposition, or if $(M,\sigma,\nu)$
is never hermitian for $\nu>0$, set
$\CO'=\CO.$ (This notation will be used in section \ref{sec:4.1a}.)

\begin{proof} We give an outline of the argument. Complete details for
  type $E_8$ are presented in section \ref{sec:maxpar}.

 If $\fz(\CO)=T_1$, and $(M,\sigma,\nu)$, $\nu>0$,
 is hermitian, then $A(\CO)\neq 1$ (\cite{Ca}). 
 The standard module $X(M,\sigma,\nu)$ has two lowest $W-$types
$\mu(\CO,triv)$ and $\mu(\CO,\psi)$ both appearing with multiplicity 
one and having lowest harmonic degrees of opposite parity. At $\nu=0$, 
$X(M,\sigma,0)$ is reducible and each factor is a tempered module,
therefore unitary. If $\nu>0$, $\mu(\CO,triv)$ and $\mu(\CO,\psi)$
always occur in $X(M,\sigma,\nu)$. Having opposite signature at
$\infty$, they have opposite signature for all $\nu>0.$

\medskip
 
If $\fz(\CO)=A_1$, then  $\ovl X(M,\sigma,\nu)$ has a
unique lowest $W-$type $\mu(\CO,triv)$ (\cite{Ca}). The module $X(M,\sigma,0)$ is irreducible
and tempered.  At $\nu=\frac 1{i(\sigma)}$, all factors other than
$\ovl X(M,\sigma,\nu)$ are parametrized by strictly
larger nilpotent orbits. One of the factors
corresponds to the orbit $\CO'$  and
 lowest W-type $\mu(\CO',triv).$ We verify in every case that $\mu(\CO',triv)$ satisfies the conditions of lemma
\ref{l:3.1}. Moreover $\mu(\CO',triv)$ has harmonic degree of
opposite parity to $\mu(\CO,triv)$. The claim follows then from the fact that for
$\nu>\nu_0$, the two $W-$types $\mu(\CO,triv)$ and $\mu(\CO',triv)$ occur in the
factor $\ovl X(M,\sigma,\nu)$.

\end{proof}

In summary,
\begin{corollary}\label{c:3.1}  
\begin{itemize}
\item[(1)] If $\fz(\CO)$ is of type $T_1$, 
then $\ovl X(M,\sigma,\nu)$ is not unitary  for $\nu>0$.\notag

\item[(2)] If  $\fz(\CO)$ is of type
  $A_1$, then $X(M,\sigma,\nu)$ is unitary and irreducible if
  and only if $0\le\nu<\nu_0$, where $\nu_0$ is the first
  reducibility point of $X(M,\sigma,\nu)$ on the half-line
  $\nu>0$.\notag

\end{itemize}
\end{corollary}

\begin{proof}
Part (1) follows directly from proposition \ref{p:3.1}. For part (2),
we also immediately have that $X(M,\sigma,\nu)$ can only be unitary in
the interval $[0,\frac 1{i(\sigma)}).$ Since  $X(M,\sigma,\nu)$ is irreducible
and unitary at $\nu=0,$ it stays unitary until the first point of
reducibility $\nu_0.$ When $\CO\neq A_4+A_2+A_1,$ we have
$\nu_0=\frac 1{i(\sigma)}$ (see proposition \ref{p:3.3}), so this
completes the argument.  
For $\CO=A_4+A_2+A_1$ in $E_8,$ we need an extra argument to rule out
the segment $(\nu_0,\frac 1{i(\sigma)})=(\frac 3{10},\frac 12).$ The
details of this case appear in section \ref{sec:A4A2A1}.
\end{proof}

{
\begin{remark}Notice that equation (\ref{eq:3.1.2}) only tells us
  that, at the  reducibility point $\nu=1/i(\sigma)$,
  the order of the zero for the operator
  $\CA_{w,\mu(\CO',triv)}(\sigma,\nu)$ is an odd integer $\ell.$ This
  is of course sufficient to conclude that corollary \ref{c:3.1}
  holds. However, it is natural to conjecture that $\ell=1$ for all $(M,\sigma),$ where $M$ is a Levi of a maximal
  parabolic, and $\sigma$ is generic. We
verified this conjecture by computing
$\CA_{w,\mu(\CO',triv)}(\sigma,\nu)$ explicitly in all cases
$(M,\sigma)$ as above, when $G$ is simply-laced of rank at most $7.$  
\end{remark}
}

\subsection{}\label{sec:4.1a} Fix a nilpotent orbit $\CO$, and, as
before,  let
$\sigma$ be the generic discrete series of $\bH_{M_0}$ parametrized by
$\CO.$  

\begin{definition}\label{d:4.1a} If
  $\fm_{BC}$ is maximal Levi subalgebra, recall the orbit $\CO'$ constructed 
  in proposition \ref{p:3.1}. We say that the $W-$type $\mu$ is
  {\it $\sigma$-petite} if $\mu$ is a lowest $W-$type for $\CO$ or for
  $\CO'.$   

\smallskip

If $\fm_{BC}$ is not maximal, let $\fm_1,\dots,\fm_k$ be all the Levi
subalgebras, not necessarily of a standard parabolic subalgebra, such
that $\fm_{BC}\subset\fm_j$, and $\fm_{BC}$ is a maximal Levi
subalgebra of $\fm_j$, $j=1,\dots ,k.$ For every $j$, let $\{\mu_{ij}\}_i$
denote the set of $W(M_j)$-types which are $\sigma$-petite in
$\fm_j$. We say that the $W-$type $\mu$ is {\it $\sigma$-petite} (in
$\fg$) if for every $j$, the only 
$W(M_j)$-types of $\sigma$ contained in the restriction
$\mu|_{W(M_j)}$ are the petite $W(M_j)$-types $\mu_{ij}.$  

\end{definition}

Clearly every lowest $W-$type of $X(M,\sigma,\nu)$ which contains
$\mu_M(\CO,triv)$ in its restriction to $W(M)$ is
$\sigma$-petite. 

\medskip

\noindent{\bf Example.} For the spherical principal series, \ie
$\sigma=triv$, $M=H,$ this definition is a tautology: every $W-$type is 
$\sigma$-petite. The other extremal case is when $M$ is maximal
parabolic; then there are exactly two $\sigma$-petite $W-$types, those
from proposition \ref{p:3.1}.

An intermediate example (for $\bH$ of simply-laced type) is when $\CO=A_1$, the
minimal nilpotent orbit. Then $\sigma=St$ for $M=A_1,$ and a $W-$type
$\mu$ is $\sigma$-petite if and only if $\mu|_{W(A_2)}$ does not
contain the sign representation.

\subsection{}\label{sec:4.2a}
The following lemma should be compared with corollary \ref{c:2.8b}.

\begin{lemma}\label{l:4.2a} 
\begin{enumerate}
\item The Langlands parameter $(M,\sigma,\nu)$ is hermitian if and
  only if $\nu$ is a hermitian (spherical) parameter for $\bH(Z(\CO)).$
\item If $(M,\sigma,\nu)$ is hermitian, but either $\nu$ is not a hermitian
  (spherical) parameter for $\bH(\fz(\CO))$ or $\nu$ is not in the
  semisimple part of $\fz(\CO)$, then $\overline
  X(M,\sigma,\nu)$ is not unitary.
\end{enumerate}
\end{lemma}

\begin{proof} We verify these assertions in section \ref{sec:5}. For
  part (2), the method is the same as in proposition \ref{p:3.1}(1):
  we find two lowest $W-$types $\mu(\CO,triv)$ and $\mu(\CO,\psi),$
  $\psi\neq triv$ of $\ovl
  X(M,\sigma,\nu),$ occurring with multiplicity one, such that the operators
  $\CA_{\mu(\CO,triv)}(\sigma,\nu)$ and $\CA_{\mu(\CO,\psi)}(\sigma,\nu)$
  have opposite   signatures. 
\end{proof}

\subsection{}\label{sec:4.2}
The main result, theorem \ref{t:4.1}, is a consequence of
the construction in this section, which also provides an explanation of why
such a result should hold.  The method of calculation is uniform, but
the details need to be checked in each case. (In sections
\ref{sec:5} and \ref{exceptions}, we will only present the detailed
calculations in type $E_8.$)  
To help orient the  reader, we give an outline of the method.

Recall that $X(M,\sigma,\nu)$ is an induced module, where $\sigma$ is
a generic discrete series parameterized by a Lie triple
$\{e,h,f\}\subset\fm.$
Also from section \ref{sec:2.10},
recall that $\fa$ denotes a Cartan subalgebra of $\fz(\CO)$ with $\nu\in \fa$
and $C(\fa,M)\subset W$ is defined by (\ref{eq:2.9.1}). For simplicity
we drop here the subscript $BC.$ By proposition
\ref{p:2.10}, $C(\fa,M)$ is the
image of a homomorphism of $W(Z(\CO))$ to $W$.  If $w$ is an element of $W(Z(\CO))$, we will denote by
$\ovl w$ its image in $W$ under this homomorphism.

 By lemma \ref{l:4.2a}, we may assume that
$(M,\sigma,\nu)$ is hermitian and that $\nu$ is hermitian (spherical)
for $\bH(\fz(\CO))$ and in the semisimple part of $\fz(\CO).$ This
means that there exists 
\begin{equation}
\label{eq:4.2.0}
w_Z\in W(\fz(\CO)) \text{ such that }w_Z\nu=-\nu. 
\end{equation}
Let $\CA_{\ovl w_Z}(\sigma,\nu)$ be the $\bH$-intertwining operator
(see section \ref{sec:disc}) which induces
the operators $\CA_{\ovl w_Z,\mu}(\sigma,\nu)$, $\mu\in \widehat W.$

The element $w_Z\in \fz(\CO)$ defines a spherical $\bH(\fz(\CO))$-intertwining
operators (equation (\ref{eq:2.11.2a})) $a_{\rho(\mu)}(\nu),$
$\rho(\mu)\in \widehat{W(\fz(\CO))}.$ We would like to show that for
$\mu$ a $\sigma$-petite $W-$type, these two operators defined by
$w_Z$, actually {\it coincide}.

The idea is to decompose
$\CA_{\ovl w_Z,\mu}(\sigma,\nu)$ into a product of factors similar to the usual
decomposition of the spherical long intertwining operator (as in
section \ref{sec:2.11}) for
$\bH(\fz(\CO))$, such that  each factor in
$\CA_{\ovl w_Z,\mu}(\sigma,\nu)$ is identical to the corresponding simple factor in the
spherical intertwining operator of  $\bH(\fz(\CO))$. 

\smallskip 

For each simple root $\bar\al\in \Pi(\fz(\CO),\fa)$, we find an element $\bar
s_\alpha\in C(\fa,M)$, which induces the corresponding simple reflection
on $\fa$. Then the $\bar s_\al$'s generate a subgroup of $C(\fa,M)$
isomorphic to $W(\fz(\CO))$. Let  $\overline w_Z$ be the
image in $C(\fa,M)$ of $w_Z.$ 

We apply the construction in section
\ref{sec:disc}. First, the operators $\CA_{\ovl w_Z,\mu}(\sigma,\nu)$ decompose into a product
of the form 
\begin{equation}\label{4.2.1}
\CA_{\bar s_{\al_1},\mu}(\sigma,\nu_{\bar\al_1})\dotsb \CA_{\bar
  s_{\al_k},\mu}(\sigma,\nu_{\bar\al_k}), 
\end{equation}
corresponding to a decomposition $\ovl w_Z=\bar s_{\al_1}\dotsb \bar
s_{\al_k}$.  

Fix an $\bar\al.$  The reflection $\bar s_\al$ preserves
$(M,\sigma)$. By lemma \ref{l:factoring} there exists a chain of
adjacent Levi components $\fm=\fm_0,\dots,\fm_k=\fm,$ such that $\bar
s_\al$ can be decomposed into a product 
\begin{equation}\label{eq:4.2.1a}
\bar s_\al=w_{k}\dotsb w_{1},
\end{equation}
as in equation (\ref{eq:2.4.8}). 

The operator $\CA_{\bar
  s_\al,\mu}(\sigma,\nu_{\bar\al})$ acquires a decomposition accordingly
into a
product of maximal parabolic factors of the form
$\CA_{w_{m_j},\mu}(\wti w_{j}\sigma,\wti w_{j}\nu_{\bar\al}),$ where
  $\wti w_j=w_{m_{j+1}}\dots w_{m_1}.$   

\medskip

Recall from definition \ref{d:2.10b} that the space
$\Hom_{W(M)}[\mu,\sigma]$ has a natural structure of
a $W(\fz(\CO))$-type, which is denoted $\rho(\mu),$ and a structure
of $W(Z(\CO)$-type, which is denoted $\rho'(\mu).$ 

\begin{lemma}\label{l:4.2} With the notation above, if $\mu$ is a
  $\sigma$-petite $W-$type (definition \ref{sec:4.1a}), and $\bar\al$
  is a simple root of $\fz(\CO),$ then $$\CA_{\bar s_\al,\mu}(\nu)=a_{\rho(\mu),\bar\al}(\nu),$$
where $a_{\rho(\mu),\bar\al}(\nu)$ is given by equation (\ref{2.12.3}).
\end{lemma}

\begin{proof}
In the discussion above, we have decomposed $\CA_{\bar
  s_\al,\mu}(\sigma,\nu_{\bar\al})$ into a product of factors,
$\CA_{w_{m_j},\mu}(\wti w_j\sigma,\wti w_j\nu_{\bar\al}),$ each
  induced from some maximal parabolic case $\fk m_j\subset \Sigma_j.$ As such, for every $j$, the
  discrete series $\wti w_j\sigma$ is parametrized in $\bH_{M_{j,0}}$
  by a nilpotent element whose reductive centralizer $\fz_{\Sigma_j}$ in $\Sigma_j$ is
  either an $sl(2)$ or a one-dimensional torus. 

By inspection, in section \ref{sec:5}, we find that in the
decomposition induced by (\ref{eq:4.2.1a}), there exists $j_0$ such that
$\fz_{\Sigma_{j_0}}=sl(2)$, and if $j\neq j_0,$ then $\fz_{\Sigma_j}$
is a torus. By the definition of $\sigma$-petite in the maximal
parabolic case, and proposition \ref{p:3.1}, the factors
$j\neq j_0$ do not contribute, while the factor $j=j_0$ is identical
with $a_{\rho(\mu),\bar\al}(\nu) $.   
\end{proof}

We summarize the construction in the following proposition. Retain the
previous notation, and let $X_{\bH(\fz(\CO))}(\nu)$ denote the
spherical principal series for $\bH(\fz(\CO)).$ 

\begin{proposition}\label{p:4.2}  
Assume $(M,\sigma,\nu)$ (where $M=M_{BC})$ is hermitian with $\nu$
hermitian for $\bH(\fz(\CO))$ and $w_Z$ as in (\ref{eq:4.2.0}). 
If $\mu$ is a $\sigma$-petite $W-$type (definition \ref{sec:4.1a}), let $\rho(\mu)$ be the
corresponding $W(\fz(\CO))$-type (section \ref{sec:2.10b}). 

The
$\bH$-intertwining operator $\CA_{\ovl w_Z,\mu}(\sigma,\nu)$ on the space
$$\Hom_W[\mu:X(M,\sigma,\nu)]$$ coincides with the spherical
$\bH(\fz(\CO))$-intertwining operator 
$a_{\rho(\mu)}(\nu)$ on the space
$$\Hom_{W(\fz(\CO))}[\rho(\mu):X_{\bH(\fz(\CO))}(\nu)]=\rho(\mu)^*.$$

In this matching, the generic lowest $W-$type $\mu(\CO,triv)$ corresponds to the
trivial $W(\fz(\CO))$-type.
\end{proposition}

\subsection{}\label{sec:4.3} Assume we are in the setting of
proposition \ref{p:4.2}. If the parameter $\nu$ is such that
$A_G(e,\nu)\neq A_M(e,\nu),$ then the image of the intertwining operator
$\CA_{\ovl w_Z}(\sigma,\nu)$ is {\it not} irreducible. In this case,
by proposition \ref{p:2.10}, (2) (see also the remark after
(\ref{eq:2.10.7})), we have a decomposition under the action of $A_G(e,\nu)$
\begin{equation}\label{eq:4.3.1}
X(M,\sigma,\nu)=\bigoplus_{(\psi,V_\psi)\in \widehat{A_G(e,\nu)}}
X(M,\sigma,\nu,\psi)\otimes V_\psi,
\end{equation}
which induces
\begin{equation}\label{eq:4.3.2}
\Hom_W[\mu:X(M,\sigma,\nu)]=\bigoplus_{(\psi,V_\psi)\in
  \widehat{A_G(e,\nu)}} \Hom_W[\mu:X(M,\sigma,\nu,\psi)]\otimes V_\psi.
\end{equation}

Recall that the intertwining operators are normalized so that the operator
on the generic lowest $W-$type is identically $1$. Then, as in section
\ref{sec:disc}, $\CA_{\ovl w_z}(\sigma,\nu)$ induces operators
\begin{equation}\label{eq:4.3.1a}
\CA_{\ovl w_Z,\mu}(\sigma,\nu,triv):\Hom_W[\mu:X(M,\sigma,\nu,triv)]\to \Hom_W[\mu:X(M,\sigma,-\nu,triv)].
\end{equation}

Recall that in section \ref{sec:extended} we constructed the spherical
principal series $X'_{\bH(Z(\CO))}(\nu)$ (equation (\ref{eq:2.8b.3})) for the extended Hecke
algebra $\bH(Z(\CO)),$ as well as the operators $a'_{\rho(\mu')}(\nu)$
(equation (\ref{eq:2.8.9})). 

\begin{corollary}\label{p:4.3}Retain the notation from proposition
  \ref{p:4.2} and equation (\ref{eq:4.3.1a}). 

The $\bH$-intertwining operator $\CA_{\ovl w_Z,\mu}(\sigma,\nu,triv)$ on
$$\Hom_W[\mu:X(M,\sigma,\nu,triv)]$$ is identical with the
$\bH(Z(\CO))$-intertwining operator $a'_{\rho'(\mu)}(\nu)$ (defined in
\ref{eq:2.8.9}) on the space
$$\Hom_{W(Z(\CO))}[\rho'(\mu):X'_{\bH(Z(\CO))}(\nu)],$$ 
which in turn is equivalent with the $\bH(\fz(\CO))$-intertwining operator $a_{\rho(\mu)}(\nu)$ restricted to
the subspace
$$((\rho'(\mu)^*)^{A_G(e,\nu)}.$$     
\end{corollary}

\begin{proof} Follows from proposition \ref{p:4.2} and section \ref{sec:2.8c}.

\end{proof}

\subsection{}\label{sec:4.4}
Fix $\CO$ a nilpotent orbit in $\fg$, and let $M=M_{BC}$, $\{e,h,f\},$ and $\sigma$
be as before. Let $\fk S(\CO)$
denote the set of $\sigma$-petite $W-$types (definition
\ref{d:4.1a}). Set
\begin{equation}\label{4.4.1}
\rho(\fk S(\CO))=\{\rho(\mu)\in \widehat{W(\fz(\CO))}: \mu\in \fk S(\CO)\},
\end{equation}
where $\rho(\mu)$ is defined in \ref{d:2.10b}.

 By comparison with the spherical
intertwining operators in $\bH(\fz(\CO))$, the matching of intertwining operators in
 sections  \ref{p:4.2} and \ref{p:4.3} tells us the signature of
the hermitian form on the $\sigma$-petite $W-$types. 

By section \ref{sec:relevant},  one knows a very small subset of $W\widehat{(\fz(\CO))}$, the
$0$-relevant $W(\fz(\CO))$-types (definition \ref{d:relevant}), which are
sufficient to detect the unitarity of the $0$-complementary
series. Call this set $\fk B(\fz(\CO)).$

\begin{definition}\label{d:4.4}
We say that $\CO$ satisfies the {\it signature criterion} if $\fk
B(\fz(\CO))\subset \rho(\fk S(\CO)).$
\end{definition}

Our main criterion of nonunitarity follows from this discussion.

\begin{corollary}\label{c:4.4}
\begin{enumerate}
\item If $\CO$ satisfies the signature criterion, then necessarily a parameter
$\chi=h/2+\nu$ is in the $\CO$-complementary series of $\bH$ {\it only
  if} $\nu$ is the $0$-complementary series for $\bH(\fz(\CO)).$
\item If $\bH$ is of type $E,$ the only nilpotent orbits which do {\it
    not} satisfy the signature criterion are $4A_1$ in $E_7,$
and $D_4+A_1$, $2A_2+2A_1$, $4A_1$ in $E_8.$
\end{enumerate}
\end{corollary}

\begin{proof}
Part (1) is clear. Part (2) is established by computing the
$\sigma$-petite $W-$types. The calculations for type $E_8$ are in
section 
\ref{sec:5}. 
\end{proof}

Note that the nilpotent $4A_1$ in
$E_8$ is one of the exceptions in theorem \ref{t:4.1}, and in fact the
complementary series turns out to be larger than the $0$-complementary
series for the centralizer $\fz(\CO)=C_4.$ 

For the other cases, $4A_1$ in $E_7$, $D_4+A_1$,
$2A_2+2A_1$ in $E_8,$ we use ad-hoc additional arguments involving the signature of some other $W-$types which appears with
small multiplicity (and by Springer's correspondence belong to
nilpotent orbits close to $\CO$ in the closure ordering), to prove the
inclusion of the $\CO$-complementary series of $\bH$ into the
$0$-complementary series of $\bH(\fz(\CO)).$ 

\subsection{}\label{sec:4.5} {\bf Case $\mathbf{\CO\neq 4A_1}$ in
  $\mathbf{E_8}$}. Let us assume that, by
the previous discussion, we know that the complementary series of $\CO$
is included in the $0$-complementary series of $\bH(\fz(\CO)).$

 Using the method of decomposing intertwining
operators into factors coming from maximal parabolic cases (section
\ref{sec:disc}), and the reducibility points for maximal parabolic
cases (section \ref{sec:3.3}), we can determine the hyperplanes of
reducibility of standard modules $X(M,\sigma,\nu,triv)$.

We check if any of these hyperplanes of reducibility intersects the
$0$-complementary series of $\bH(\fz(\CO)).$
When this happens, we are in one of the exceptions of theorem
\ref{t:4.1}. In these cases, we need some extra arguments involving
the signature of operators on $W-$types which are not $\sigma$-petite,
but they rule out the nonunitary parameters
$\chi=h/2+\nu$, with $\nu$ inside the  $0$-complementary series of
$\bH(\fz(\CO)).$ 
The  details are in sections \ref{sec:A4A2A1}-\ref{sec:A2+2A1}.

\smallskip

 {We consider the cases when the reducibility
hyperplanes do not intersect the $0-$complementary series; in this
case we need to show that the parameters in the $0-$ complementary
series for $\bH(\fz(\CO))$ are unitary for $\bH$.  Every
parameter $\chi=\frac h2+\nu$ in this set can be deformed continuously and
irreducibly to a parameter $\chi_0=\frac h2+\nu_0,$ for which the
corresponding standard module is unitarily and irreducibly induced
from a unitary module on a Levi subgroup. The unitarity follows from
the following well known result.

\begin{lemma}\label{l:4.5}
For $0\le t\le 1,$ let $\xi_t\in \fz(\fm)$ be a family of characters which
depend continuously on $t$, and $\xi_0$ is unitary. Assume that $\Ind_M^G[\C V\otimes\xi_t]$ is
irreducible, where $\C V$ is a module for $\bH_M.$ If $\C V\otimes \xi_t$ is
hermitian for all $0\le t\le 1,$ then $\Ind_M^G[\C V\otimes\xi_1]$ is
unitary if and only if $\C V$ is unitary.
\end{lemma}
}

\medskip
\noindent{\bf Case $\mathbf{\CO=4A_1}$ in $\mathbf{E_8}$}. Here  $\fz(\CO)=C_4,$ $M=4A_1$ and
$\sigma=St$. The  details are in
section \ref{sec:4A1}. 
Using the signature of the $\sigma$-petite $W-$types, we
find that the $4A_1$-complementary series is formed of parameters
$\chi=h/2+\nu$, where $\nu$ must lie in one of two regions. 

The first region corresponds to $\nu$ in the $0$-complementary series of the
$\bH(\fz(\CO))$, and we can show that $\chi$ is
unitary by the same deformation argument as in lemma \ref{l:4.5}.

 If $\chi$
is in the second region, called $\C R$ in section \ref{sec:4A1},  a more delicate argument is needed. First we
analyze the signature of other $W-$types, which are not
$\sigma$-petite, and find that there exists only one possible unitary
subregion $\C R_3$ of $\C R.$ (The notation and explicit description
are in (\ref{eq:4A1}).)  Now assume $\nu\in \C R_3.$ We deform $\nu$
continuously to $\nu_0,$ such that $X(4A_1,St,\nu)$ is irreducible for
$\nu\neq \nu_0,$ but $X(4A_1,St,\nu_0)$ is reducible. We find that
$X(4A_1,St,\nu_0)$ has two composition factors, and that they are both
unitary. Then we use a signature filtration type of argument to
conclude that $X(4A_1,St,\nu)$ must be unitary.

\section{Explicit calculations for type $E_8$}\label{sec:5}

The simple roots $\al_i$ and coweights $\check\om_i$, $i=1,8$ in type
$E_8$ are as in \cite{Bou}.
 The $W-$types for  $E_8$ were classified in \cite{Fr}, and we will
  use  the same labeling of  the irreducible
  characters. (See also \cite{Ca}.)
  The $W-$structure of
  standard modules is given by the Green polynomials calculated in
  \cite{BS}; we also used the (unpublished) tables in
  \cite{Al}. 
For restrictions of $W-$types and for the computation
  of the associated $W(\fz(\CO))$-type $\rho(\mu)$ to a given W-type
  $\mu$ (notation as in \ref{sec:4.2}), we used the software
  ``GAP''. 
For some of the explicit
  computations with intertwining operators in the maximal parabolic
  cases for exceptional groups (see the remark after proposition \ref{p:3.1}),
 we used 
  integer matrix models of $W-$types, 
  and the software ``Mathematica''. 
The classification and labeling of nilpotent orbits is as in \cite{Ca}.

\subsection{} If a nilpotent orbit is
distinguished, it parametrizes discrete series, and in particular,
exactly one generic discrete series. The corresponding infinitesimal
characters are in the tables of section \ref{tables}.

For the explicit calculations of intertwining operators that we need
(see remark \ref{sec:3.1} for example),
when the standard module is not induced from a Steinberg
representation on a Levi subalgebra, we embed it into an induced from
the Steinberg representation from a smaller subalgebra, such that the
generic lowest W-type appears with multiplicity one. This is possible
because the rank is small. Below is the table of embeddings for
discrete series. We give the distinguished non-principal nilpotent orbit $\CO$, the
lowest W-type $\mu_0$ corresponding to the trivial representation
in $\widehat {A(\CO)}$, and a Levi component $M$ such that $\dim\Hom_W[\mu_0:\text{Ind}_M^G(St)]=1.$

\begin{tiny}
\begin{center}
\begin{longtable}{|c|c|c|c|}
\caption{Embeddings of discrete series}
\label{table:embedding}\\
\hline
\multicolumn{1}{|c|}{\bf Type } &\multicolumn{1}{|c|}{\bf Nilpotent } &\multicolumn{1}{c|}{\bf Lowest
  W-type} &\multicolumn{1}{c|}{\bf Levi component}\\\hline
\endfirsthead
\multicolumn{4}{c}%
{{ \tablename\ \thetable{} -- continued from previous page}}
\\
\endhead
\hline
\endlastfoot

$D_4$ &$(5,3)$ &$1^3\times 1$ &$D_3$\\ 
\hline
$D_5$ &$(7,3)$ &$1^4\times 1$ &$D_4$\\
\hline
$D_6$ &$(9,3)$ &$1^5\times 1$ &$D_5$\\
      &$(7,5)$ &$1^4\times 1^2$ &$A_5$\\
\hline
$D_7$ &$(11,3)$ &$1^6\times 1$ &$D_6$\\
      &$(9,5)$ &$1^5\times 1^2$ &$A_6$\\
\hline
$E_6$ &$E_6(a_1)$ &$6_p'$ &$D_5$\\
      &$E_6(a_3)$ &$30_p'$ &$A_5$\\
\hline
$E_7$ &$E_7(a_1)$ &$7_a$ &$E_6$\\
      &$E_7(a_2)$ &$27_a'$ &$E_6$\\
      &$E_7(a_3)$ &$56_a$ &$D_6$\\
      &$E_7(a_4)$ &$189_b$ &$D_5+A_1$\\
      &$E_7(a_5)$ &$315_a$ &$A_5+A_1$

\end{longtable}
\end{center}
\end{tiny}

\subsection{}\label{sec:maxpar}
 For the maximal parabolic cases, we verify all the details of the argument
outlined in the proof of proposition \ref{p:3.1}. Depending on the
details of the discussion, there are three types of arguments that we
consider. For each type, we present the details in one example, then
list the other nilpotents for which the same argument applies. The
only exception is the nilpotent $A_4+A_2+A_1$, which we treat separately.

To simplify notation, we will denote by $\mu_0,$ $\mu_0',\dots,$  the
lowest $W-$types $\mu(\CO,triv),$ $\mu(\CO,\psi),\dots$, and by
$\mu_1,\mu_2,\dots$, the $W-$types of the form $\mu(\CO',triv).$

\subsubsection{$\mathbf {E_7}$} The centralizer is $\fz(\CO)=A_1$, the
lowest W-type is $\mu_0=84_x',$ and the infinitesimal character is
$\chi=(0,1,2,3,4,5,-17/2,17/2)+\nu \check\omega_8,$ with $\nu\ge 0.$  

The standard module corresponding to $\CO=E_7$ is $X(E_7,St,\nu).$ The first
reducibility point is at $\nu_0=\frac 12,$ where the generic factor is
parametrized by the nilpotent orbit $\CO'=E_8(a_3)$ and lowest W-type
$\mu_1=112_z'$. The $W-$types $\mu_0$ and $\mu_1$ have opposite signs at
$\infty$. Since the nilpotent $\CO'$ is distinguished, there cannot be
another factor with lowest W-type $\mu_1$ for $\nu>\nu_0.$
Therefore $\mu_0$ and $\mu_1$ stay in the same factor for
$\nu>\nu_0.$ 
The {complementary series} is $0\le\nu<\frac 12$.

\begin{tiny}
\begin{center}
\begin{longtable}{|c|c|c|c|c|}\caption {Maximal parabolic cases,
    $\fz(\CO)=A_1$, type 1}\\
\hline
\multicolumn{1}{|c|}{$\mathbf{\CO}$} &\multicolumn{1}{c|}{$\chi$}  &\multicolumn{1}{c|}{$\mathbf{\mu_0}$} &\multicolumn{1}{c|}{$\mathbf{\CO'}$}&\multicolumn{1}{c|}{$\mathbf{\mu_1}$}\\
\hline
\endfirsthead
\multicolumn{5}{c}%
{{  \tablename\ \thetable{} -- continued from previous page}}
\\
\hline
\multicolumn{1}{|c|}{$\mathbf{\CO}$} &\multicolumn{1}{c|}{$\chi$}  &\multicolumn{1}{c|}{$\mathbf{\mu_0}$} &\multicolumn{1}{c|}{$\mathbf{\CO'}$}
&\multicolumn{1}{c|}{$\mathbf{\mu_1}$}\\
\hline
\endhead

$E_7$ &$(0,1,2,3,4,5,-\frac{17}2,\frac{17}2)+\nu \check\omega_8$  &$84_x'$
&$E_8(a_3)$ &$112_z'$\\
\hline
$E_7(a_1)$ &$(0,1,1,2,3,4,-\frac{13}2,\frac{13}2)+\nu\check\omega_8$ &$567_x'$ &$E_8(b_4)$ &$560_z'$\\
\hline
$D_7$ &$(0,1,2,3,4,5,6,0)+\nu\check\omega_1$   &$400_z'$
&$E_8(a_5)$ &$700_x'$\\
\hline
$E_7(a_2)$ &$(0,1,1,2,2,3,-\frac {11}2,\frac {11}2)+\nu\check\omega_8$
 &$1344_x'$ &$E_8(b_5)$ &$1400_z'$\\
\hline
$A_7$ &$(-\frac
{17}4,-\frac {13}4,-\frac 94,-\frac 54,-\frac 14,\frac 34,\frac
74,\frac 74)+\nu\check\omega_2$  &$1400_{zz}'$ &$E_8(b_6)$
&$2240_x'$\\
\hline
$E_7(a_5)$ &$(0,0,1,1,1,2,-\frac 52,\frac 52)+\nu\check\omega_8$ 
&$7168_w$ &$E_8(a_7)$ &$4480_y$\\

\hline
\end{longtable}
\end{center}
\end{tiny}

\subsubsection{$\mathbf {E_6+A_1}$}\label{sec:E6A1} The centralizer is $\fz(\CO)=A_1$,
the lowest W-type is $\mu_0=448_z'$, and the \ic is
$(0,1,2,3,4,-\frac 92,-\frac 72,4)+\nu\check\omega_7$. 

The standard module is $X(E_6+A_1,St,\nu),$ $\nu\ge 0.$ 
The first reducibility
point is at $\nu_0=\frac 12$, where the generic factor is parametrized
by the nilpotent orbit $E_8(b_5)$ and lowest W-type $1400_z'$. (But the
argument from the nilpotent $E_7$  does not apply here since $448_z'$
and $1400_z'$ have the same signature at $\infty$.) 

At $\nu=\nu_0$ there  may also be
a factor parametrized by the nilpotent orbit $\CO'=E_7(a_2)$ with
lowest W-type $\mu_1=1344_x'$. The $W-$types $\mu_0$ and
$\mu_1$ may only be separate for $\nu=\nu_0$. The reason is that
for $\nu>\nu_0$, any irreducible factor with lowest W-type $\mu_1$ must also contain
the W-type $\mu_2=1008_z'$. But $\mu_2$ does not appear in
$X(E_6+A_1,St,\nu)|_W$ at all. Moreover, since $\mu_0$ and $\mu_1$ have
opposite signs at $\infty$, they must be separate at least once, so
they are separate exactly at $\nu=\nu_0$.
The {complementary series} is $0\le\nu<\frac 12$. 

\begin{tiny}
\begin{center}
\begin{longtable}{|c|c|c|c|c|c|}\caption {Maximal parabolic cases,
    $\fz(\CO)=A_1$, type 2}\\

\hline
\multicolumn{1}{|c|}{$\mathbf{\CO}$} &\multicolumn{1}{c|}{$\chi$}
&\multicolumn{1}{c|}{$\mathbf{\mu_0}$}
&\multicolumn{1}{c|}{$\mathbf{\CO'}$}&\multicolumn{1}{c|}{$\mathbf{\mu_1}$}&\multicolumn{1}{c|}{$\mathbf{\mu_2}$}\\
\hline
\endfirsthead
\multicolumn{6}{c}%
{{  \tablename\ \thetable{} -- continued from previous page}}
\\
\hline
\multicolumn{1}{|c|}{$\mathbf{\CO}$} &\multicolumn{1}{c|}{$\chi$}  &\multicolumn{1}{c|}{$\mathbf{\mu_0}$} &\multicolumn{1}{c|}{$\mathbf{\CO'}$}
&\multicolumn{1}{c|}{$\mathbf{\mu_1}$}&\multicolumn{1}{c|}{$\mathbf{\mu_2}$}\\
\hline
\endhead

$E_6+A_1$ &$(0,1,2,3,4,-\frac 92,-\frac 72,4)+\nu\check\omega_7$ 
&$448_z'$ &$E_7(a_2)$ &$1344_x'$ &$1008_z'$ \\
\hline
$E_7(a_3)$ &$(0,0,1,1,2,3,-\frac 92,\frac 92)+\nu\check\omega_8$ 
&$2268_x'$ &$D_7(a_1)$ &$3240_z'$ &$1050_x'$\\
\hline
$E_7(a_4)$ &$(0,0,1,1,1,2,-\frac 72,\frac 72)+\nu\check\omega_8$ 
&$6075_x'$ &$D_5+A_2$ &$4536_z'$ &$840_x'$\\
\hline
$A_6+A_1$ &$(\frac
{13}4,-\frac 94,-\frac 54,-\frac 14,\frac 34,\frac 74,\frac
{11}4,\frac 14)+\nu\check\omega_3$  &$2835_x'$ &$D_5+A_2$
&$4536_z'$ &$840_x'$\\
\hline
$E_6(a_3)A_1$ &$(0,0,1,1,2,-\frac 52,-\frac 32,2)+\nu\check\omega_7$
 &$3150_y$ &$E_7(a_5)$ &$7168_w$ &$1680_y$\\
\hline
$D_5(a_1)A_2$ &$(0,1,1,2,-\frac
52,-\frac 32,-\frac 12,\frac 32)+\nu\check\omega_6$  &$1344_w$
&$E_6(a_3)A_1$ &$1134_y$ &$448_w$\\
\hline
$A_4+A_3$ &$(0,1,2,-\frac 52,-\frac 32,-\frac 12,\frac
12,1)+\nu\check\omega_5$  &$420_y$ &$D_5(a_1)A_2$ &$1344_w$ &$1134_y$\\
\hline
\end{longtable}
\end{center}
\end{tiny}

\subsubsection{$\mathbf {D_7(a_1)}$} The centralizer is
$\fz(\CO)=T_1$, and the \ic is  $(0,1,1,2,3,4,5,0)+\nu\check\omega_1.$

The standard module $X(D_7,\sigma,\nu)$, where $\sigma$ is the generic
discrete series parametrized by the nilpotent orbit $(11,3)$ in the
Hecke algebra of type $D_7$, is reducible at $\nu=0,$ and it has two
lowest $W-$types for $\nu>0,$ $\mu_0=3240_z'$ and $\mu_0'=1050_x'.$    

At $\nu=0$, $X$ breaks into the sum of tempered modules, each
containing one lowest W-type, which are unitary. 
For $\nu>0,$ $\mu_0$ and $\mu_0'$ stay in the same factor, and
they have opposite signs at $\infty.$ 
There is {no complementary series}. The generic module is unitary only
at $\nu=0$.

\begin{tiny}
\begin{center}
\begin{longtable}{|c|c|c|c|}\caption {Maximal parabolic cases,
    $\fz(\CO)=T_1$}\\

\hline
\multicolumn{1}{|c|}{$\mathbf{\CO}$} &\multicolumn{1}{c|}{$\chi$} &
\multicolumn{1}{c|}{$\mathbf{\mu_0}$} &\multicolumn{1}{c|}{$\mathbf{\mu_0'}$}\\
\hline
\endfirsthead
\multicolumn{4}{c}%
{{  \tablename\ \thetable{} -- continued from previous page}}
\\
\hline
\multicolumn{1}{|c|}{$\mathbf{\CO}$} &\multicolumn{1}{c|}{$\chi$} &
\multicolumn{1}{c|}{$\mathbf{\mu_0}$}
&\multicolumn{1}{c|}{$\mathbf{\mu_0'}$} \\
\hline
\endhead

$D_7(a_1)$ &$ (0,1,1,2,3,4,5,0)+\nu\check\omega_1$ &$3240_z'$ &$1050_x'$\\
\hline
$E_6(a_1)+A_1$ &$(0,1,1,2,3,-\frac 72,-\frac 52,3)+\nu\check\omega_7$
&$4096_z'$ &$4096_x'$\\
\hline
$D_7(a_2)$ &$(0,1,1,2,2,3,4,0)+\nu\check\omega_1$ &$4200_x'$ &$3360_z'$\\
\hline
$D_5+A_2$ &$(0,1,2,3,-3,-2,-1,2)+\nu\check\omega_6$ &$4536_z'$ &$840_x'$\\
\hline
\end{longtable}
\end{center}
\end{tiny}

\subsubsection{$\mathbf{ A_4+A_2+A_1}$}\label{sec:A4A2A1}
 The centralizer is $A_1$, the
lowest W-type is $\mu_0=2835_x$, and the \ic is $(0,1,-\frac
52,-\frac 32,-\frac 12,\frac 12,$ $\frac 32,\frac 12)+\nu\check\omega_4$. 
The standard module is $X(A_4+A_2+A_1,St,\nu)$, $\nu\ge 0.$ The first
 reducibility point is at
 $\nu_0=\frac 3{10}$, where
the generic factor is parametrized by the nilpotent $\CO'=A_4+A_3$ and
 the W-type
$\mu_1=420_y$. There are exactly two composition factors at this
 point, one parametrized by $\CO$ (with lowest W-type $\mu_0$), and
 the generic factor. Then, either all the $W-$types in the generic
 factor change sign at $\nu=\nu_0$, or none of them do. A direct
 calculation shows that the determinant of the operator on the W-type
 $35_x'$ has opposite sign to the scalar on the sign representation
 $1_x'$, in this interval. It follows that
on the interval $(\frac 3{10},\frac 12),$ also $\mu_1$ has negative sign.

The next reducibility point is at $\nu=\frac 12$. A similar
argument as for the nilpotent $E_6+A_1$ (section \ref{sec:E6A1}),
shows that $\overline X(A_4+A_2+A_1,St,\nu)$ is not unitary for $\nu>1/2.$
The {complementary series} is $0\le\nu<\frac 3{10}$.

\medskip

For the rest of the nilpotents in $E_8$,  we check the details of the
argument outlined in the proof of proposition \ref{p:4.2} in
every case, and determine the correspondences between intertwining
operators on $W-$types and spherical operators on $W(\fz(\CO))$-types.
The exceptions (i.e., the nilpotent orbits for which the
complementary series is not the same as the $0$-complementary series
of the centralizer) are discussed separately. 
If $\Delta_1$ is a root system, and
$\Delta_2\subset \Delta_1$ is a subsystem, we denote by
$w_m(\Delta_1,\Delta_2)$, the element $w_0(\Delta_1)\cdot w_0(\Delta_2).$

\subsection{Single lowest $W$-type orbits} We begin with two
representative examples.

\subsubsection{$\mathbf{E_6}$}\label{sec:E6} The centralizer is $\fz(\CO)=G_2,$ the
lowest W-type is $\mu_0=525_x'$, and the \ic is
$(0,1,2,3,4,-4,-4,4)+\nu_2(0,0,0,0,0,0,1,1)+\nu_1(0,0,0,0,0,1,1,2),$
with $\nu_1\ge 0,\ \nu_2\ge 0$.
The standard module is $X(E_6,St,\nu),$ $\nu=(\nu_1,\nu_2).$ 

The subgroup $W(\fz)\cong W(G_2)\subset W$ is generated by:
\begin{equation}
\bar s_1=w_m(E_7,E_6), \quad\bar s_2=s_8.
\end{equation}
The intertwining operator $A(E_6,St,\nu)$ decomposes according to the
decomposition $w_m=(\bar s_1\cdot \bar s_2)^3.$

The restrictions of $W-$types are:

\begin{small}\begin{tabular} {lllllll}
Nilpotent &$E_6$ &$E_6A_1$ &$E_7(a_2)$ &$E_8(b_5)$ &$E_8(b_5)$ &$E_8(a_5)$ \\
W-type &$525_x'$ &$448_z'$ &$1344_x'$ &$1008_z'$ &$1400_z'$ &$700_x$'\\
Multiplicity &$1$ &$1$ &$2$ &$1$ &$2$ &$1$\\
$E_6\subset E_7$ &$21_b$ &$21_b$ &$27_a',\ 21_b$ &$27_a'$ &$27_a',\ 21_b$ &$27_a'$\\
$A_1$ &$(2)$ &$(11)$ &$(2),\ (11)$ &$(2)$ &$(2),\ (11)$ &$(11)$\\
$W(G_2)$ &$1_1$   &$1_4$ &$2_2$ &$1_3$ &$2_1$ &$1_2$
\end{tabular}\end{small}

On the factor corresponding to $\bar s_2$, the root $\al_8$ takes
values $3\nu_1+2\nu_2,\ 3\nu_1+\nu_2,$ and $\nu_2.$

The factor corresponding to $\bar s_1$ is induced from an intertwining
operator for the Hecke algebra of type $E_7,$ with the nilpotent orbit
$E_6$ in $E_7$, and infinitesimal character
$(0,1,2,3,4,-4,-4,4)+\bar\nu(0,0,0,0,0,1,-\frac 12,\frac 12),$ where
$\bar\nu$ takes the values $\nu_1,\ 2\nu_1+\nu_2$ and $\nu_1+\nu_2.$

The reducibility hyperplanes for $X(E_6,St,\nu)$ are:
$\nu_1,\nu_2+2\nu_1,\nu_1+\nu_2=1$ and
$2\nu_2+3\nu_1,\nu_2+3\nu_1,\nu_2=1$ (as in the centralizer $G_2$),
and $\nu_1,\nu_2+2\nu_1,\nu_1+\nu_2=5,9$.

The operators match as follows:
\begin{small}\begin{tabular}{l|llllll}
$\widehat W$ &$525_x'$
&$448_z'$
&$1344_x'$
&$1008_z'$
&$1400_z'$
&$700_x'$\\
\hline
$\widehat{W(G_2)}$
& $1_1$
&$1_4$
&$2_2$
&$1_3$
&$2_1$
&$1_2$,
\end{tabular}\end{small}
and all the relevant $W(G_2)$-types are matched.

\subsubsection{$\mathbf{D_4+A_1}$} The centralizer is $\fz(\CO)=C_3$,
the lowest W-type is $\mu_0=700_{xx}$, and the \ic is
$(0,1,2,3,-\frac 12,\frac
12,0,0)+(0,0,0,0,\nu_1,\nu_1,-\nu_2+\nu_3,\nu_2+\nu_3)$. 
The standard module is $X(D_4+A_1,St,\nu),$ where
$\nu=(\nu_1,\nu_2,\nu_3).$

The hyperplanes of reducibility are $\nu_i=\frac 12,$
$\nu_i\pm\nu_j=1$, as for the centralizer $C_3$, and  $\nu_i=\frac 32,\frac 72,\frac 92$,
$\pm\nu_i\pm\nu_j=4$ and $\pm\nu_1\pm\nu_2\pm\nu_3=\frac 32$.

The operators match as follows:
\begin{small}\begin{tabular}{l|llll}
$\widehat W$
&$700_{xx}$
&$2800_z$
&$6075_x$
&$5600_z$\\
\hline
$\widehat {W(C_3)}$
&$3\times 0$
&$ 0\times 3$
&$1\times 2$
&$ 0\times 12.$
\end{tabular}\end{small}

The $W-$types that match operators from $C_3$ are not sufficient for
concluding that the generic complementary series is included in the
one for $C_3$. They are positive in the unitary region for $C_3$:
$\{0\le\nu_3\le\nu_2\le\nu_1<\frac 12\}$, but also in the region
$\mathcal R=\{\nu_1+\nu_3>1,\ \nu_2>\frac 12,\ \nu_1-\nu_3<1,\
\nu_2+\nu_3<1\}.$

We also need to use the signature of the operator on the $W$-type
$4200_z$, which has multiplicity four. 

Among the extra hyperplanes of reducibility, $\nu_1+\nu_2-\nu_3=\frac
32$ cuts the region $\mathcal R$ into two open subregions:

$\mathcal R_1:$ $\nu_1+\nu_2-\nu_3<\frac 32$, sample point $(1,\frac
58,\frac 3{16})$. The determinant of $4200_z$ is negative in this
region.

$\mathcal R_2:$ $\nu_1+\nu_2-\nu_3>\frac 32$, and the determinant of
$4200_z$ is positive. We choose a point on the boundary of the
region, which is unitarily induced: $(\frac 78,\frac 78,
\frac 18)$. The corresponding parameter is induced from $D_7$, with
$(0,1,2,3,-\frac  38,\frac 58,\frac 74)$ in $D_7$. Furthermore, this
can be deformed irreducibly to $(0,1,2,3,-\frac 12,\frac 12,\frac
74)$, which unitarily induced from $(0,1,2,3,\frac 74)$ in $D_5$ and
the signatures are induced from $D_5$. In $D_5$, for this parameter,
$21^3\times 0$ and $1^4\times 1$ have opposite signs. 

Since $4200_z$ contains in his restriction $21^3\times 0$ and $1^4\times
1$, it follows that the form is indefinite on it, so the lowest W-type
factor is not unitary at this boundary point. But then, the entire
region $\mathcal R_2$ must be nonunitary.


\subsubsection{}\label{sec:match1} We list the matching of $W$-types for the other nilpotent orbits in
$E_8$ of similar kind. The infinitesimal characters $\chi=\frac
h2+\nu$ are in the tables in section \ref{tables}.

\begin{tiny}
\begin{center}
\begin{longtable}{|c|c|l|}\caption {Nilpotent orbits $\CO$ with single
  lowest $W$-type}\label{t:match1}\\

\hline
\multicolumn{1}{|c|}{$\mathbf{\CO}$}  
&\multicolumn{1}{c|}{$\mathbf{\fz(\CO)}$}
&\multicolumn{1}{c|}{\bf matching}\\
\hline
\endfirsthead
\multicolumn{3}{c}%
{{  \tablename\ \thetable{} -- continued from previous page}}
\\
\hline
\multicolumn{1}{|c|}{$\mathbf{\CO}$}  
&\multicolumn{1}{c|}{$\mathbf{\fz(\CO)}$}
&\multicolumn{1}{c|}{\bf matching}\\
\hline
\endhead

$\mathbf{E_6}$ &$G_2$ &{\begin{tabular}{llllll}
$525_x'$
&$448_z'$
&$1344_x'$
&$1008_z'$
&$1400_z'$
&$700_x'$\\
\hline
 $1_1$
&$1_4$
&$2_2$
&$1_3$
&$2_1$
&$1_2$
\end{tabular}}\\
\hline

$\mathbf{D_6}$ &$B_2$ &\begin{tabular}{llll}
$972_x'$
&$2268_x'$
&$3240_z'$
&$1050_x'$\\
\hline
$2\times 0$
&$11\times 0$
&$1\times 1$
&$0\times 2$
\end{tabular}\\
\hline

$\mathbf{A_6}$ &$2A_1$ &\begin{tabular}{lll}
$4200_z'$ &$6075_x'$ &$2835_x'$ \\
\hline
$(2)\otimes (2)$ &$(11)\otimes (2)$ &$(2)\otimes (11)$ 
\end{tabular}\\
\hline

$\mathbf{D_5+A_1}$ &$2A_1$ &\begin{tabular} {lll}
$3200_x'$  &$5600_z'$   &$6075_x'$\\
\hline
$(2)\otimes (2)$ &$(2)\otimes (11)$  &$(11)\otimes (11)$ 
\end{tabular}\\
\hline

$\mathbf{A_5+A_1}$ &$2A_1$ &\begin{tabular}{lll}
$2016_w$ &$3150_y$ &$4200_y$\\
\hline 
$(2)\otimes (2)$  &$(2)\otimes (11)$ &$(11)\otimes (2)$ 
\end{tabular}\\
\hline

$\mathbf{D_5}$ &$B_3$ &\begin{tabular}{lllll}
$2100_y$
&$3200_x'$
&$5600_z'$
&$2400_z'$
&$6075_x'$\\
\hline
$ 3\times 0$
&$ 12\times 0$
&$ 2\times 1$
&$ 0\times 3$
&$ 1\times 2$
\end{tabular}\\
\hline

$\mathbf{D_5(a_1)+A_1}$ &$2A_1$ &\begin{tabular}{lll}
$6075_x$ &$4200_z $ &$2400_z$\\
\hline
$(2)\otimes (2)$  &$(2)\otimes (11)$  &$(11)\otimes (2)$ 
\end{tabular}\\
\hline

$\mathbf{E_6(a_3)}$ &$G_2$ &\begin{tabular}{llllll}
 $5600_z$
& $3150_y$
& $7168_w$
& $5600_w$
& $4480_y$
& $1680_y$\\
\hline
 $1_1$
 &$1_4$
 &$2_2$
 &$2_2$
 &$2_1$
 &$1_3$.
\end{tabular}\\
\hline

$\mathbf{A_5}$ &$G_2+A_1$ &\begin{tabular}{lllllll}
$3200_x$ 
&$2016_w$ 
&$5600_z$ 
&$4200_y$ 
&$3150_y$ 
&$4480_y$ 
&$1680_y$ \\
\hline
$(1_1),(2)$
&$ (1_4),(2)$
&$ (1_1),(11)$
&$(2_2),(2)$
&$(1_4),(11)$
&$(2_1),(11)$
&$(1_3), (11)$
\end{tabular}\\
\hline

$\mathbf{A_4+A_2}$ &$2A_1$ &\begin{tabular}{lll}
$4536_z$ &$2835_x$ &$6075_x$\\
\hline
$(2)\otimes (2)$ &$(2)\otimes (11)$ &$(11)\otimes (2)$ 
\end{tabular}\\
\hline

$\mathbf{2A_3}$ &$B_2$ &\begin{tabular}{lll}
$840_x$ 
&$4200_x$ 
&$4536_z$ \\
\hline
$2\times 0$
&$ 11\times 0$
&$ 1\times 1$
\end{tabular}\\
\hline

$\mathbf {D_4+A_1}$ &$C_3$ &
\begin{tabular}{llll}
 $700_{xx}$
&$2800_z$
&$6075_x$
&$5600_z$\\
\hline
 $3\times 0$
&$ 0\times 3$
&$1\times 2$
&$ 0\times 12 $
\end{tabular}\\
\hline

$\mathbf{A_3+A_2+A_1}$ &$2A_1$ &\begin{tabular}{lll}
$1400_{zz}'$ &$4096_x'$ &$2240_x'$\\
\hline
$(2)\otimes (2)$ &$(11)\otimes (2)$  &$(2)\otimes (11)$
\end{tabular}\\
\hline

$\mathbf{A_3+2A_1}$ &$B_2+A_1$ &\begin{tabular}{llll}
$1050_x$ &
$1400_x$ &
$972_x$ &
$3240_z$ \\
\hline
$(2\times 0)\otimes (2)$
&$(0\times 2)\otimes (2)$
&$(11\times 0)\otimes (2)$
&$(0\times 2)\otimes (11)+(1\times 1)\otimes (2)$.
\end{tabular}\\
\hline

$\mathbf{2A_2+2A_1}$ &$B_2$ &\begin{tabular}{llll}
$175_x$ 
&$1050_x$ 
&$972_x$ 
&$3240_z$
\\
\hline
$2\times 0$
&$11\times 0$
&$0\times 11$
&$1\times 1$\footnotemark
\end{tabular}\\
\hline

$\mathbf{D_4}$ &$F_4$ &\begin{tabular}{lllllll}
$525_x$ &
$700_{xx}$ &
$2800_z$ &
$2100_x$ &
$6075_x$ &
$4200_z$ &
$5600_z$ \\
\hline
$1_1$
&$2_3$
&$4_2$
&$2_1$
&$9_1$
&$8_3$
&$8_1$.
\end{tabular}\\
\hline

$\mathbf{A_3+A_1}$ &$B_3+A_1$ &\begin{tabular}{lllll}
$1344_x$ &
$1400_z$ &
$1050_x$ &
$1400_x$ &
$350_x$\\
\hline
$(3\times 0)\otimes (2)$
&$(3\times 0)\otimes (11)$
&$(21\times 0)\otimes (2)$
&$(2\times 1)\otimes (11)$
&$(0\times 3)\otimes (11)$.
\end{tabular}\\
\hline

$\mathbf{2A_2+A_1}$ &$G_2+A_1$ &\begin{tabular}{llll}
$448_z$ &
$1344_x$ &
$175_x$ &
$1050_x$\\
\hline
$1_1\otimes (2)$
&$1_1\otimes (11)$
&$1_4\otimes (2)$
&$2_2\otimes (11)$.
\end{tabular}\\
\hline

$\mathbf{A_3}$ &$B_5$ &\begin{tabular}{lllllll}
$567_x$ &
$1344_x$ &
$1400_z$ &
$56_z$ &
$1050_x$ &
$1400_x$ &
$350_x$ \\
\hline
$5\times 0$
&$41\times 0$
&$4\times 1$
&$0\times 5$
&$32\times 0$
&$3\times 2$
&$1\times 4$.
\end{tabular}\\
\hline

$\mathbf{A_2+3A_1}$ &$G_2+A_1$ &\begin{tabular}{llllll}
$400_z$ &
$700_x$ &
$448_z$ &
$1344_x$ &
$1008_z$ &
$1400_z$ \\
\hline
$1_1\otimes (2)$
 &$1_1\otimes (11)$
 &$1_4\otimes (11)$
 &$2_2\otimes (11)$
 &$1_3\otimes (11)$
 &$2_1\otimes (11)$.
\end{tabular}\\
\hline

$\mathbf{A_2+2A_1}$ &$B_3+A_1$ &\begin{tabular}{lllll}
$560_z$ &
$567_x$ &
$400_z$ &
$700_x$ &
$300_x$ \\
\hline
$(3\times 0)\otimes (2)$
&$(3\times 0)\otimes (11)$
&$(12\times 0)\otimes (2)$
 &$(2\times 1)\otimes (2)$
 &$(0\times 3)\otimes (2)$.
\end{tabular}\\
\hline

$\mathbf{4A_1}$ &$C_4$ &\begin{tabular}{lllll}
$50_x$ &$210_x$ &$560_z$ &$567_x$ &$300_x$\\
\hline
$4\times 0$ &$0\times 4$ &$1\times 3$ &$0\times
13$ &$0\times 22$.
\end{tabular}\\
\hline

$\mathbf{3A_1}$ &$F_4+A_1$ &\begin{tabular}{lllll}
$84_x$ &
$112_z$ &
$50_x$ &
$210_x$ &
$160_z$ \\
\hline
$1_1\otimes (2)$
&$1_1\otimes (11)$
&$2_3\otimes (2)$
&$4_2\otimes (11)$
&$2_1\otimes (11)$.
\end{tabular}\\
\hline

$\mathbf{2A_1}$ &$B_6$ &\begin{tabular}{lllll}
$35_x$&
$84_x$&
$112_z$ &
$28_x$&
$50_x$\\
\hline
$6\times 0$
&$15\times 0$
&$5\times 1$
&$0\times 6$
&$33\times 0$.
\end{tabular}\\
\hline

$\mathbf{A_1}$ &$E_7$
&\begin{tabular}{llll}
$8_z$ &$35_x$ &$84_x$ &$50_x$\\
\hline
$1_a$ &$7_a$ &$21_b'$ &$15_a'$
\end{tabular}\\
\hline

\hline
\end{longtable}
\end{center}
\end{tiny}

\footnotetext{$\CO=2A_2+2A_1$: although not identical with $1\times 1$, the
      operator on $3240_z$ has the same signature as $1\times 1$ in
      the open regions of $B_2.$}

\subsection{Exceptions}\label{exceptions}

\subsubsection{ $\mathbf{4A_1}$}\label{sec:4A1} We present the
case of the complementary series for the nilpotent orbit $4A_1$
in detail. This is the only case in which the complementary series is
{larger} than the 0-complementary series of the centralizer, which is
of type $C_4$. 

The standard module is $X(4A_1,St,\nu)$,
$\nu=(\nu_1,\nu_2,\nu_3,\nu_4),$ and it has lowest $W-$type 
$\mu_0=50_x$. The infinitesimal character is $(0,1,-\frac 12,\frac 12,-\frac
12,\frac 12,0,0)+(0,0,\nu_1,\nu_1,\nu_2,\nu_2,-\nu_3+\nu_4,\nu_3+\nu_4).$

The operators match as follows:
\begin{small}\begin{tabular}{l|lllll}
$\widehat W$ &$50_x$ &$210_x$ &$560_z$ &$567_x$ &$300_x$\\
\hline
$\widehat {W(C_4)}$ &$4\times 0$ &$0\times 4$ &$1\times 3$ &$0\times
13$ &$0\times 22$
\end{tabular}\end{small}.

These $W-$types only change sign when passing a hyperplane as in $C_4$:
$\nu_i=\frac 12$ and $\pm\nu_j+\nu_i=1$. We know that the region
$\nu_4<\frac 12$ is the only unitary 0-complementary series in
$C_4$. The $W-$types above  are not sufficient however to rule out all
other (4-dimensional) open regions in $C_4$. They are all positive
semidefinite also in the region
$\CR=\{\nu_1+\nu_4<1,\nu_2+\nu_3<1,\nu_2+\nu_4>1,\nu_3>\frac 12\}$.

The hyperplanes of reducibility $-\nu_2+\nu_3+\nu_4=\frac 32$,
$-\nu_1+\nu_3+\nu_4=\frac 32$ and $\nu_1+\nu_3+\nu_4=\frac 32$ cut the
region $\C R$  into the following open regions:

\begin{enumerate}

\item $\CR_1=\{\nu_1+\nu_4<1,\nu_2+\nu_3<1,\nu_2+\nu_4>1,\nu_3>\frac
  12,-\nu_2+\nu_3+\nu_4>\frac 32\}$,

\item $\CR_2=\{\nu_1+\nu_4<1,\nu_2+\nu_3<1,\nu_2+\nu_4>1,\nu_3>\frac
  12,-\nu_2+\nu_3+\nu_4<\frac 32<-\nu_1+\nu_3+\nu_4\}$,

\item $\CR_3=\{\nu_1+\nu_4<1,\nu_2+\nu_3<1,\nu_2+\nu_4>1,\nu_3>\frac
  12,-\nu_1+\nu_3+\nu_4<\frac 32<\nu_1+\nu_3+\nu_4\}$,

\item $\CR_4=\{\nu_1+\nu_4<1,\nu_2+\nu_3<1,\nu_2+\nu_4>1,\nu_3>\frac
  12,\nu_1+\nu_3+\nu_4>\frac 32\}$.

\end{enumerate}

In $\CR_1$, $\CR_2$ and $\CR_4$,  one
can deform the parameter to $\nu_1=0$, where the module is
unitarily induced irreducible from a nonunitary module attached to
$4A_1$ in $E_7$.

\begin{proposition} The open region $\C R_3$ is unitary:
\begin{equation}\label{eq:4A1}\{\nu_1+\nu_4<1,\ \nu_2+\nu_3<1,\
\nu_2+\nu_4>1,\ -\nu_1+\nu_3+\nu_4<\frac 32<\nu_1+\nu_3+\nu_4\}.
\end{equation}  
\end{proposition}

\begin{proof} We divide the proof into four parts.

\smallskip

\noindent{\bf Step 1.} {\it The generic modules are unitary on the walls of $\CR_3$.}

For each wall, we find the nilpotent orbit $\CO'$ parameterizing the
generic module:  $\nu_1+\nu_4=1$, $\nu_2+\nu_3=1$, $\nu_2+\nu_4=1$
correspond to $A_2+2A_1$, and $\nu_1+\nu_3+\nu_4=\frac 32$,
$-\nu_1+\nu_3+\nu_4=\frac 32$ correspond to $A_2+3A_1.$ 
The claim follows thenby comparison with the complementary series
attached to the nilpotent orbits $A_2+2A_1$ and $A_2+3A_1.$

\smallskip

We deform the parameter to a particular point on the walls: $p=(\frac
1{12},\frac 3{12},\frac 9{12},\frac 9{12})$. The point $p$ lies at the
intersection of the walls $\nu_2+\nu_3=1$ and $\nu_2+\nu_4=1$. The
corresponding point $\bar p=\frac 12 h+p$, in $E_8$-coordinates is
$\bar p=(0,1,-\frac 5{12},\frac 7{12},-\frac 3{12},\frac 9{12},\frac
{18}{12},0)$.

\noindent{\bf Step 2.} {\it The standard module $X(4A_1,St,p)$ has two
  composition factors: $\overline X(4A_1,St,p)$ and $X(A_2+2A_1,St,p)$.}

The standard module $X(4A_1,St,p)$ is reducible. A necessary condition
for a nilpotent $\CO'>\CO$ to parameterize a composition factor is 
that $w\bar p=\frac 12 h'+\nu'$, for $h'$ the middle
element of a Lie triple $\{e',h',f'\}$ of $\CO'$, $w\in W$ and
$\nu'\in\fz(\CO').$  We check that the
nilpotent $\CO'$ satisfying the condition are $A_2+A_1$
and $A_2+2A_1$, so potentially there are three factors. $A_2+2A_1$
parameterizes the generic factor. The lowest $W$-type
of $A_2+A_1$ is $210_x$, and the operator on $210_x$ matches $0\times
4$ in $C_4$, so it is invertible at $p.$

\smallskip

\noindent{\bf Step 3.} {\it The non-generic factor $\overline
  X(4A_1,St, p)$ is unitary.} 

The point $\bar p$ is unitarily induced reducible from
$D_7$. The corresponding nilpotent in $D_7$ is $(1^32^43)$ and the
infinitesimal character is of the form $(0,1,-\frac 12+\bar\nu_1,
\frac 12+\bar\nu_1,-\frac 12+\bar\nu_2,\frac 12+\bar\nu_2,\bar\nu_3)$,
with $(\bar\nu_1,\bar\nu_2,\bar\nu_3)=(\frac 1{12},\frac 3{12},\frac
{18}{12})$. Moreover, the parameter in $D_7$ can be deformed
irreducibly to a unitarily induced from $D_3\times GL(4)$, where the
parameter on $D_3$ is $(0,1,\frac {18}{12})$ (nilpotent $(1^33)$ in
$D_3$) and on $GL(4)$, it is $(-\frac 12,\frac 12,-\frac 12,\frac
12)+(\frac 1{12},-\frac 1{12},-\frac 1{12},\frac 1{12})$ (nilpotent
$(22)$). Therefore, the signature of the form on $E_8$ can be computed
from the signatures on $D_3$ and $GL(4)$:

\begin{small}\begin{tabular}{l|l|l|ccc}

$D_3$ &$(0,1,\nu)$ &$\nu=\frac {18}{12}$ &$111\times 0$ &$11\times 1$ &$12\times 0$\\ 
&&&$+$ &$+$ &$-$\\
\hline
$GL(4)$ &$(-\frac 12+\nu,\frac 12+\nu,-\frac 12-\nu,\frac 12-\nu)$
&$\nu=\frac 1{12}$ &$(22)$ &$(211)$ &$(1^4)$\\
&&&$+$&$+$&$+$
\end{tabular}\end{small}

The signature of the hermitian form on $D_3\times GL(4)$ will
therefore be $(24,12)$. The unitarily induced form on $D_7$ will have
signature $(13440,6720)$ and the unitarily induced form in $E_8$ has
signature $(29030400,14515200)$. Since the W-dimension of
$X(A_2+2A_1)$ is $|W|/|W(A_2)W(A_1)^2|=29030400$ and
$X(A_2+2A_1)$ is unitary, it follows that the induced form on the
factor $\overline X(4A_1)$ is (negative) definite, so after the appropriate
normalization, the factor $\overline X(4A_1)$ is also unitary.

\smallskip

\noindent{\bf Step 4.} {\it In the interior of $\CR_3$, the
  intertwining operator $A_\mu(4A_1,St,\nu)$ is positive definite
  for all $\mu \in \widehat W$.}

It is sufficient to calculate the intertwining
operator on a single W-type which appears in both factors $\overline X(4A_1,St,p)$
and $X(A_2+2A_1,St,p)$. The W-type $560_z$ has this property, and  $A_{560_z}(4A_1,St,\nu)=a_{1\times
  3}(\nu)$ which is positive definite inside $\CR_3$. This concludes the proof. 

\end{proof}

\subsubsection{
  $\mathbf{D_5(a_1)+A_1}$}\label{sec:D5(a1)A1} The centralizer is $\fz(\CO)=2A_1$, the lowest
  W-type is $\mu_0=6075_x$, and the \ic is $(0,1,1,2,3,$ $-\frac
  12+\nu_2,\frac 12+\nu_2,2\nu_1).$
The standard module is $X(D_5+A_1,\sigma\otimes St,\nu)$, where
$\nu=(\nu_1,\nu_2)$ and $\sigma$ is the discrete series parametrized by
the nilpotent $(73)$ in $D_5.$

The matching of operators in table \ref{t:match1} imply that the
complementary series is included in
$\{0\le\nu_1<1,0\le\nu_2<1/2\}$, the complementary series of the
  centralizer. 
There are hyperplanes of reducibility $2\nu_1\pm\nu_2=\frac 32$ which
cut this region. We need to use the scalar operator on $1344_w$ (a $W$-type
with multiplicity one). This is negative
in the region $\{2\nu_1-\nu_2<\frac 32<2\nu_1+\nu_2, \nu_2<\frac 12\}$.
It follows that the {complementary series} is $\{0\le\nu_2<\frac 12$,
$2\nu_1+\nu_2<\frac 32\}$ and  $\{0\le\nu_1<1,2\nu_1-\nu_2>\frac
32\}$.

\subsubsection{ $\mathbf{A_4+A_2}$}\label{sec:A4A2}
The centralizer is  $\fz(\CO)=2A_1$, the lowest W-type is
$4536_z$, and the
infinitesimal character is $s=(-\frac 12,\frac 12,-\frac 52,-\frac
32,-\frac 12,\frac 12,\frac 32,\frac 12)+$ $\nu_2(1,1,0,0,0,0,0,0)+$
$\nu_1(0,0,1,1,1,1,1,5),$ with $\nu_1\ge 0$ and $\nu_2\ge 0.$ 
The standard module is $X(A_4+A_2,St,\nu),$ $\nu=(\nu_1,\nu_2).$ 

The matching of operators in table \ref{t:match1} imply that the
complementary series is included in
$\{0\le\nu_1<1/2,0\le\nu_2<1/2\}$, the complementary series of the
  centralizer. 
There are hyperplanes of reducibility $5\nu_1\pm\nu_2=2$ which
cut this region. We need to use the scalar operator on $420_y$ (a $W$-type
with multiplicity one). This is negative
in the region $\{5\nu_1-\nu_2<2<5\nu_1+\nu_2, \nu_2<\frac 12\}$.
It follows that the {complementary series} is 
$\{0\le\nu_2<\frac 12, 5\nu_1+\nu_2<2\}
  \cup\{0\le\nu_1<\frac 12, 5\nu_1-\nu_2>2\}.$

\subsubsection{ $\mathbf{A_2+3A_1}$}\label{sec:A2+3A1} The centralizer is
  $\fz(\CO)=G_2+A_1$, the lowest W-type is $\mu_0=400_z$, and the
  \ic is 
  $(0,1,-1,0,-1,0,-\frac 12,\frac 12)+\nu_1 (0,0,1,1,1,1,-2,2)+\nu_2
  (0,0,0,0,1,1,-1,1)+\nu_3 (0,0,0,0,0,0,1,1)$. 
The standard module is $X(A_2+3A_1,St,\nu),$
$\nu=(\nu_1,\nu_2,\nu_3).$ 

The matching of operators in table \ref{t:match1} imply that the
complementary series is included in $\{3\nu_1+2\nu_2<1,\nu_3<\frac
12\}$ and $\{3\nu_1+\nu_2>1>2\nu_1+\nu_2,\nu_3<\frac 12\}$, the complementary series of the
  centralizer. 
There are hyperplanes of reducibility $3\nu_1+2\nu_2+\nu_3=\frac 32$, $3\nu_1+2\nu_2-\nu_3=\frac
32$ and $3\nu_1+\nu_2+\nu_3=\frac 32$ which
cut the second region into five (open) subregions. We need to use the determinant of the operator on $175_x$ (a $W$-type
with multiplicity two). This is negative in two of the five
subregions.

It follows that the complementary series is the union of
$\{3\nu_1+\nu_2<1,0\le\nu_3<\frac 12\}$,  $3\nu_1+2\nu_2+\nu_3<\frac
32,$ $3\nu_1+2\nu_2-\nu_3>\frac 32,$ and $\mathcal R_4$: $3\nu_1+2\nu_2-\nu_3<\frac 32.$

\subsubsection{$\mathbf{A_2+2A_1}$}\label{sec:A2+2A1} The centralizer
is $\fz(\CO)=B_3+A_1$, the lowest W-type is $\mu_0=560_z$, and the
\ic is $(0,1,-1,0,1,0,0,0)+$ $(0,0,\nu_1,\nu_1,\nu_1,\nu_2,\nu_3,\nu_4)$.
The standard module is $X(A_2+A_1,St,\nu),$
$\nu=(\nu_1,\nu_2,\nu_3,\nu_4).$ 

The matching of operators in table \ref{t:match1} imply that the
complementary series is included in
the complementary series for the centralizer
$B_3+A_1:$ $\mathcal R_1=\{0\le\nu_1<1,\ \nu_3+\nu_4<1\}$ and
$\mathcal R_2=\{0\le\nu_1<1,\nu_2+\nu_4>1,\nu_2+\nu_3<1,\nu_4<1\}.$
There are hyperplanes of reducibility 
$3\nu_1+\nu_2+\nu_3-\nu_4=3$,
$3\nu_1-\nu_2-\nu_3+\nu_4=3$,
$3\nu_1+\nu_2-\nu_3+\nu_4=3$,
$3\nu_1-\nu_2+\nu_3+\nu_4=3$,
$3\nu_1+\nu_2+\nu_3+\nu_4=3$, which cut $\mathcal R_1$ and $\mathcal
R_2$ into twelve open subregions. 
We need to use the determinant of the operator on $448_z$ (a $W$-type
with multiplicity four). This is negative in five subregions, the
other seven forming the complementary series (see section \ref{tables}
for the explicit description).

\subsection{Multiple lowest $W$-types orbits} We begin with two
typical examples.

\subsubsection{$\mathbf{D_4(a_1)}$}\label{sec:D4(a1)} The centralizer is $\fz(\CO)=D_4$,
with component group $A(\CO)=S_3.$ The \ic is $(0,1,1,2,0,0,0,0)+$
$(0,0,0,0,\nu_4,\nu_3,\nu_2,\nu_1)$. 

The standard module $X(D_4,\sigma,\nu),$
$\nu=(\nu_1,\nu_2,\nu_3,\nu_4),$ with $\sigma$ the discrete series
parametrized by the nilpotent $(53)$ in $D_4,$ has three lowest
$W-$types, $\mu_0=1400_z,$ $\mu_0'=1008_z,$ and $\mu_0''=56_z.$
Note that $\mu_0'$ has multiplicity two. They have the same
signature at $\infty$, and stay in the same factor unless the
parameter satisfies $\nu_4=0$ or
$\nu_1-\nu_2-\nu_3-\nu_4=0.$

If, for example, $\nu_4=0,$ the standard
module corresponding to the generic case is $X(D_5,\sigma',\nu'),$
$\nu'=(\nu_1,\nu_2,\nu_3),$ where $\sigma'$ is the generic limit of
discrete series parametrized by the nilpotent $(5311)$ in $D_5$, and it
contains two lowest $W-$types, $\mu_0$ and $\mu_0'.$

If, $\nu_4=0$ and $\nu_1-\nu_2-\nu_3-\nu_4=0,$ the
standard module corresponding to the generic case is $X(E_6,\sigma'',\nu''),$ $\nu''=(\nu_1,\nu_2)$, where $\sigma''$ is the generic
limit of discrete series module parametrized by the nilpotent orbit $D_4(a_1)$ in
$E_6$, and it contains a single lowest W-type, $\mu_0.$

The subgroup $C(\fa,M)\cong W(F_4)$ is generated by: 
\begin{equation}
\bar s_1=s_8,\quad
\bar s_2=s_7,\quad
\bar s_3=w_m(D_5(2),D_4),\quad
\bar s_4=w_m(D_5(1),D_4),
\end{equation}
and the subgroup $W(\fz)\cong W(D_4)$ by $\{\bar s_3\cdot\bar s_2\cdot
\bar s_3,\  \bar s_2,\ \bar s_1,\ \bar s_4\cdot \bar s_3\cdot \bar s_2\cdot
\bar s_3\cdot \bar s_4\}.$ 

The restrictions of $W-$types are:

\begin{small}\begin{tabular}{llllll}
Nilpotent &$D_4(a_1)$ &$D_4(a_1)$ &$D_4(a_1)$  &$D_4(a_1)A_1$ &$A_3A_2$\\
W-type &$1400_z'$ &$1008_z'$ &$56_z'$ &$1400_x'$ &$3240_z'$  \\
Multiplicity &$1+0+0$ &$0+1+0$ &$0+0+1$  &$2+1+0$  &$3+3+0$  \\
$D_5$ &$211\times 1$ &$211\times 1$ &$1^3\times 2$  &$2\cdot 211\times 1$  &$6\cdot 211\times 1$ \\
&&$1^3\times 2$ &&$2\cdot 1^3\times 2$ 
&$3\cdot 1^3\times 2$ \\

$A_1$ &$(2)$ &$2\cdot (2)$ &$(2)$    &$3\cdot (2),(11)$    &$6\cdot
(2),3\cdot (11)$   \\

$W(F_4)$ &$1_1$ &$2_1$ &$1_2$ &$4_2$  &$9_1$ \\
$D_4$ &$4\times 0$ &$2\cdot 4\times 0$ &$4\times 0$ &$3\times 1$
&$2\times 2, 31\times 0$. \\

\end{tabular}\end{small}

In addition to the hyperplanes of reducibility as in $D_4$, there are
the following reducibility hyperplanes:
$\nu_i=2,3$, $i=1,4$,
$\pm\nu_1\pm\nu_2\pm\nu_3\pm\nu_4=4,6$.

The operators (normalized by the scalar on $\mu_0$) match operators
for the Hecke algebra of type $F_4$ with parameter $0$ on the long
roots, or equivalently operators for the Hecke algebra of type $D_4$
(see section \ref{sec:2.8d}):

\begin{small}\begin{tabular}{l|lllll}
$\widehat W$ &
$1400_z$ &
$1008_z$ &
$56_z$  &
$1400_x$ &
$3240_z$ \\
\hline
$\widehat{W(F_4)}$
 &$1_1$  
&$2_1$  
&$1_2$   
&$4_2$  
&$9_1$\\
\hline  
$\widehat{W(D_4)}$
&$4\times 0$
&$2\cdot 4\times 0$
&$4\times 0$
&$3\times 1$
&$2\times 2+31\times 0$.
\end{tabular}\end{small}

\subsubsection{$\mathbf{A_4}$} We realize the Bala-Carter Levi
subalgebra $\fk m=\{\al_5,\al_6,\al_7,\al_8\}.$ The centralizer of the
nilpotent orbit is $\fz(\CO)=A_4,$ and it is realized by
$\{\al_3,\al_1,$ $(-\frac 12,-\frac 12,\frac 12,\frac 12,\frac 12,\frac
12,\frac 12,\frac 12),\al_2\}.$ The infinitesimal character is
$\chi=(0,0,-2,-1,$ $0,1,2,0,0)+(\nu_4,-\nu_1+\nu_2,\nu_3,\nu_3,\nu_3,\nu_3,\nu_3,\nu_1+\nu_2).$ 

The standard module $X(A_4,St,\nu),$ $\nu=(\nu_1,\nu_2,\nu_3,\nu_4)$
has two lowest $W-$types, $\mu_0=2268_x$ and $\mu_0'=1296_z.$ They
have opposite signs at infinity, and they are separate if and only is
$\nu_3=\nu_4=0.$ We will assume that this is the case; therefore,
$\chi=(0,-\nu_1+\nu_2,-2,-1,0,1,2,\nu_1+\nu_2).$ The standard module
corresponding to the generic case is $X(D_6,\sigma,\nu),$
$\nu=(\nu_1,\nu_2),$ where $\sigma$ is the generic limit of discrete
series parametrized by the nilpotent $(5511)$ in $D_6.$

The subgroup $W(\fz)\cong W(A_4)\subset W$ is generated by:
\begin{equation}
\bar s_1=s_3,\ 
\bar s_2=s_1,\ 
\bar s_3=w_m(A_5,A_4)\cdot w_m(D_6,A_4)\cdot s_1\cdot
  w_m(D_6,A_4)\cdot w_m(A_5,A_4),\ 
\bar s_4=s_2.
\end{equation}
The intertwining operator decomposes according to the decomposition
$w_m=\bar s_1\cdot \bar s_2\cdot \bar s_3\cdot \bar s_4\cdot \bar
s_1\cdot \bar s_2\cdot \bar s_3\cdot \bar s_1\cdot \bar s_2\cdot \bar s_1.$

We compute the restrictions of $W$-types as in section
\ref{sec:D4(a1)}. The operators on  $W-$types in the generic factor of
$X(A_4,St,\nu)$ match
hermitian spherical operators in $A_4$ as follows:
\begin{small}\begin{tabular}{l|lllll}
$W$-type &$2268_x$ &$4096_x$ &$4096_z$ &$4200_x$ &$3360_z$\\
\hline
$W(A_4)$-type &$(5)$ & $(41)$ & $(41)$ & $(32)$ &$(32)$\\
eigenspace of $w_0(A_4)$ &$+1$-eig. &$+1$-eig. &$-1$-eig. &$+1$-eig. &$-1$-eig.
\end{tabular}\end{small}.

\subsubsection{} We list the matching of $W$-types for the other
nilpotent orbits of similar kind. The infinitesimal characters are in
the tables in section \ref{tables}.

\begin{tiny}
\begin{center}
\begin{longtable}{|c|c|l|}\caption {Nilpotent orbits $\CO$ with multiple
  lowest $W$-type}\\

\hline
\multicolumn{1}{|c|}{$\mathbf{\CO}$}  
&\multicolumn{1}{c|}{$\mathbf{\fz(\CO),A(\CO)}$}
&\multicolumn{1}{c|}{\bf matching}\\
\hline
\endfirsthead
\multicolumn{3}{c}%
{{  \tablename\ \thetable{} -- continued from previous page}}
\\
\hline
\multicolumn{1}{|c|}{$\mathbf{\CO}$}  
&\multicolumn{1}{c|}{$\mathbf{\fz(\CO),A(\CO)}$}
&\multicolumn{1}{c|}{\bf matching}\\
\hline
\endhead

$\mathbf {E_6(a_1)}$ &$A_2,\bZ_2$ &\begin{tabular}{lll}
$2800_z'$ &$4096_z$ &$4096_x$\\
\hline
$(3)$   &$(21)$ &$(21)$ \\
$+1$-eig. &$+1$-eig. &$-1$-eig.
\end{tabular}\\
\hline

$\mathbf{D_6(a_1)}$ &$2A_1,\bZ_2$ &\begin{tabular}{lll}
$5600_z'$ &$2400_z'$ &$6075_x'$\\
\hline
$(2)\otimes (2)$ &$(2)\otimes (2)$ &$(2)\otimes (11)+(11)\otimes (2).$
\end{tabular}\\
\hline

$\mathbf{D_6(a_2)}$ &$2A_1,\bZ_2$   &\begin{tabular}{lll}
$4200_y$  &$2688_y$  &$7168_w$\\
\hline
$(2)\otimes (2)$  &$(2)\otimes (2)$ &$(2)\otimes (11)+(11)\otimes (2)$
\end{tabular}\\
\hline

$\mathbf{D_4+A_2}$  &$A_2,\bZ_2$  &\begin{tabular}{lll}
$4200_z$  &$1344_w$  &$3150_y$\\
\hline
$(3)$ &$(21)$ &$(111)$ \\
$+1$-eig. &$+1$-eig &$+1$-eig.
\end{tabular}\\
\hline

$\mathbf{A_4+2A_1}$  &$A_1+T_1,\bZ_2$ &\begin{tabular}{ll}
$4200_x$ &$4536_z$\\
\hline
$(2)$ &$(11)$
\end{tabular}\\
\hline

$\mathbf{D_5(a_1)}$ &$A_3,\bZ_2$ &\begin{tabular}{lll}
$2800_z$ &$6075_x$ &$4200_z$\\
\hline
$(4)$ &$(31)$ &$(22)$ \\
$+1$-eig. &$+1$-eig. &$+1$-eig.
\end{tabular}\\
\hline

$\mathbf{A_4+A_1}$ & $A_2+T_1,\bZ_2$ &\begin{tabular}{lll}
$4096_x$ &$4200_x$ &$3360_z$\\
\hline
$(3)$ &$(21)$&$(21)$\\
$+1$-eig. &$+1$-eig. &$-1$-eig.
\end{tabular}\\
\hline

$\mathbf{A_4}$ &$A_4,\bZ_2$ &\begin{tabular}{lllll}
$2268_x$ &$4096_x$ &$4096_z$ &$4200_x$ &$3360_z$\\
\hline
$(5)$ & $(41)$ & $(41)$ & $(32)$ &$(32)$\\
$+1$-eig. &$+1$-eig. &$-1$-eig. &$+1$-eig. &$-1$-eig.
\end{tabular}\\
\hline

$\mathbf{D_4(a_1)+A_2}$ &$A_2,\bZ_2$ &\begin{tabular}{lll}
$2240_x$ &$4096_x$ &$4096_z$\\
\hline
$(3)$&$(21)$&$(21)$\\
$+1$-eig.  &$+1$-eig. &$-1$-eig.
\end{tabular}\\
\hline

$\mathbf{A_3+A_2}$ &$B_2+T_1,\bZ_2$ &\begin{tabular}{llll}
$3240_z$ &
$1400_{zz}$ &
$2240_x$ &
$840_z$ \\
\hline
$ 2\times 0$
&$ 11\times 0$
&$1\times 1$
&$0\times 2$.
\end{tabular}\\
\hline

$\mathbf{D_4(a_1)+A_1}$ &$3A_1,S_3$ &\begin{tabular}{ll}
$1400_x$ &$3240_z$\\
\hline
$(2)\otimes(2)\otimes(2)$ &$(11)\otimes (2)\otimes (2)+(2)\otimes
  (11)\otimes (2)+(2)\otimes (2)\otimes (11)$.
\end{tabular}\\
\hline

$\mathbf{D_4(a_1)}$ &$D_4,S_3$ &\begin{tabular}{lllll}
$1400_z$ &
$1008_z$ &
$56_z$  &
$1400_x$ &
$3240_z$ \\
\hline  
$4\times 0$
&$2\cdot 4\times 0$
&$4\times 0$
&$3\times 1$
&$2\times 2+31\times 0$
\end{tabular}\\
\hline

$\mathbf{2A_2}$  &$2G_2,\bZ_2$ &\begin{tabular}{llllll}
$700_x$ &$300_x$ &$448_z$ &$1344_x$ &$1400_z$ &$1008_z$\\
\hline
$1_1\otimes 1_1$ & &$1_4\otimes
1_1$ &$2_2\otimes 1_1$&$2_1\otimes 1_1$ &$1_3\otimes 1_1$\\
 &$1_1\otimes 1_1$ &$1_1\otimes 1_4$ &$1_1\otimes 2_2$ &$1_1\otimes 2_1$ &$1_1\otimes
1_3$.
\end{tabular}\\
\hline

$\mathbf{A_2+A_1}$ &$A_5$ &\begin{tabular}{lll}
$210_x$ &$560_z$ &$400_z$\\
\hline
$(6)$ &$(51)$ &$(33)$\\
$+1$-eig. &$+1$-eig. &$+1$-eig.\\
\end{tabular}\\
\hline

$\mathbf{A_2}$ &$E_6$ &\begin{tabular}{llllllll}

$112_z$ &
$210_x$ &
$160_z$ &
$560_z$ &
$400_z$ &
$700_x$ &
$300_x$ \\
\hline
$1_p$ 
&$6_p$ 
&$6_p$ 
&$20_p$
&$15_q$
&$30_p$
&$15_p$\\
 $+1$-eig.
 &$+1$-eig.
 &$-1$-eig.
 &$+1$-eig.
 &$+1$-eig.
 &$+1$-eig.
 &$-1$-eig.
\end{tabular}\\

\hline
\end{longtable}
\end{center}
\end{tiny}

\section{Tables of generic unitary parameters}\label{tables}


\subsection{Parameters for $\CO\neq 0$}
We give tables which contain the nilpotent orbits
(see \cite{Ca}), the hermitian infinitesimal character, and the
coordinates and type of the centralizer. 


The nilpotent orbits which are exceptions are marked with $*$ in the
tables. The description of the complementary series for them is
recorded after the tables.
For the rest of the nilpotents, an \ic $\chi$ is in the
complementary series if and only if the corresponding parameter $\nu$
is in the $0$-complementary series for $\fz(\CO).$ The parameter $\nu$
is given by a string $(\nu_1,\dots,\nu_k)$, and the order agrees with
the way the centralizer $\fz(\CO)$ is written in the tables. The parts of
$\nu$ corresponding to a torus $T_1$ or $T_2$  in $\fz(\CO)$ must be
$0$, in order for $\chi$ to be unitary. In addition, if $\nu$
corresponds to $A_1$, the complementary series is $0\le\nu<\frac 12,$
while the notation $A_1^\ell$ means that it is $0\le\nu<1.$ If a string
$(\nu_1,\dots,\nu_k)$ of $\nu$ corresponds to type $A_k,$ the last
$k-[\frac k2]$ entries must be $0$ in order for $\chi$ to be
unitary. For example, in the table for $E_8$, for the nilpotent
$A_4+A_1$, the $\nu$-string is $(\nu_1,\nu_2,\nu_3)$ and the
centralizer is $A_2+T_1.$ This means that the unitary parameters are
those for which $\nu_3=0$ (this is the $T_1$-piece), $\nu_2=0$ and
$0\le\nu_1<\frac 12$ (this is the $0$-complementary series of $A_2$).

There is one difference in $E_6$ due to the fact that we only consider
hermitian spherical infinitesimal characters $\chi.$ In this table,
the $\nu$-string already refers to the semisimple and hermitian
spherical parameter of the centralizer. For example, the nilpotent
$A_2+A_1$ in $E_6$ has centralizer $A_2+T_1,$ and the corresponding
$\chi$ has a single $\nu$. This $\nu$ corresponds to the hermitian
parameter in the $A_2$ part of $\fz(\CO)$, so it must satisfy
$0\le\nu<\frac 12.$

\begin{tiny}
\begin{center}
\begin{longtable}{|c|c|c|}
\caption{Table of hermitian parameters $(\CO,\nu)$ for $E_6$}\label{table:E6}\\
\hline
\multicolumn{1}{|c|}{$\mathbf\CO$} &
\multicolumn{1}{c|}{$\mathbf\chi$} 
& \multicolumn{1}{c|}{$\mathbf{\fz(\CO)}$} \\ \hline 
\endfirsthead

\multicolumn{3}{c}%
{{  \tablename\ \thetable{} -- continued from previous page}} \\
\hline
\multicolumn{1}{|c|}{$\mathbf\CO$} &
\multicolumn{1}{c|}{$\mathbf\chi$} 
& \multicolumn{1}{c|}{$\mathbf{\fz(\CO)}$}
 \\ \hline 
\endhead


\hline \hline
\endlastfoot


$E_6$ &$(0,1,2,3,4,-4,-4,4)$  &$1$\\
\hline
$E_6(a_1)$ &$(0,1,1,2,3,-3,-3,3)$  &$1$\\
\hline
$D_5$ &$(\frac 12,\frac 12,\frac 32,\frac 32,\frac 52,-\frac 52,-\frac
  52,\frac 52)$  &$T_1$\\
\hline
$E_6(a_3)$ &$(0,0,1,1,2,-2,-2,2)$  &$1$\\
\hline
$D_5(a_1)$ &$(\frac 14,\frac 34,\frac 34,\frac 54,-\frac 74,-\frac
  74,\frac 74)$ 
  &$T_1$\\
\hline
$A_5$ &$(-\frac {11}4,-\frac 74,-\frac 34,\frac 14,\frac 54,-\frac
  54,-\frac 54,\frac 54)+\nu(\frac 12,\frac 12,\frac 12,\frac 12,\frac 12,-\frac 12,-\frac
  12,\frac 12)$  &$A_1$\\

\hline
$A_4+A_1$ &$(0,\frac 12,\frac 12,1,\frac 32,-\frac 32,-\frac 32,\frac
  32)$  &$T_1$\\
\hline
$D_4$ &$(0,1,2,3,\nu,-\nu,-\nu,\nu)$  &$A_2$\\
\hline
$A_4$ &$(-2,-1,0,1,2,0,0,0)+\nu(\frac 12,\frac 12,\frac 12,\frac 12,\frac 12,-\frac
  12,-\frac 12,\frac 12)$  &$A_1T_1$\\
\hline
$D_4(a_1)$ &$(0,0,1,1,1,-1,-1,1)$  &$T_2$\\
\hline
$A_3+A_1$ &$(-\frac 54,-\frac 14,\frac 34,-\frac 54,-\frac 14,-\frac
  34,-\frac 34,\frac 34)+\nu(\frac 12,\frac 12,\frac 12,\frac 12,\frac 12,-\frac 12,-\frac
  12,\frac 12)$  &$A_1T_1$\\
\hline
$2A_2+A_1$ &$(0,1,-\frac 32,-\frac 12,\frac 12,-\frac 12,-\frac
  12,\frac 12)+\nu(0,0,1,1,1,-1,-1,1)$  &$A_1$\\
\hline
$A_3$ &$(-\frac 32,-\frac 12,\frac 12,\frac
  32,0,0,0,0)+(\frac {\nu_1}2,\frac {\nu_1}2,\frac {\nu_1}2,\frac {\nu_1}2,\frac
  {\nu_2}2,-\frac {\nu_2}2,-\frac {\nu_2}2,\frac {\nu_2}2)$  &$B_2T_1$\\

\hline
$A_2+2A_1$ &$(\frac 54,-\frac 14,\frac 34,-\frac 34,\frac 14,-\frac
  14,-\frac 14,\frac 14)+\nu(-\frac 12,\frac 12,\frac 12,\frac 32,\frac 32,-\frac 32,-\frac 32,\frac 32)$   &$A_1T_1$\\

\hline
$2A_2$ &$(-\frac 12,\frac 12,-\frac 32,-\frac 12,\frac 12,-\frac
  12,-\frac 12,\frac 12)+$  &$G_2$\\
&$(\frac {\nu_2}2,\frac {\nu_2}2,\frac{2\nu_1+\nu_2}2,\frac
  {2\nu_1+\nu_2}2,\frac{2\nu_1+\nu_2}2,-\frac
  {2\nu_1+\nu_2}2,-\frac{2\nu_1+\nu_2}2,\frac {2\nu_1+\nu_2}2)$ &\\
\hline

$A_2+A_1$ &$(-\frac 12,\frac 12,-1,0,-\frac 12,-\frac 12,-\frac
  12,\frac 12)+\nu(\frac 12,\frac 12,\frac 12,\frac 12,\frac 12,-\frac 12,-\frac
  12,\frac 12)$  &$A_2T_1$\\

\hline
$A_2$ &$(0,-1,0,1,0,0,0,0)+(\frac {-\nu_1+\nu_2}2,\frac {-\nu_1+\nu_2}2,$  &$2A_2$\\
      &$\frac
  {-\nu_1+\nu_2}2,\frac {-\nu_1+\nu_2}2,\frac {\nu_1+\nu_2}2,\frac
  {-\nu_1+\nu_2}2,\frac {-\nu_1+\nu_2}2,\frac {\nu_1+\nu_2}2)$ &\\

\hline
$3A_1$ &$(0,1,-\frac 12,\frac 12,0,0,0,0)+(0,0,\nu_1,\nu_2,\nu_1,-\nu_1,-\nu_1,\nu_1)$  &$A_2A_1$\\
     \hline
$2A_1$ &$(-\frac 12,\frac 12,-\frac 12,\frac 12,0,0,0,0)+$
   &$B_3T_1$\\
&$(\frac {-\nu_1+\nu_2}2,\frac {-\nu_1+\nu_2}2,\frac
  {\nu_1+\nu_2}2,\frac {\nu_1+\nu_2}2,\frac {\nu_1}2,-\frac
  {\nu_1}2,-\frac {\nu_1}2,\frac {\nu_1}2)$ &\\
 
\hline
$A_1$ &$(\frac 12,\frac 12,0,0,0,0,0,0)+(\frac {-\nu_1+\nu_2}2,\frac {\nu_1-\nu_2}2,$
  &$A_5$\\
      &$\frac
  {-\nu_1+\nu_2}2+\nu_3,\frac {\nu_1-\nu_2}2+\nu_3,\frac {\nu_1+\nu_2}2,-\frac
  {\nu_1+\nu_2}2,-\frac {\nu_1+\nu_2}2,\frac {\nu_1+\nu_2}2)$ &\\

\end{longtable}
\end{center}
\end{tiny}

\begin{tiny}
\begin{center}
\begin{longtable}{|c|c|c|}
\caption{Table of parameters $(\CO,\nu)$ for $E_7$}\label{table:E7}\\
\hline
\multicolumn{1}{|c|}{$\mathbf\CO$} &
\multicolumn{1}{c|}{$\mathbf\chi$} 
& \multicolumn{1}{c|}{$\mathbf{\fz(\CO)}$} \\ \hline 
\endfirsthead

\multicolumn{3}{c}%
{{ \tablename\ \thetable{} -- continued from previous page}} \\
\hline
\multicolumn{1}{|c|}{$\mathbf\CO$} &
\multicolumn{1}{c|}{$\mathbf\chi$} 
& \multicolumn{1}{c|}{$\mathbf{\fz(\CO)}$}
 \\ \hline 
\endhead


\hline \hline
\endlastfoot


$E_7$ &$(0,1,2,3,4,5,-\frac {17}2,\frac {17}2)$   &$1$\\
\hline
$E_7(a_1)$  &$(0,1,1,2,3,4,-\frac{13}2,\frac {13}2)$   &$1$\\
\hline
$E_7(a_2)$  &$(0,1,1,2,2,3,-\frac {11}2,\frac {11}2)$   &$1$\\
\hline
$E_7(a_3)$  &$(0,0,1,1,2,3,-\frac 92,\frac 92)$   &$1$\\
\hline
$E_6$       &$(0,1,2,3,4,-4,-4,4)+\nu(0,0,0,0,0,1,-\frac 12,\frac 12)$
  &$A_1^\ell$\\
\hline
$D_6$  &$(0,1,2,3,4,5,0,0)+\nu(0,0,0,0,0,0,-1,1)$     &$A_1$\\
\hline
$E_6(a_1)$ &$(0,1,1,2,3,-3,-3,3)+\nu(0,0,0,0,0,1,-\frac 12,\frac 12)$
 &$T_1$\\
\hline
$E_7(a_4)$ &$(0,0,1,1,1,2,-\frac 72,\frac 72)$  &$1$\\
\hline
$D_6(a_1)$ &$(0,1,1,2,3,4,0,0)+\nu(0,0,0,0,0,0,-1,1)$  &$A_1$\\
\hline
$A_6$ &$(-\frac 72,-\frac 52,-\frac 32,-\frac 12,\frac 12,\frac
32,-\frac 32,\frac 32)+\nu(\frac 12,\frac 12,\frac 12,\frac 12,\frac
12,\frac 12,-1,1)$  &$A_1^\ell$\\
\hline
$D_5+A_1$ &$(0,1,2,3,-\frac 52,-\frac 32,-2,2)+\nu(0,0,0,0,1,1,-1,1)$
 &$A_1$\\
\hline
$E_7(a_5)$ &$(0,0,1,1,1,2,-\frac 52,\frac 52)$   &$1$\\
\hline
$D_6(a_2)$ &$(0,1,1,2,2,3,0,0)+\nu(0,0,0,0,0,0,-1,1)$  &$A_1$\\
\hline
$A_5+A_1$ &$(\frac {11}4,-\frac 74,-\frac 34,\frac 14,\frac 54,\frac
94,-\frac 14,\frac 14)+\nu(-\frac 12,\frac 12,\frac 12,\frac 12,\frac
12,\frac 12,-\frac 32,\frac 32)$  &$A_1$\\
\hline
$D_5$ &$(0,1,2,3,-2,-2,-2,2)+\nu_1(0,0,0,0,1,1,-1,1)$ 
&$2A_1$\\
      &$+\nu_2(0,0,0,0,-1,1,0,0)$ &\\
\hline
$E_6(a_3)$ &$(0,0,1,1,2,-2,-2,2)+\nu(0,0,0,0,0,1,-\frac 12,\frac 12)$
 &$A_1^\ell$\\
\hline
$D_5(a_1)A_1$ &$(0,1,1,2,-2,-1,-\frac 32,\frac
32)+\nu(0,0,0,0,1,1,-1,1)$  &$A_1^\ell$\\
\hline
$(A_5)'$ &$(-\frac 52,-\frac 32,-\frac 12,\frac 12,\frac 32,\frac
52,0,0)+\nu_1(0,0,0,0,0,0,-1,1)$  &$A_1A_1^\ell$\\
&$+\nu_2(\frac 12,\frac 12,\frac 12,\frac 12,\frac 12,\frac 12,0,0)$ 
&\\
\hline
$A_4+A_2$ &$(0,1,2,-2,-1,0,-1,1)+\nu(0,0,0,1,1,1,-\frac 32,\frac 32)$
 &$A_1^\ell$\\
\hline
$(A_5)''$ &$(\frac 52,-\frac 32,-\frac 12,\frac 12,\frac 32,\frac
52,0,0)+\nu_2(0,0,0,0,0,0,-1,1)$  &$G_2$\\
&$+\nu_1(-\frac 12,\frac 12,\frac 12,\frac 12,\frac 12,\frac 12,-\frac
32,\frac 32)$ &\\
\hline
$D_5(a_1)$ &$(0,1,1,2,3,0,0,0)+\nu_1(0,0,0,0,0,0,-1,1)$
 &$A_1T_1$\\
&$+\nu_2(0,0,0,0,0,1,-\frac 12,\frac 12)$ &\\
\hline
$A_4+A_1$ &$(\frac 94,-\frac 54,-\frac 14,\frac 34,\frac 74,-\frac
14,-\frac 12,\frac 14)+\nu_1(\frac 12,\frac 12,\frac 12,\frac 12,\frac
12,\frac 12,-1,1)$  &$T_2$\\
&$+\nu_2(0,0,0,0,0,1,-\frac 12,\frac 12)$  &\\
\hline
$D_4+A_1$ &$(0,1,2,3,-\frac 12,\frac 12,0,0)+\nu_1(0,0,0,0,-\frac
12,-\frac 12,-\frac 12,\frac 12)$  &$B_2$\\
&$+\nu_2(0,0,0,0,\frac 12,\frac 12,-\frac 12,\frac 12)$  &\\
\hline
$A_3A_2A_1$ &$(0,1,-2,-1,0,1,-\frac 12,\frac
12)+\nu(0,0,1,1,1,1,-2,2)$  &$A_1$\\
\hline
$A_4$ &$(0,-2,-1,0,1,2,0,0)+\nu_1(0,0,0,0,0,0,-1,1)$
 &$A_2T_1$ \\
&$+\nu_2(-\frac 12,\frac 12,\frac 12,\frac 12,\frac 12,\frac 12,-\frac
32,\frac 32)+\nu_3(\frac 12,\frac 12,\frac 12,\frac 12,\frac 12,\frac
12,-1,1)$ &\\
\hline
$A_3+A_2$ &$(0,1,2,-1,0,1,0,0)+\nu_1(0,0,0,0,0,0,-1,1)$
 &$A_1T_1$\\
&$+\nu_2(0,0,0,1,1,1,0,0)$ &\\
\hline
$D_4$ &$(0,1,2,3,\nu_2-\nu_1,\nu_2+\nu_1,-\nu_3,\nu_3)$
 &$C_3$\\
\hline
$D_4(a_1)A_1$ &$(0,1,1,2,-\frac 12,\frac
12,0,0)+(0,0,0,0,\nu_2,\nu_2,-\nu_1,\nu_1)$  &$2A_1$\\
\hline
$A_3+2A_1$ &$(0,1,-\frac 32,-\frac 12,\frac 12,\frac
32,0,0)+(0,0,\nu_2,\nu_2,\nu_2,\nu_2,-\nu_1,\nu_1)$ 
&$2A_1$\\
\hline
$D_4(a_1)$ &$(0,1,1,2,\nu_2-\nu_3,\nu_2+\nu_3,-\nu_1,\nu_1)$
 &$3A_1$\\
\hline
$(A_3+A_1)'$ &$(0,1,2,0,-\frac 12,\frac
12,0,0)+(0,0,0,2\nu_2,\nu_3,\nu_3,-\nu_1,\nu_1)$
 &$3A_1$\\
\hline
$2A_2+A_1$ &$(\frac 54,-\frac 14,\frac 34,-\frac 54,-\frac 14,\frac
34,-\frac 14,\frac 14)+\nu_1(1,-1,-1,1,1,1,0,0)$ 
&$2A_1$\\
&$+\nu_2(-\frac 12,\frac 12,\frac 12,\frac 12,\frac 12,\frac 12,-\frac
32,\frac 32)$ &\\
\hline
$(A_3+A_1)''$ &$(\frac 32,-\frac 12,\frac 12,\frac 32,-\frac 12,\frac
12,0,0)+$  &$B_3$\\
&$(-\frac {\nu_1}2,\frac {\nu_1}2,\frac {\nu_1}2,\frac {\nu_1}2,\frac
{\nu_3-\nu_2}2,\frac {\nu_3-\nu_2}2,-\frac {\nu_3+\nu_2}2,\frac
{\nu_3+\nu_2}2)$  &\\
\hline
$A_2+3A_1$ &$(0,1,-1,0,-1,0,-\frac 12,\frac
12)+\nu_1(0,0,1,1,1,1,-2,2)$  &$G_2$\\
&$+\nu_2(0,0,0,0,1,1,-1,1)$  &\\
\hline


$2A_2$ &$(-\frac 12,\frac 12,-\frac 32,-\frac 12,\frac 12,-\frac
12,-\frac 12,\frac 12)+\nu_1(0,0,1,1,1,-1,-1,1)$
 &$G_2A_1$\\
&$+\nu_2(\frac 12,\frac 12,\frac 12,\frac 12,\frac 12,-\frac 12,-\frac
12,\frac 12)+\nu_3(0,0,0,0,0,1,-\frac 12,\frac 12)$ & \\
\hline
$A_3$ &$(0,1,2,\nu_1,\nu_2,\nu_3,\nu_4)$ 
&$B_3A_1$\\
\hline
$*A_2+2A_1$
&$(0,1,-1,0,1,0,0,0)+(0,0,\nu_2,\nu_2,\nu_2,\nu_3,-\nu_1,\nu_1)$
 &$A_12A_1^\ell$\\
\hline
$A_2+A_1$ &$(1,0,1,0,-\frac 12,\frac 12,0,0)+(0,0,0,0,\nu_2,\nu_2,-\nu_1,\nu_1)$
 &$A_3T_1$\\
&$+\nu_3(0,0,0,1,1,1,-\frac 32,\frac 32)+\nu_4(-\frac 12,\frac
12,\frac 12,\frac 12,\frac 12,\frac 12,-\frac 32,\frac 32)$ &\\
\hline
$4A_1$ &$(0,1,-\frac 12,\frac 12,-\frac 12,\frac
12,0,0)+(0,0,\nu_3,\nu_3,\nu_2,\nu_2,-\nu_1,\nu_1)$
 &$C_3$ \\
\hline
$A_2$
&$(1,0,1,0,0,0,0,0)+(0,0,0,0,\nu_2-\nu_3,\nu_2+\nu_3,-\nu_1,\nu_1)$
 &$A_5$\\
&$+\nu_4(-\frac 12,\frac 12,\frac 12,\frac 12,\frac 12,\frac 12,-\frac
32,\frac 32)+\nu_5(0,0,0,1,1,1,-\frac 32,\frac 32)$ &\\
\hline
$(3A_1)'$ &$(-\frac 12,\frac 12,-\frac 12,\frac 12,-\frac 12,\frac
12,0,0)+(\nu_1,\nu_1,\nu_2,\nu_2,\nu_3,\nu_3,-\nu_4,\nu_4)$
 &$C_3A_1$\\
\hline
$(3A_1)''$ &$(\frac 12,\frac 12,-\frac 12,\frac 12,-\frac 12,\frac
12,0,0)+(-\nu_4,\nu_4,\nu_3,\nu_3,\nu_2,\nu_2,-\nu_1,\nu_1)$
 &$F_4$\\
\hline
$2A_1$ &$(0,1,\nu_1,\nu_2,\nu_3,\nu_4,-\nu_5,\nu_5)$
 &$B_4A_1$\\
\hline
$A_1$ &$(\frac {\nu_1+\nu_2+\nu_3-\nu_4}2,\frac
    {\nu_1+\nu_2-\nu_3+\nu_4}2,\frac {\nu_1-\nu_2+\nu_3+\nu_4}2,\frac
    {-\nu_1+\nu_2+\nu_3+\nu_4}2,$  &$D_6$\\
&$-\frac 12+\frac {-\nu_5+\nu_6}2,\frac 12+\frac
    {-\nu_5+\nu_6}2,-\frac {\nu_5+\nu_6}2,\frac {\nu_5+\nu_6}2)$ &\\
\end{longtable}
\end{center}
\end{tiny}

\noindent{\bf $\mathbf{E_7}$ exception:} 

\noindent $\mathbf{A_2+2A_1}$. Three regions: $\{0\le\nu_1<\frac
12,0\le\nu_2<1,0\le\nu_3<1, \nu_1+\frac {3\nu_2}2+\frac{\nu_3}2<\frac
32\}$, $\{0\le\nu_1<\frac 12,0\le\nu_2<1,0\le\nu_3<1, -\nu_1+\frac
{3\nu_2}2+\frac{\nu_3}2<\frac 32,\nu_1+\frac
{3\nu_2}2-\frac{\nu_3}2>\frac 32 \}$, and  $\{0\le\nu_1<\frac
12,0\le\nu_2<1,0\le\nu_3<1, \frac{3\nu_2}2+\frac{\nu_3}2>\frac 32,\nu_1+\frac
{3\nu_2}2-\frac{\nu_3}2<\frac 32 \}$.

\begin{tiny}
\begin{center}
\begin{longtable}{|c|c|c|}
\caption{Table of parameters $(\CO,\nu)$ for $E_8$}\label{table:E8}\\
\hline
\multicolumn{1}{|c|}{$\mathbf\CO$} &
\multicolumn{1}{c|}{$\mathbf\chi$} 
& \multicolumn{1}{c|}{$\mathbf{\fz(\CO)}$} \\ \hline 
\endfirsthead

\multicolumn{3}{c}%
{{  \tablename\ \thetable{} -- continued from previous page}} \\
\hline
\multicolumn{1}{|c|}{$\mathbf\CO$} &
\multicolumn{1}{c|}{$\mathbf\chi$} 
& \multicolumn{1}{c|}{$\mathbf{\fz(\CO)}$}
 \\ \hline 
\endhead


\hline \hline
\endlastfoot

$E_8$ &$(0,1,2,3,4,5,6,23)$   &$1$\\
\hline
$E_8(a_1)$ &$(0,1,1,2,3,4,5,18)$ &$1$\\
\hline
$E_8(a_2)$ &$(0,1,1,2,2,3,4,15)$ &$1$\\
\hline
$E_8(a_3)$ &$(0,0,1,1,2,3,4,13)$ &$1$\\
\hline
$E_8(a_4)$ &$(0,0,1,1,2,2,3,11)$ &$1$\\
\hline
$E_7$   &$(0,1,2,3,4,5,-\frac{17}2,\frac {17}2)+\nu(0,0,0,0,0,0,1,1)$
&$A_1$\\
\hline
$E_8(b_4)$ &$(0,0,1,1,1,2,3,10)$ &$1$\\
\hline
$E_8(a_5)$ &$(0,0,1,1,1,2,2,9)$ &$1$\\
\hline
$E_7(a_1)$
&$(0,1,1,2,3,4,-\frac{13}2,\frac{13}2)+\nu(0,0,0,0,0,0,1,1)$ &$A_1$\\
\hline
$E_8(b_5)$ &$(0,0,1,1,1,2,3,8)$ &$1$\\
\hline
$D_7$ &$(0,1,2,3,4,5,6,0)+\nu(0,0,0,0,0,0,0,2)$ &$A_1$\\
\hline
$E_8(a_6)$ &$(0,0,1,1,1,2,2,7)$ &$1$\\
\hline
$E_7(a_2)$ &$(0,1,1,2,2,3,-\frac{11}2,\frac
{11}2)+\nu(0,0,0,0,0,0,1,1)$ &$A_1$\\
\hline
$E_6+A_1$ &$(0,1,2,3,4,-\frac 92,-\frac 72,4)+\nu(0,0,0,0,0,0,1,1,2)$
&$A_1$\\
\hline
$D_7(a_1)$ &$(0,1,1,2,3,4,5,0)+\nu(0,0,0,0,0,0,0,2)$ &$T_1$\\
\hline
$E_8(b_6)$ &$(0,0,1,1,1,1,2,6)$ &$1$\\
\hline
$E_7(a_3)$ &$(0,0,1,1,2,3,-\frac 92,\frac 92)+\nu(0,0,0,0,0,0,1,1)$
&$A_1$\\
\hline
$E_6(a_1)A_1$  &$(0,1,1,2,3,-\frac 72,-\frac
52,3)+\nu(0,0,0,0,0,1,1,2)$  &$T_1$\\
\hline
$A_7$ &$(-\frac{17}4,-\frac {13}4,-\frac 94,-\frac 54,-\frac 14,\frac
34,\frac 74,\frac 74)+\nu(\frac 12,\frac 12,\frac 12,\frac 12,\frac
12,\frac 12,\frac 12,\frac 52)$ &$T_1$\\
\hline
$E_6$ &$(0,1,2,3,4,-4,-4,4)+\nu_1(0,0,0,0,0,1,1,2)$ &$G_2$\\
      &$+\nu_2(0,0,0,0,0,0,1,1)$ &\\
\hline
$D_6$ &$(0,1,2,3,4,5,\nu_1,\nu_2)$ &$B_2$\\
\hline
$D_5+A_2$ &$(0,1,2,3,-3,-2,-1,2)+\nu(0,0,0,0,1,1,1,3)$ &$T_1$\\
\hline
$E_6(a_1)$ &$(0,1,1,2,3,-3,-3,3)+\nu_2(0,0,0,0,0,1,1,2)$ &$A_2$\\
           &$+\nu_1(0,0,0,0,0,0,1,1)$ &\\
\hline
$E_7(a_4)$ &$(0,0,1,1,1,2,-\frac 72,\frac 72)+\nu(0,0,0,0,0,0,1,1)$
&$A_1$\\
\hline
$A_6+A_1$ &$(\frac {13}4,-\frac 94,-\frac 54,-\frac 14,\frac34,\frac
74,\frac {11}4,\frac 14)+\nu(-\frac 12,\frac 12,\frac 12,\frac
12,\frac 12,\frac 12,\frac 12,\frac 72)$ &$A_1$\\
\hline
$D_6(a_1)$ &$(0,1,1,2,3,4,0,0)+\nu_1(0,0,0,0,0,0,-1,1)$
&$2A_1$\\
           &$+(0,0,0,0,0,0,1,1)$ &\\
\hline
$A_6$ &$(-3,-2,-1,0,1,2,3,0)+\nu_2(\frac 12,\frac 12,\frac 12,\frac 12,\frac 12,\frac
12,\frac 12,\frac 12)$ &$2A_1$\\
&$+\nu_1(-\frac 12,-\frac 12,-\frac 12,-\frac
12,-\frac 12,-\frac 12,\frac 72)$ &\\
\hline
$E_8(a_7)$ &$(0,0,0,1,1,1,1,4)$ &$1$\\
\hline
$D_5+A_1$ &$(0,1,2,3,4,-\frac 12,\frac 12,0)+\nu_1(0,0,0,0,0,0,0,2)$
&$2A_1$\\
         &$+\nu_2(0,0,0,0,0,1,1,0)$ &\\
\hline
$E_7(a_5)$ &$(0,0,1,1,1,2,-\frac 52,\frac 52)+\nu(0,0,0,0,0,0,1,1)$
&$A_1$\\
\hline
$E_6(a_3)A_1$ &$(0,0,1,1,2,-\frac52,-\frac 32,2)+\nu(0,0,0,0,0,1,1,2)$
&$A_1$\\
\hline
$D_6(a_2)$ &$(0,1,1,2,2,3,-\nu_1+\nu_2,\nu_1+\nu_2)$ &$2A_1$\\
\hline
$D_5(a_1)A_2$ &$(0,1,1,2,-\frac 52,-\frac 32,-\frac 12,\frac
32)+\nu(0,0,0,0,1,1,1,3)$ &$A_1$\\
\hline
$A_5+A_1$ &$(\frac 14,-\frac{11}4,-\frac 74,-\frac 34,\frac 14,\frac
54,\frac 94,\frac 14)+\nu_2(-1,0,0,0,0,0,0,1)$ &$2A_1$\\
        &$+\nu_1(\frac 32,\frac 12,\frac 12,\frac 12,\frac 12,\frac
        12,\frac 12,\frac 32)$ &\\
\hline
$A_4+A_3$ &$(0,1,2,-\frac 52,-\frac 32,-\frac 12,\frac
12,1)+\nu(0,0,0,1,1,1,1,4)$ &$A_1$\\
\hline
$D_5$ &$(0,1,2,3,4,\nu_1,\nu_2,\nu_3)$ &$B_3$\\
\hline
$E_6(a_3)$ &$(0,0,1,1,2,-2,-2,2)+\nu_1(0,0,0,0,0,1,1,2)$
&$G_2$\\
&$+\nu_2(0,0,0,0,0,0,1,1)$ &\\
\hline
$D_4+A_2$ &$(0,1,2,3,-1,0,1,0)+\nu_2(0,0,0,0,1,1,1,3)$ &$A_2$\\
         &$+\nu_1(0,0,0,0,0,0,0,2)$ &\\
\hline
$*A_4A_2A_1$ &$(0,1,-\frac 52,-\frac 32,-\frac 12,\frac 12,\frac
32,\frac 12)+\nu(0,0,1,1,1,1,1,5)$ &$A_1$\\
\hline
$*D_5(a_1)A_1$ &$(0,1,1,2,3,-\frac 12+\nu_2,\frac 12+\nu_2,2\nu_1)$
&$A_1^\ell A_1$\\
\hline
$A_5$ &$(\frac 52,-\frac 32,-\frac 12,\frac 12,\frac 32,\frac
52,0,0)+\nu_1(-\frac 12,\frac 12,\frac 12,\frac 12,\frac 12,\frac
12,-\frac 32,\frac 32)$ &$G_2A_1$\\
&$+\nu_2(0,0,0,0,0,0,-1,1)+\nu_3(0,0,0,0,0,0,1,1)$ &\\
\hline
$*A_4+A_2$ &$(-\frac 12,\frac 12,-\frac 52,-\frac 32,-\frac 12,\frac
12,\frac 32,\frac 12)+\nu_2(1,1,0,0,0,0,0,0)$ &$2A_1$\\
&$+\nu_1(0,0,1,1,1,1,1,5)$ &\\
\hline
$A_4+2A_1$ &$(0,1,-2,-1,0,1,2,0)+\nu_1(0,0,0,0,0,0,0,2)$
&$A_1T_1$\\
&$\nu_2(0,0,1,1,1,1,1,0)$ &\\
\hline
$D_5(a_1)$ &$(0,1,1,2,3,\nu_3,\nu_2,\nu_1)$ &$A_3$\\
          &                                &\\
\hline
$2A_3$ &$(0,1,2,-\frac 32,-\frac 12,\frac 12,\frac
32,0)+\nu_2(0,0,0,\frac 12,\frac 12,\frac 12,\frac 12,1)$
&$B_2$\\
\hline
$A_4+A_1$ &$(0,1,2,-\frac 32,-\frac
12,-1,-1,1)+\nu_2(0,0,0,0,0,1,1,2)$ &$A_2T_1$\\
& $+\nu_1(0,0,0,0,0,0,1,1)+\nu_3(0,0,0,1,1,1,1,4)$
 &\\
\hline
$D_4(a_1)A_2$ &$(0,1,1,2,-1,0,1,0)+\nu_1(0,0,0,0,1,1,1,3)$
&$A_2$\\
             &$+\nu_2(0,0,0,0,0,0,0,2)$ &\\
\hline
$D_4+A_1$ &$(0,1,2,3,-\frac 12,\frac 12,0,0)+$ &$C_3$\\
&$(0,0,0,0,\nu_1,\nu_1,-\nu_2+\nu_3,\nu_2+\nu_3)$ &\\
\hline
$A_3A_2A_1$ &$(0,1,-2,-1,0,1,-\frac 12,\frac
12)+\nu_1(0,0,1,1,1,1,-2,2)$ &$2A_1$\\
&$+\nu_2(0,0,0,0,0,0,1,1)$ &\\
\hline
$A_4$ &$(0,-2,-1,0,1,2,0,0)+$ &$A_4$\\
&$(\nu_4,-\nu_1+\nu_2,\nu_3,\nu_3,\nu_3,\nu_3,\nu_3,\nu_1+\nu_2)$ &\\
\hline
$A_3+A_2$ &$(0,1,2,-1,0,1,0,0)+(0,0,0,\nu_3,\nu_3,\nu_3,\nu_1,\nu_2)$
&$B_2T_1$\\
\hline
$D_4(a_1)A_1$ &$(0,1,1,2,-\frac 12,\frac
12,0,0)+(0,0,0,0,\nu_1,\nu_1,-\nu_2+\nu_3,\nu_2+\nu_3)$
&$3A_1$\\
\hline
$A_3+2A_1$ &$(0,1,-\frac 32,-\frac 12,\frac 12,\frac
32,0,0)+(0,0,\nu_1,\nu_1,\nu_1,\nu_1,\nu_2,\nu_3)$
&$A_1B_2$\\
\hline
$2A_2+2A_1$ &$(0,1,-\frac 32,-\frac 12,\frac 12,-1,0,\frac
12)+\nu_1(0,0,-\frac 12,-\frac 12,-\frac 12,1,1,\frac 12)$ &$B_2$\\
&$+\nu_2(0,0,\frac 12,\frac 12,\frac 12,0,0,\frac 32)$ &\\
\hline
$D_4$
&$(0,1,2,3,\nu_3-\nu_4,\nu_3+\nu_4,\nu_1-\nu_2,\nu_1+\nu_2)$
&$F_4$\\
\hline
$D_4(a_1)$ &$(0,1,1,2,\nu_4,\nu_3,\nu_2,\nu_1)$  &$D_4$\\
\hline
$A_3+A_1$ &$(0,1,2,-\frac 12,\frac
12,0,0,0)+(0,0,0,\nu_1,\nu_1,\nu_2,\nu_3,\nu_4)$
&$A_1B_3$\\
\hline
$2A_2+A_1$ &$(0,1,-\frac 32,-\frac 12,\frac 12,-\frac 12,-\frac
12,\frac 12)+\nu_1(0,0,1,1,1,-1,-1,1)$ &$A_1G_2$\\
&$+\nu_2(0,0,0,0,0,1,1,2)+\nu_3(0,0,0,0,0,0,1,1)$ &\\
\hline
$2A_2$ &$(-\frac 12,\frac 12,-\frac 32,-\frac 12,\frac 12,-\frac
12,-\frac 12,\frac 12)+\nu_1(0,0,1,1,1,-1,-1,1)$
&$2G_2$\\
&$+\nu_2(\frac 12,\frac 12,\frac 12,\frac 12,\frac 12,-\frac 12,-\frac
12,\frac 12)+\nu_3(0,0,0,0,0,1,1,2)$ &\\
&$+\nu_4(0,0,0,0,0,0,1,1)$ &\\
\hline
$*A_2+3A_1$ &$(0,1,-1,0,-1,0,-\frac 12,\frac
12)+\nu_1(0,0,1,1,1,1,-2,2)$ &$G_2A_1$\\
&$+\nu_2(0,0,0,0,1,1,-1,1)+\nu_3(0,0,0,0,0,0,1,1)$ &\\
\hline
$A_3$ &$(0,1,2,\nu_1,\nu_2,\nu_3,\nu_4,\nu_5)$ &$B_5$\\
\hline
$*A_2+2A_1$
&$(0,1,-1,0,1,0,0,0)+(0,0,\nu_1,\nu_1,\nu_1,\nu_2,\nu_3,\nu_4)$
 &$A_1B_3$\\
\hline
$A_2+A_1$ &$(1,0,1,0,-\frac 12,\frac 12,0,0)+$ & $A_5$\\
&$(-\nu_5,\nu_5,\nu_5,\nu_4,\nu_3,\nu_3,-\nu_2+\nu_1,\nu_2+\nu_1)$
&\\
\hline
$*4A_1$ &$(0,1,-\frac 12,\frac 12,-\frac 12,\frac 12,0,0)+$
&$C_4$\\
&$(0,0,\nu_1,\nu_1,\nu_2,\nu_2,-\nu_3+\nu_4,\nu_3+\nu_4)$&\\
\hline
$A_2$ &$(\frac
{\nu_1-\nu_2-\nu_3+\nu_4}2,\frac{-\nu_1+\nu_2-\nu_3+\nu_4}2,\frac{-\nu_1-\nu_2+\nu_3+\nu_4}2,$
&$E_6$\\
&$\frac{\nu_1+\nu_2+\nu_3+\nu_4}2,-1+\frac
{\nu_5-\nu_6}2,\frac{\nu_5-\nu_6}2,1+\frac{\nu_5-\nu_6}2,\frac{\nu_5+3\nu_6}2)$
&\\
\hline
$3A_1$ &$(\frac 12,\frac 12,-\frac 12,\frac 12,-\frac 12,\frac
12,0,0)+(-\nu_4,\nu_4,\nu_3,\nu_3,\nu_2,\nu_2,-\nu_1,\nu_1)$
&$F_4A_1$\\
&$+\nu_5(0,0,0,0,0,0,1,1)$ &\\
\hline
$2A_1$ &$(0,1,\nu_1,\nu_2,\nu_3,\nu_4,\nu_5,\nu_6)$
&$B_6$\\
\hline
$A_1$
&$(\frac{\nu_1+\nu_2+\nu_3-\nu_4}2,\frac{\nu_1+\nu_2-\nu_3+\nu_4}2,\frac
{\nu_1-\nu_2+\nu_3+\nu_4}2,\frac {-\nu_1+\nu_2+\nu_3+\nu_4}2,$ &$E_7$\\
&$\frac
{-\nu_5-\nu_6+2\nu_7}2,-\frac 12+\frac {-\nu_5+\nu_6}2,\frac
12+\frac{-\nu_5+\nu_6}2, \frac {\nu_5+\nu_6+2\nu_7}2)$ &\\

\end{longtable}
\end{center}
\end{tiny}

\noindent{\bf $\mathbf {E_8}$ exceptions:}

\noindent$\mathbf{A_4+A_2+A_1}.$ $\{0\le\nu<\frac 3{10}\}.$

\noindent$\mathbf{D_5(a_1)+A_1}.$ Two regions: $\{0\le\nu_2<\frac 12,2\nu_1+\nu_2<\frac
32\},$ $\{0\le\nu_1<1,2\nu_1-\nu_2>\frac 32\}.$

\noindent$\mathbf{A_4+A_2}.$ Two regions: $\{0\le\nu_2<\frac
12,5\nu_1+\nu_2<2\},$ and $\{0\le\nu_1<\frac 12,5\nu_1-\nu_2>2\}.$

\noindent$\mathbf{A_2+3 A_1}$. Four regions: $\{3\nu_1+2\nu_2<1,
0\le\nu_3<\frac 12\}$, $\{2\nu_1+\nu_2<1<3\nu_1+\nu_2,0\le\nu_3<\frac
12,3\nu_1+2\nu_2+\nu_3<\frac 32\}$,
$\{2\nu_1+\nu_2<1<3\nu_1+\nu_2,0\le\nu_3<\frac
12,3\nu_1+\nu_2+\nu_3<\frac 32<3\nu_1+2\nu_2-\nu_3\}$, and
$\{2\nu_1+\nu_2<1<3\nu_1+\nu_2,0\le\nu_3<\frac 12,
3\nu_1+2\nu_2-\nu_3<\frac 32<3\nu_1+\nu_2+\nu_3\}.$

\noindent$\mathbf{A_2+2A_1}.$ Seven regions:
$\{0\le\nu_1<1,\nu_3+\nu_4<1, 3\nu_1+\nu_2+\nu_3+\nu_4<3\}$, 
$\{0\le\nu_1<1,\nu_3+\nu_4<1,3\nu_1+\nu_2-\nu_3+\nu_4<3<3\nu_1-\nu_2+\nu_3+\nu_4\}$,
$\{0\le\nu_1<1,\nu_3+\nu_4<1, 3\nu_1-\nu_2-\nu_3+\nu_4>3\}$,
$\{0\le\nu_1<1,\nu_3+\nu_4<1, 3\nu_1+\nu_2+\nu_3-\nu_4>3 \}$,
$\{0\le\nu_1<1,\nu_2+\nu_4>1,\nu_2+\nu_3<1,\nu_4<1,3\nu_1+\nu_2+\nu_3+\nu_4<3 \}$,
$\{0\le\nu_1<1,\nu_2+\nu_4>1,\nu_2+\nu_3<1,\nu_4<1,
3\nu_1-\nu_2-\nu_3+\nu_4>3 \}$, and
$\{0\le\nu_1<1,\nu_2+\nu_4>1,\nu_2+\nu_3<1,\nu_4<1,
3\nu_1+\nu_2+\nu_3-\nu_4>3\}$.

\noindent$\mathbf{4A_1}$. Two regions:
$\{0\le\nu_1\le\nu-2\le\nu_3\le\nu_4<\frac 12\}$ and
$\{\nu_1+\nu_4<1,\nu_2+\nu_3<1,\nu_2+\nu_4>1, -\nu_1+\nu_3+\nu_4<\frac
32<\nu_1+\nu_3+\nu_4\}.$

\subsection{$0$-complementary series}

 We record next the precise description of the $0$-complementary series
(that is, the generic spherical unitary
parameters) for types $E_6,\ E_7,\ E_8$. This answer is obtained inductively
from corollary \ref{c:7.4}.

\subsubsection{$E_6$}\label{sec:genspherE6} In $W(E_6)$, the longest Weyl
group element $w_0$ does not act by $-1.$ Therefore, we only consider
dominant parameters $\chi$ such that $w_0\chi=-\chi:$
\begin{small}\begin{align}\notag
\left(\frac {\nu_1-\nu_2}2-\nu_3,\frac {\nu_1-\nu_2}2-\nu_4,\frac {\nu_1-\nu_2}2+\nu_4,\frac {\nu_1-\nu_2}2+\nu_3,\frac {\nu_1+\nu_2}2, -\frac {\nu_1+\nu_2}2,-\frac {\nu_1+\nu_2}2,\frac {\nu_1+\nu_2}2\right).
\end{align}\end{small}
The $0$-complementary series is:

\begin{itemize}
\item[(1)] $\al_{36}<1,$ and $\al_1,\al_2,\al_3,\al_4,\al_5,\al_6\ge 0.$
\item[(2)] $\al_{34}<1,$ $\al_{35}>1,$ and $\al_1,\al_2,\al_3,\al_5,\al_6\ge 0.$
\end{itemize}

\subsubsection{$E_7$}\label{sec:genspherE7}

The  parameters are 
$\chi=(\nu_1,\nu_2,\nu_3,\nu_4,\nu_5,\nu_6,-\nu_7,\nu_7),$ assumed
dominant.
The $0$-complementary series is:

\smallskip

\begin{enumerate}

\item $\alpha_{63}<1,$ and $\al_1,\al_2,\al_3,\al_4,\al_5,\al_6,\al_7\ge 0.$

\item $\alpha_{61}<1,\ \ \alpha_{62}>1$ and
  $\al_1,\al_2,\al_4,\al_5,\al_6,\al_7\ge 0.$

\item $\alpha_{58}<1,\ \alpha_{59}<1,\ \ \alpha_{60}>1$ and
  $\al_1,\al_3,\al_4,\al_6,\al_7\ge 0.$

\item $\alpha_{53}<1,\ \alpha_{54}<1,\ \alpha_{55}<1,\ \
  \alpha_{56}>1,\ \alpha_{57}>1$ and $\al_1,\al_3,\al_5\ge 0.$

\item $\alpha_{46}<1,\ \alpha_{47}<1,\ \alpha_{48}<1,\
  \alpha_{49}<1,\ \  \alpha_{50}>1,\   \alpha_{51}>1,\  \alpha_{52}>1$
  and $\al_2\ge 0.$

\item $\alpha_{53}<1,\ \alpha_{59}<1,\ \ \alpha_{56}>1$ and
  $\al_1,\al_3,\al_4,\al_5,\al_6\ge 0.$

\item $\alpha_{49}<1,\ \alpha_{53}<1,\ \alpha_{54}<1,\ \
  \alpha_{52}>1,\ \alpha_{56}>1$ and $\al_3,\al_4,\al_5\ge 0.$

\item $\alpha_{47}<1,\ \alpha_{48}<1,\ \alpha_{49}<1,\ \alpha_{53}<1,\
  \ \alpha_{51}>1,\ \alpha_{52}>1$ and $\al_2,\al_4\ge 0.$

\end{enumerate}

\subsubsection{$E_8$}\label{sec:genspherE8}

The parameters are
$\chi=(\nu_1,\nu_2,\nu_3,\nu_4,\nu_5,\nu_6,\nu_7,\nu_8),$  assumed
dominant.
The $0$-complementary series is:

\smallskip

\begin{enumerate}

\item $\alpha_{120}<1$ and $\alpha_1,\alpha_2,\alpha_3,\alpha_4,\alpha_5,\alpha_6,\alpha_7,\alpha_8\ge0$.

\item $\alpha_{113}<1,\alpha_{114}<1$; $\alpha_{115}>1$ and $\alpha_1,\alpha_4,\alpha_5,\alpha_6,\alpha_7,\alpha_8\ge 0$.

\item $\alpha_{109}<1,\alpha_{110}<1$; $\alpha_{111}>1,\alpha_{112}>1$ and $\alpha_3,\alpha_5,\alpha_6,\alpha_7,\alpha_8\ge 0$.

\item $\alpha_{91}<1,\alpha_{92}<1,\alpha_{97}<1,\alpha_{98}<1$; $\alpha_{95}>1,\alpha_{96}>1,\alpha_{101}>1$ and $\alpha_3,\alpha_4\ge 0$.

\item $\alpha_{90}<1,\alpha_{91}<1,\alpha_{92}<1,\alpha_{97}<1$; $\alpha_{94}>1,\alpha_{95}>1,\alpha_{96}>1$  and $\alpha_1,\alpha_3\ge 0$.

\item $\alpha_{89}<1,\alpha_{90}<1,\alpha_{91}<1,\alpha_{92}<1$; $\alpha_{93}>1,\alpha_{94}>1,\alpha_{95}>1,$ $\alpha_{96}>1$ and $\alpha_1\ge 0$.

\item $\alpha_{104}<1,\alpha_{110}<1$; $\alpha_{107}>1,\alpha_{112}>1$ and $\alpha_3,\alpha_4,\alpha_5,\alpha_7,\alpha_8\ge 0$.

\item $\alpha_{104}<1,\alpha_{105}<1,\alpha_{106}<1$; $\alpha_{107}>1,\alpha_{108}>1$ and $\alpha_2,\alpha_4,\alpha_7,\alpha_8\ge 0$.

\item $\alpha_{118}<1$; $\alpha_{119}>1$ and $\alpha_1,\alpha_2,\alpha_3,\alpha_4,\alpha_5,\alpha_6,\alpha_8\ge 0$.

\item $\alpha_{97}<1,\alpha_{110}<1$; $\alpha_{101}>1,\alpha_{112}>1$ and $\alpha_3,\alpha_4,\alpha_5,\alpha_6,\alpha_7\ge 0$.

\item $\alpha_{97}<1,\alpha_{105}<1,\alpha_{106}<1$; $\alpha_{101}>1,\alpha_{108}>1$ and $\alpha_2,\alpha_4,\alpha_6,\alpha_7\ge 0$.

\item $\alpha_{116}<1$; $\alpha_{117}>1$ and $\alpha_1,\alpha_2,\alpha_3,\alpha_4,\alpha_6,\alpha_7,\alpha_8\ge 0$.

\item $\alpha_{97}<1,\alpha_{98}<1,\alpha_{106}<1$; $\alpha_{101}>1,\alpha_{102}>1$ and $\alpha_2,\alpha_4,\alpha_5,\alpha_6\ge 0$.

\item $\alpha_{97}<1,\alpha_{98}<1,\alpha_{99}<1$; $\alpha_{96}>1,\alpha_{101}>1,\alpha_{102}>1$ and $\alpha_2,\alpha_4,\alpha_5\ge 0$.

\item $\alpha_{97}<1,\alpha_{98}<1,\alpha_{99}<1,\alpha_{100}<1$; $\alpha_{101}>1,\alpha_{102}>1,\alpha_{103}>1$ and $\alpha_2,\alpha_5\ge 0$.

\item $\alpha_{114}<1$; $\alpha_{112}>1$ and $\alpha_1,\alpha_3,\alpha_4,\alpha_5,\alpha_6,\alpha_7,\alpha_8\ge 0$.

\end{enumerate}

\subsubsection{Roots for type $E$}\label{roots} The notation for the
positive roots which appeared in the lists of 
$0$-complementary series for $E_6,E_7,E_8$ is as follows.

\medskip

\begin{tiny}

\noindent\begin{longtable}{|lll|}
\hline
\begin{small}E6\end{small} &&\\
$\al_{34}=\frac 12(-1,1,-1,1,1,-1,-1,1)$ &$\al_{35}=\frac
  12(-1,-1,1,1,1,-1,-1,1)$
&$\al_{36}=\frac 12(1,1,1,1,1,-1,-1,1)$\\

\hline

\begin{small}E7\end{small} && \\

$\al_{46}=\frac 12(-1,1,-1,1,1,-1,-1,1)$ &$\al_{47}=\frac
  12(-1,1,1,-1,-1,1,-1,1)$
&$\al_{48}=\frac 12(1,-1,-1,1,-1,1,-1,1)$ \\
$\al_{49}=\ep_5+\ep_6$
& $\al_{50}=\frac 12(-1,-1,1,1,1,-1,-1,1)$ &$\al_{51}=\frac
  12(-1,1,-1,1,-1,1,-1,1)$\\
$\al_{52}=\frac 12(1,-1,-1,-1,1,1,-1,1)$ &$\al_{53}=\frac 12(1,1,1,1,1,-1,-1,1)$
&$\al_{54}=\frac 12(-1,-1,1,1,-1,1,-1,1)$\\
$\al_{55}=\frac
  12(-1,1,-1,-1,1,1,-1,1)$
& $\al_{56}=\frac 12(1,1,1,1,-1,1,-1,1)$ &$\al_{57}=\frac
  12(-1,-1,1,-1,1,1,-1,1)$\\

$\al_{58}=\frac 12(1,1,1,-1,1,1,-1,1)$ &$\al_{59}=\frac
  12(-1,-1,-1,1,1,1,-1,1)$
&$\al_{60}=\frac 12(1,1,-1,1,1,1,-1,1)$\\
$\al_{61}=\frac
  12(1,-1,1,1,1,1,-1,1)$
&$\al_{62}=\frac 12(-1,1,1,1,1,1,-1,1)$ &$\al_{63}=-\ep_7+\ep_8$\\




\hline
{\begin{small}E8\end{small}} &&\\ 

$\al_{89}=\frac 12(1,-1,1,1,1,1,-1,1)$ &$\al_{90}=\frac
  12(1,1,-1,1,1,-1,1,1)$
&$\al_{91}=\frac 12(1,1,1,-1,-1,1,1,1)$\\
$\al_{92}=\frac
  12(-1,-1,-1,1,-1,1,1,1)$ &$\al_{93}=\frac 12(-1,1,1,1,1,1,-1,1)$ &$\al_{94}=\frac
  12(1,-1,1,1,1,-1,1,1)$\\
$\al_{95}=\frac 12(1,1,-1,1,-1,1,1,1)$ &$\al_{96}=\frac
  12(-1,-1,-1,-1,1,1,1,1)$
&$\al_{97}=-\ep_7+\ep_8$\\
$\al_{98}=\frac 12(-1,1,1,1,1,-1,1,1)$
&$\al_{99}=\frac 12(1,-1,1,1,-1,1,1,1)$ &$\al_{100}=\frac
  12(1,1,-1,-1,1,1,1,1)$\\
$\al_{101}=-\ep_6+\ep_8$ &$\al_{102}=\frac 12(-1,1,1,1,-1,1,1,1)$
&$\al_{103}=\frac 12(1,-1,1,-1,1,1,1,1)$\\ 
$\al_{104}=-\ep_5+\ep_8$
&$\al_{105}=\frac 12(-1,1,1,-1,1,1,1,1)$ &$\al_{106}=\frac
  12(1,-1,-1,1,1,1,1,1)$\\

$\al_{107}=-\ep_4+\ep_8$ &$\al_{108}=\frac 12(-1,1,-1,1,1,1,1,1)$
&$\al_{109}=-\ep_3+\ep_8$\\
$\al_{110}=\frac 12(-1,-1,1,1,1,1,1,1)$
&$\al_{111}=-\ep_2+\ep_8$ &$\al_{112}=\frac 12(1,1,1,1,1,1,1,1)$\\

$\al_{113}=\ep_1+\ep_8$ &$\al_{114}=-\ep_1+\ep_8$
&$\al_{115}=\ep_2+\ep_8$\\
$\al_{116}=\ep_3+\ep_8$ &$\al_{117}=\ep_4+\ep_8$ &$\al_{118}=\ep_5+\ep_8$\\

$\al_{119}=\ep_6+\ep_8$ &$\al_{120}=\ep_7+\ep_8$ &\\

\hline
\end{longtable}

\end{tiny}

\ifx\undefined\bysame
\newcommand{\bysame}{\leavevmode\hbox to3em{\hrulefill}\,}
\fi

\end{document}